\documentclass[10pt]{amsart}
\usepackage{a4wide}
\usepackage{amsmath,amssymb,amsfonts,amsthm}
\usepackage{mathrsfs}
\usepackage{esint,url}
\usepackage{hyperref}
\usepackage{marginnote}
\usepackage[shortlabels]{enumitem} 
\numberwithin{equation}{section}

\usepackage{accents}
\usepackage{dsfont} 

\usepackage{color}
\def\R{\mathds{R}}

\def\wto{\rightharpoonup}

\def\Div{\textup{div}}
\def\dist{\textup{dist}}

\def\Chi#1{\hbox{{\large $\chi$}{\Large $_{_{#1}}$}}}

\renewcommand{\H}{\mathcal{H}}
\newcommand{\beq}{\begin{equation}}
	\newcommand{\eeq}{\end{equation}}
\newcommand{\pa}{\partial}
\newcommand{\divr}{\mathfrak{Div}^r}
\newcommand{\divs}{\textup{div}_s}
\newcommand{\vphi}{\varphi}
\newcommand{\noop}[1]{} 
\newcommand{\osc}{\operatorname*{osc}}

\newcommand{\Lip}{\mathrm{Lip}}

\newcommand{\J}{J_r}
\renewcommand{\P}{\textnormal{Per}_r}
\newcommand{\Ps}{\mathcal{P}_s}

\newcommand{\td}{\tilde{d}}
\theoremstyle{plain}
\newtheorem{theorem}{Theorem}[section]
\newtheorem{proposition}[theorem]{Proposition}

\newtheorem{lemma}[theorem]{Lemma} 
\theoremstyle{definition}
\newtheorem{definition}[theorem]{Definition}
\newtheorem{remark}[theorem]{Remark}

\newcommand{\dd}{d}
\newcommand{\N}{\mathds N}

\newcommand{\res}
{\mathop{\hbox{\vrule height 7pt width .5pt depth 0pt
			\vrule height .5pt width 5pt depth 0pt}}\nolimits}
\newcommand{\medint}{-\kern -,375cm\int}
\newcommand{\medintinrigo}{-\kern -,315cm\int}

\usepackage[foot]{amsaddr}

\allowdisplaybreaks

\author{F. Cagnetti $^{1}$}
\email{filippo.cagnetti@unipr.it}

\author{M. Morini $^{1}$}
\email{massimiliano.morini@unipr.it}

\author{D. Reggiani $^{2}$}
\email{dario.reggiani@uni-muenster.de}

\address{$^1$Universit\`a di Parma, Dipartimento di Scienze Matematiche, Fisiche e Informatiche, Parma, Italy}
\address{$^2$Applied Mathematics M\"unster, University of M\"unster, Einsteinstrasse 62, 48149 M\"unster, Germany}

\begin{document}
	
	\title{A distributional approach to nonlocal curvature flows}
%
%
	
	\maketitle

	\begin{abstract}
		{\rm  In \cite{CMP17} a novel distributional approach has been introduced to provide a well-posed formulation of a class of crystalline mean curvature flows. In this paper, such an approach is  extended  to the nonlocal setting. Applications include the fractional mean curvature flow and the Minkowski flow; i.e., the geometric flow generated by the $(N-1)$-dimensional Minkowski pre-content.}
	\end{abstract}
	
	\tableofcontents
	
			\section{Introduction}
		
Nonlocal geometric functionals and the associated geometric flows have gained increasing attention in recent years, as replacements for the classical (local) perimeter in various applications. Unlike classical approaches, which extract boundary information from infinitesimal neighborhoods, nonlocal frameworks can integrate interactions over finite (or even large) regions, thus capturing more global structural details. This broadened perspective proves valuable in contexts ranging from Image Processing to Machine Learning. As an example, in Image Processing, Gilboa and Osher \cite{GO} introduced nonlocal methods to overcome some of the limitations related to local filtering and segmentation.  

In this paper we will consider two notable instances of such functionals:
a Minkowski pre-content type \emph{nonlocal perimeter}, and the \textit{$s$-fractional perimeter}.

\subsection{Minkowski pre-content type nonlocal perimeter} 
The Minkowski pre-content type \emph{nonlocal perimeter} is given by the (normalised) Lebesgue measure of the $r$-tubular neighborhood of the boundary, where $r>0$ is a fixed scale.  Specifically, for a measurable set $E \subset \R^N$ and a fixed radius $r>0$, we consider the  (oscillation) nonlocal perimeter  defined by
\begin{equation}\label{eq:nonlocal_perimeter_def}
    \P (E) := \frac{1}{2r} \int_{\R^N}\osc_{B_r(w)} \chi_E \, dw, \,
\end{equation}
and the corresponding \emph{nonlocal total variation}  
\begin{equation}\label{eq:nonlocal_TV_def}
    J_r (u) := \frac{1}{2r} \int_{\R^N}\osc_{B_r(w)} u \, dw, \,
\end{equation}
defined for any $u\in L^1_{loc}(\R^N)$.
Here, $B_r(w)$ is the ball of radius $r$ centered at $w$, $\chi_E$ denotes the characteristic function of $E$, and $\osc_{B_r(w)}  u$ stands for the   (essential) oscillation of $u$ over  $B_r(w)$. Intuitively, this framework incorporates scale-dependent effects by assigning negligible energy cost to boundary irregularities occurring at small (compared to $r$) scales, while reproducing classical perimeter-like behavior when $E$ is sufficiently smooth and significantly larger than 
$r$. As a consequence,   variational models for denoising based on this type of energy  can selectively smooth out noise while leaving finer, yet well-structured, oscillations almost intact.

These features have been effectively leveraged in image processing tasks, where one frequently desires to remove random noise without destroying small-scale yet meaningful details such as textures, edges, or patterns, see \cite{BKLMP} (see also  \cite{CDNV18} for some rigorous analysis on (local) minimisers of the nonlocal perimeter). For similar reasons, this nonlocal perimeter has also been adopted in medical imaging—for example, to segment retinal blood vessels —leading to a method commonly referred to as the \emph{Infinite Perimeter Active Contour Model} (IPAC), see for instance \cite{IPAC} and the references therein. The same nonlocal perspective has also proved effective in adversarial learning contexts:   by penalizing transitions that occur over extended regions, one can promote spatial coherence in classification models, thereby enhancing their resilience against adversarial attacks and contributing  to more robust decision boundaries, see for instance \cite{BS24, GTM22}.

Associated with the nonlocal perimeter \eqref{eq:nonlocal_perimeter_def} is the notion of nonlocal mean curvature, obtained formally by taking the first variation of $\P(E)$. As computed in \cite[Lemma~2.1]{CMP12}, it  is given by
$\mathcal{K}_r(x, E) = (\mathcal{K}_r^+ +\mathcal{K}_r^-)(x, E)$, where 
\begin{equation} \label{our curvature}
\mathcal{K}_r^\pm(x,E): =\begin{cases}
 \pm \frac{1}{2 r}\det (I\pm r \nabla \nu_E(x))\  &
\textrm{ if }  \dist(x\pm r \nu_E(x),\partial E)= r\,, \\
0 & \textrm{ otherwise,}
\end{cases}
\end{equation}
provided  that $E$ is of class $C^2$ and  that $r_{max}(x):=\max\{t>0:  \dist(x+ t \nu_E(x),\partial E)= t \text{ and } \dist(x- t \nu_E(x),\partial E)= t\}\neq r$ for $\H^{N-1}$-a.e. $x\in \pa E$.  Here $\nu_E(x)$ denotes the outward unit normal to $\pa E$ at $x$.

In this paper we consider the associated  nonlocal geometric motion 
\beq\label{lei}
V=-\mathcal{K}_r,
\eeq
referred to as the \emph{Minkowski curvature flow}, in which the outward normal velocity $V$ of the evolving set boundary is given by minus the nonlocal curvature $\mathcal{K}_r$. The expected behavior is comparable to that of the classical curvature flow for large, smooth sets, because in those cases the nonlocal curvature behaves similarly to the classical one (and indeed can coincide for sufficiently regular sets in low dimensions). In contrast, small or isolated sets, which are perceived as high-curvature regions, tend to shrink and disappear rapidly. On the other hand, if these smaller sets are part of a larger structure (e.g., a repeating texture like fingerprint ridges), local modifications to these subsets have minimal influence on the overall energy, implying that such subsets often persist as long as the surrounding configuration remains intact. We refer to \cite{DNV19} for some interesting qualitative rigorous results on the the behaviour of the Minkowski flow in the plane. 

In \cite{CMP15} a unified theory for a large class of (translation invariant) nonlocal curvature motions  has  been developed, including  existence and uniqueness of level-set formulations via a viscosity solution approach. There, a crucial hypothesis is that the nonlocal curvature is well-defined for sets with $C^2$ boundaries, since the viscosity theory  relies on $C^2$ test functions and requires the nonlocal curvature to be defined on their  level sets. However, as shown in \cite{CMP12},  for certain $C^2$ domains the oscillation  perimeter does not yield a well-defined curvature, thus preventing direct application of the general theory of \cite{CMP15}. To circumvent this issue, a regularized version of the nonlocal perimeter (recalled in Section~\ref{sec:variant}) has been proposed, ensuring that the corresponding curvature remains well-defined for sufficiently regular sets while retaining many of the essential nonlocal features, see \cite{CMP12, CMP15}. 
We refer to \cite{BarCar08, BarCar09, BarImb08, BarSou98, Car00, Car01, Car07, Car08, Car06, Imbert_2009, Slep03} for some important results concerning nonlocal geometric motions.

In this paper, we address the original (unregularized) equation \eqref{lei} by developing a nonlocal variant of the \emph{distance distributional approach}, originally introduced in \cite{CMP17} and \cite{CMNP19} (see also \cite{CMNP19-}) to establish the well-posedness of non-smooth (potentially crystalline) anisotropic mean curvature flows, where the viscosity approach fails. 

\subsection{The distance function approach.} 
To motivate the aforementioned distance distributional approach, let us recall a classical result by Soner \cite{S93} (see also \cite{BarSonSou93, EvSonSou92}).
If $\psi$ is a smooth anisotropy and $m$ is a norm (playing the role of a mobility), then (non-fattening) weak solutions $E(t)$ to the  anisotropic mean curvature flow 
\beq\label{flow}
V=-m(\nu) H_\psi\,,
\eeq
with $H_\psi$ denoting the anisotropic mean curvature associated with $\psi$, are characterized  by the inequalities
\beq\label{soner}
\begin{cases}
\partial_td\geq \Div(\nabla \psi(\nabla d)) & \text{in }\{d>0\}\\
\partial_td\leq \Div(\nabla \psi(\nabla d)) & \text{in }\{d<0\}
\end{cases}
\qquad\text{in the viscosity sense,}
\eeq
where $d(x, t):=\inf_{y\in E(t)}m^o(y-x)-\inf_{y\not \in E(t)}m^o(y-x)$ is the signed distance function to $E(t)$ induced by the dual norm $m^o$ of $m$. Here $\nabla$ and  $\Div$ clearly stand for the spatial gradient and divergence, respectively. 

Note now that when $\psi$ is non-smooth or even crystalline the inequalities \eqref{soner} do not  make sense anymore, due to the obvious fact that for some directions $\nu$ the \emph{Cahn-Hoffman} vector $\nabla\psi(\nu)$ might not be well-defined. In this situation the subdifferential $\pa\psi(\nu)$ becomes multi-valued. In order to recover \eqref{soner} one would have to look at  suitable Cahn-Hoffman selections $z(x,t)$ of the subdifferential map $\pa \psi(\nabla d(x,t))$; i.e, vector fields $z$ such that $\|z\|_\infty\leq 1$ and $z(x,t)\in \pa \psi(\nabla d(x,t))$ for a.e. $(x,t)$. It is not difficult to see that this is equivalent to 
\beq\label{CH}
\|z\|_\infty\leq 1\quad\text{and}\quad \int\int\eta \nabla d\cdot z\, dxdt=\int\int \eta  \psi(\nabla d) \, dxdt,
\eeq
for all smooth $\eta$  with compact support in space-time. 
Any such $z$ will be referred to as a \emph{Cahn-Hoffman field for $d$}. The idea of \cite{CMP17, CMNP19} is to reformulate \eqref{soner} in a \emph{distributional sense}; namely, to define weak solutions of the flow  \eqref{flow} essentially by requiring 
\beq\label{soner2}
\begin{cases}
\partial_td\geq \Div z \quad \text{ in }\{d>0\}\vspace{-0.3cm}\\
& \hspace{-2cm}\text{in the distributional sense,}\vspace{-0.3cm}\\
\partial_td\leq \Div z \quad \text{in }\{d<0\}\vspace{0.5pt}\\
z \text{ is a suitable  Cahn-Hoffman field for }d,
\end{cases}
\eeq
see  \cite{CMP17, CMNP19} for the precise details. See also \cite{CaCha}, where the above inequality is derived in the context of anisotropic (possibly crystalline) mean curvature flow of convex sets.  In \cite{CMP17, CMNP19} it is shown that such a notion yields global-in-time existence and uniqueness up to fattening and well-posedness of level set flows, provided a suitable compatibility condition between $m$ and $\psi$ is satisfied. Moreover, such distributional solutions can be obtained, up to fattening, as the limit of the Almgren-Taylor-Wang minimising movement scheme introduced in \cite{ATW, LS}. In particular, they coincide with the weak solution produced by the viscosity methods whenever the latter can be applied.  We refer to \cite{GP16, GP18} and to the nice survey \cite{GP22} for a completely different \emph{crystalline  viscosity approach} to crystalline flows. The latter method is restricted to purely crystalline anisotropies but in this context it is very powerful as it allows for instance to consider also geometric equations with a nonlinear dependence on the  curvature. 

\subsection{The distance function approach in the nonlocal setting.} Turning back to \eqref{lei}, the idea is to provide an analogous distributional formulation in terms of the appropriate \emph{nonlocal divergence} operator. To identify the latter, given $\eta$  smooth  with compact support in $\R^N\times(0,T)$ and  $u$ locally bounded and continuous in the space variable $w$,  denoting  $C_r:=\overline{B_r}\times \overline{B_r}$ and letting  $\mathcal{M}^{+}(C_r)$ be the set non-negative Radon measures on $C_r$,  we can  formally compute: 
\begin{align}
&\frac1{2r}\int_0^T\int_{\R^N}\eta(w,t)  \operatorname*{osc}_{B_r(w)} u \, dw  dt \nonumber \\
&= \frac1{2r}   \int_0^T\int_{\R^N}\eta(w,t)  \sup_{(x,y)\in C_r}| u(t,x+w)-u(t,y+w)| \, dwdt \nonumber \\ 
&=\frac1{2r}\int_0^T\int_{\R^N} \eta(w,t) \sup_{\substack{z\in {\mathcal M}^{+}(C_r)\\[1pt] |z|(C_r)\leq 1}} \int_{C_r}u(t,x+w)-u(t,y+w) dz(x,y)\, dwdt \label{formal}\\
&=\sup_{\substack{z\in L^\infty(\R^N\times(0,T); {\mathcal M}^{+}(C_r))\\[1pt] \|z\|_\infty\leq 1}} \frac1{2r}\int_0^T\int_{\R^N} \eta(w,t)\int_{C_r}u(t,x+w)-u(t,y+w) dz_{(w,t)}(x,y)\, dwdt, \nonumber
\end{align}
where $z_{(w,t)}$ denotes the evaluation of the measure-valued map $z$ at $(w,t)$. Therefore, it is natural to consider (measurable) maps $z: \mathbb{R}^N\times (0,T) \to \mathcal{M}^+({C}_r)$ as the nonlocal analogue of vector fields in this context. In analogy with \eqref{CH}, the optimal nonlocal fields for which the sup in the above formula is attained (for all $\eta$)  will  play the role of generalised \emph{nonlocal Cahn-Hoffman fields} for the function $u$. 
Finally, motivated by the  last integral in \eqref{formal}, given a nonlocal  Cahn-Hoffman field $z$ and any smooth test function $\vphi$ with compact support, we are led to define the \emph{nonlocal divergence} $\mathfrak{Div}^r z$ of $z$ by:
\[
\langle \mathfrak{Div}^r z, \vphi \rangle := - \frac{1}{2r} \int_0^T\int_{\mathbb{R}^N} \int_{\mathcal{C}_r} \big( \vphi(t,x+w) - \vphi(t,y+w) \big) \, d z_{(w,t)}(x,y) \, dwdt.
\]
Heuristically, when $z$ is a Cahn–Hoffman field associated with a (Lipschitz) function  $u$, the quantity $\mathfrak{Div}^{r} z$ can be interpreted as the nonlocal  curvature of the level sets of  $u$. 
With these notions in place, the idea is now to define weak solutions of \eqref{lei} by enforcing again \eqref{soner2} with $\Div$ and $z$ replaced by the nonlocal notions that we have just introduced. 
Roughly speaking, these correspond to the evolution of closed sets $E(t)$ (respectively,  open sets $A(t)$), where the associated distance function $d(\cdot, t) = \dist(\cdot, E(t))$ (respectively, $d(\cdot, t) = \dist(\cdot, A(t))$) satisfies the inequality $\partial_t d \geq \mathfrak{Div}^r z$ (respectively,  $\partial_t d \leq \mathfrak{Div}^r z$) in the sense of distributions in the region $\{d > 0\}$ 
(respectively, $\{d < 0\}$) up to the extinction time,  with $z$ being a suitable nonlocal Cahn-Hoffman field for $d$. We then  call a \emph{weak flow} for \eqref{lei} any weak  superflow whose interior part is a weak subflow. 
We refer to Section~\ref{sec:weakformulation} for the precise definition.

The main results of Sections~\ref{sec:2} can be summarized as follows: a comparison principle between superflows and subflows (Theorem~\ref{thm: wfcomparison}); the subsequential convergence of the Almgren-Taylor-Wang minimizing movement scheme (adapted to the present setting) to weak superflows and subflows (Theorem~\ref{themthm}); and the generic (up to fattening) global-in-time existence  and uniqueness of weak flows originating from closed sets with compact boundaries, along with the well-posedness of the corresponding level-set flows (Theorem~\ref{th:phiregularlevelset}). 
The general strategy for proving Theorem~\ref{themthm} and the comparison principle follows the path set forth  in \cite{CMP17}. However, many of the technical implementations differ and pose new challenges, due to the strong nonlocality of the underlying divergence operator and the very different structure of the associated nonlocal  Cahn-Hoffman fields. Among the key intermediate steps, we mention the analysis of the so-called Chambolle algorithm— a level-set reformulation of the minimizing movement scheme (see \cite{Chambolle})—within this nonlocal framework. This entails studying the nonlocal equation  
\[
\begin{cases}
-\mathfrak{Div}^r z + u = g  \text{ in $\mathbb{R}^N$ \, \, (in the sense of distributions)}, \\
z \text{ is a nonlocal Cahn-Hoffman field for } u, 
\end{cases}
\]
for a given function $g$ that is Lipschitz continuous and coercive (see Section~\ref{sec:ROF} and, in particular, Theorem~\ref{give me a name}). Additionally, a crucial step in our analysis involves estimating the discrete velocity of balls, which requires a delicate nonlocal calibration argument (see Proposition~\ref{prop no name} and Section~\ref{sec:balls}).

 \subsection{The $s$-fractional perimeter} 
We now move to the final part of the paper. In Section~\ref{sec:fractional}, to illustrate the flexibility of our method, we show that the distance distributional approach also extends to the \emph{fractional mean curvature flow}, i.e., the geometric flow associated with the $s$-fractional perimeter \eqref{pesse}, see \eqref{oeef}. Fractional perimeters have been the focus of extensive research in recent years, particularly in connection with nonlocal minimal surfaces, following the pioneering work \cite{CRS10}. The existence of weak solutions for fractional mean curvature flows has been studied in \cite{Imbert_2009, CS10, CMP15} within the framework of level set viscosity solutions. For results concerning the existence, uniqueness, and qualitative properties of smooth fractional flows, we refer to \cite{JLM20, CSV18}.

 In Section~\ref{sec:fractional}, we introduce a new distance distributional formulation in the same spirit as \eqref{soner2} and the one presented for the Minkowski flow. This formulation incorporates the appropriate notions of nonlocal divergence and the nonlocal Cahn-Hoffman field, both of which have already been identified in the literature—for instance, in the context of the fractional total variation flow (see, for example, \cite{AMRT08}). In this setting, we establish generic global-in-time existence and uniqueness, the well-posedness of the corresponding level set flows, and approximability via the minimizing movement algorithm. The latter property ensures that the new distributional formulation yields the same solutions as those obtained through the viscosity methods mentioned above. Since many of the arguments overlap with those presented in Section~\ref{sec:2}, we provide detailed proofs only for the parts that differ.

 \subsection{Final remarks} 
We conclude the introduction by highlighting some directions for future research. We expect that \eqref{lei} converges to the classical mean curvature flow as $r\to 0$, as proven in \cite{CDNP21} for the regularized version mentioned earlier. Note that a related problem is addressed in \cite{BLS24}, where the asymptotic behavior is analyzed in the context of time-discrete evolutions, considering the simultaneous vanishing of both the time step and the nonlocal scale $r$ at the same rate. It would also be interesting to explore anisotropic variants of \eqref{lei} (and their asymptotic behavior), where the underlying perimeter is defined as
\[
\P^K (E) := \frac{1}{2r} \int_{\R^N} \osc_{ w +rK} \chi_E \, dw,
\]
with $K$ being a given compact (possibly non-smooth) convex set in $\R^N$ that contains a neighborhood of the origin. These questions will be the subject of future investigations.

	\section{Preliminaries}\label{sec:preliminaries}
	
	%

In this section we introduce some notation and we recall some basic notions of measure theory. 
For more details, we direct the reader to the monographs \cite{AFP, maggiBOOK}.
Throughout the paper $N \in \mathbb{N}$ is fixed, and we use $\Lip (\R^N)$ to denote the space of Lipschitz continuous functions of $\R^N$, while $\Lip_c(\R^N)$ stands for the set of Lipschitz continuous functions of $\R^N$ with compact support. Moreover, 
we define the sets $\mathscr{L}^{+}$ and $\mathscr{L}^{-}$ as
$$
	\mathscr{L}^{\pm} := \Big\{ g \in \Lip(\R^N) : 
	\lim_{|x| \to \infty } g (x) = \pm \infty \Big\}.
$$
	Moreover, we denote by $\mathscr{M}(\R^N)$ the $\sigma$-algebra of Lebesgue measurable subsets 
	of $\R^N$. 
	For $A \in \mathscr{M}(\R^N)$, we will use both $\mathcal{L}^N (A)$ and $| A |$ to denote the $N$-dimensional Lebesgue measure of $A$, while the characteristic function $\chi_A: \R^N \to \R$ of $A$ is defined as
	\[
	\chi_A (x) = 
	\begin{cases}
	1 & \text{ if } x \in A, \\
	0 & \text{ if } x \notin A.
	\end{cases}
	\] 
The precise representative $A^{(1)}$ of $A$ is the set of points of density $1$ of $A$, 
defined by
 \begin{equation} \label{def: A1}
 \chi_{A^{(1)}} (x) = 1 \quad \Longleftrightarrow \quad
\lim_{\rho \to 0^+} \frac{| A^{(1)} \cap B_\rho (x)|}{|B_\rho (x)|} = 1. 
\end{equation} 
It turns out that $\mathcal{L}^N (A \Delta A^{(1)}) = 0$, where $\Delta$ stands for the symmetric difference of sets.
 If $A, B \in \mathscr{M}(\R^N)$, we write $A \subset_{\mathcal{L}^N} B$ and $A =_{\mathcal{L}^N} B$ 
when $\mathcal{L}^N (A \setminus B) = 0$ and $\mathcal{L}^N (A \Delta B) = 0$, 
respectively. 
If $f, g: \R^N \to \R$ are measurable functions, we define $f \vee g$ and $f \wedge g$
as the two measurable functions given by
\[
(f \vee g ) (x) := \max\{ f (x), g (x)\}, 
\quad \text{ and } \quad 
(f \wedge g ) (x) := \min\{ f (x), g (x)\}, 
\] 
for every $x \in \R^N$, and we set $f^+:= f \vee 0$, $f^-:= (-f)^+$, and
	$$
	\R + L^2 (\R^N): = \{ g: \R^N \to \R \text{ measurable such that } g - M \in L^2(\R^N) \text{ for some } M \in \R \}.
	$$
If $A \subset \R^N$, we use $\text{Int}_{\R^N}(A)$ to denote the interior of $A$ \textit{in the standard topology of $\R^N$}.
	If $E \subset \R^N \times [ 0, + \infty)$, we set $E^c = \R^N\times [0,+\infty) \setminus E$, and we denote by $\text{Int}(E)$ and $\overline{E}$ the interior and the closure of $E$ \textit{in the relative topology of $\R^N\times [0,+\infty)$}, respectively. 
If $C, D \subset \R^N$ are non empty closed sets, we use the standard notation for the distance function:
\[
\dist(w, C) = \min\left\{ | w - x | : x \in C  \right\} \, \, \, \, \forall  w \in \R^N, \, \, \text{ and } \, \,
\dist(C,D) = \min\left\{ | x - y | : x \in C, y \in D  \right\},
\] 
while $d_D$ is the signed distance function 
from $\pa D$, that is,
	$$
	d_D(w):=
	\begin{cases}
		\,\dist(w, \pa D) & \text{if }w\in \R^N \setminus D\,,\\
		- \dist(w, \pa D) & \text{if }w\in D\,.
	\end{cases}
	$$	
We also recall that a sequence of closed sets  $\{D_n\}$ of $\R^d$ converges to a closed set $E\subset \R^d$ in the Kuratowski sense, and we write $D_n \stackrel{\mathscr K}{\longrightarrow} D$, if the following two conditions hold: 
	\begin{enumerate}[label=(\roman*)]
		\item  if  $x_n\in D_n$ for all $n$, then any cluster point of $\{x_n\}_n$ belongs to $D$;
		\item any $x\in D$ is the limit of a sequence $\{x_n\}_n$, with $x_n\in D_n$ for all $n$.
	\end{enumerate}
It turns out that $D_n \stackrel{\mathscr K}{\longrightarrow} D$ if and only if 
\[
 \dist(\cdot,  D_n) \to \dist(\cdot,  D) \qquad \text{ locally uniformly in } \R^d. 
\]
From this, thanks to Arzel\`a--Ascoli Theorem, it follows that every sequence of closed sets
admits a subsequence that convergences in the Kuratowski sense.

Let  $A \subset \R^N$ be a Borel set. We denote by \(  \mathcal{M} (A) \)  (\(  \mathcal{M}^+ (A) \)) the set of (positive) Radon measures in $A$.
If $\mu \in \mathcal{M} (\R^N)$ and $V \subset \subset \R^N$ is a Borel set, 
$\mu  \res V \in \mathcal{M} (\R^N)$ is the Radon measure  of \( \R^N \)
such that $( \mu  \res V ) (B) = \mu (B \cap V)$ for every Borel set $B \subset \subset \R^N$.

For every $r > 0$ and $w \in \R^N$, $B_r (w)$ stands for the open ball of $R^N$ centred at $w$ and with radius $r$
(we use $B_r$ in the special case $w=0$), and we set 
\beq\label{eq:Cr}
C_r = \overline{B_r} \times \overline{B_r}.
\eeq
For every $E \in \mathscr{M}(\R^N)$, the $r$-neighborhood of $E$ is defined as $(E)_r:=\cup_{x\in E}B_r(x)$.
Note that, in particular, with this definition we have $(\emptyset)_r = \emptyset$.
	If $\mu \in \mathcal{M} (C_r)$, we denote by $\mu^+$ and $\mu^-$ the positive and negative parts of $\mu$, respectively, 
	while $|\mu|$ is the total variation of $\mu$. We recall that if $\mu \in \mathcal{M} (C_r)$, 
	then $\mu^+, \mu^-, |\mu| \in \mathcal{M}^+ (C_r)$, and we have $\mu = \mu^+ - \mu^-$ and $| \mu | = \mu^+ + \mu^-$. 
	For every $w \in \R^N$, the Dirac delta measure concentrated at $w$ is denoted by $\delta_w$.
If $a, b \in \overline{B_r}$, the symbol $\delta_{a} \times \delta_{b}$ stands for the positive Radon measure  	
on $C_r$ such that
\begin{equation} \label{def misure}
		\int_{C_r} f (x, y) \, d (\delta_{a} \times \delta_{b}) (x, y) = f (a, b), \quad \text{ for every } f \in C^0 (C_r),
	\end{equation}
where $C^0 (C_r)$ is the set of continuous functions on $C_r$.
Since $C_r$ is compact, by Riesz Representation Theorem (see, for instance, \cite[Remark~1.44 and Theorem~1.54]{AFP}), 
the dual space $(C^0 (C_r))'$ of $C^0 (C_r)$ can be identified with $\mathcal{M} (C_r)$.
In turn, thanks to another version of Riesz Representation Theorem (see \cite[Theorem~2.112]{FonsecaLeoniBook}), the dual of 
 $L^1\left(\R^N; C^{0}\big(C_r\big)\right)$ can be identified with the space $L^\infty_{\mathbf{w}} \left(\R^N; \mathcal{M} (C_r)\right)$ defined as follows (see \cite[Definition~2.111]{FonsecaLeoniBook})):
\[
L^\infty_{\mathbf{w}} \left(\R^N; \mathcal{M} (C_r)\right) := 
\left\{ z: \R^N \to \mathcal{M} (C_r) \text{ such that (i) and (ii) are satisfied}  \right\}, 
\]
where $z_w$ denotes the evaluation of $z$ at $w$ and

\begin{itemize}

\item[(i)] $\displaystyle w \longmapsto \int_{C_r} f (x, y) \, d z_w (x, y)$ is measurable for every $f \in C^0 (C_r)$;

\vspace{.2cm}

\item[(ii)] $w \longmapsto | z_w | (C_r)$ belongs to $L^{\infty} (\R^N)$.

\end{itemize}

\vspace{.2cm}

\noindent
The space $L^\infty_{\mathbf{w}} \left(\R^N; \mathcal{M} (C_r)\right)$ is endowed with the norm
\[
\| z \|_{L^\infty_{\mathbf{w}} \left(\R^N; \mathcal{M} (C_r)\right)}
:= \operatorname*{ess\,sup}_{w \in \R^N} | z_w | (C_r).
\] 

\medskip

We are now ready to introduce the functionals we are interested in. We will consider two types of perimeter.
	
	\begin{enumerate}
		
		\item[(a)] For every $r > 0$ fixed, a Minkowski--type perimeter  is  given by
		 \begin{equation} \label{def: Per_r}
		 \P (E) := \frac{1}{2r} \int_{\R^N}\osc_{B_r(w)} \chi_E \, dw, \,
		\qquad \textnormal{ for every } E  \in \mathscr{M}(\R^N), 
		\end{equation} 
		where, for every $A \in \mathscr{M}(\R^N)$ and $v \in L^{1}_{\textnormal{loc}} (\R^N)$, we set 
		$$
		\osc_{A} v\, := \operatorname*{ess\,sup}_A v - \operatorname*{ess\,inf}_A v  \, \in [0, +\infty].
		$$ 
		It can be checked  that $\P(E)= \frac{1}{2r} \mathcal{L}^N\big((\pa E^{(1)})_r\big)$.
		
		\noindent
		 We also consider  the corresponding non--local total variation $J_r : L^{1}_{\textnormal{loc}} (\R^N) \to [0, + \infty]$ given by:
		\beq\label{gencoarea}
		J_r (v) :=  \frac{1}{2r} \int_{\R^N} \osc_{B_r (w)} v \, d w, 
		\eeq
		for every $v \in L^{1}_{\textnormal{loc}} (\R^N)$.

		\item[(b)] For every $s \in (0, 1)$ fixed, the $s$--fractional perimeter  is  given by
		\beq\label{pesse}
		\mathcal{P}_s (E) := c_s \int_{E} \int_{E^c} \frac{1}{|x - y|^{N+s}} dx dy, \, \qquad c_s := \frac{1-s}{\omega_{N-1}},
\eeq
			 for every set  $E \in \mathscr{M}(\R^N)$, with $\omega_{N-1}$ being the volume of the unit ball in $\R^{N-1}$. The choice of the constant $c_s$ is motivated by the fact that in this way  $\mathcal{P}_s$ approaches the classical perimeter as $s\nearrow 1$, see \cite{CV, ADPM} for details. 
		Also in this case, we consider the corresponding non--local total variation $\mathcal{J}_s: L^{1}_{\textnormal{loc}} (\R^N) \to [0, + \infty]$ given by
		\beq\label{TVf}
		\mathcal{J}_s (v) :=  \frac{c_s}{2} \int_{\R^N} \int_{\R^N} \frac{| v (x) - v (y)|}{|x - y|^{N+s}} dx dy, 
		\eeq
		for every $v \in L^{1}_{\textnormal{loc}} (\R^N)$.
		
	\end{enumerate}
	It can be checked that the following generalised coarea formulas hold  (see for instance \cite{CMP12, CMP15}): 
	\begin{equation} \label{coarea osc}
		\begin{split}
			&J_r (v) = \int_{-\infty}^{\infty} \P (\{ v > t \}) \, dt
			=  \int_{-\infty}^{\infty} \P (\{ v \leq t \}) \, dt, \\
			&\mathcal{J}_s (v) 
			= \int_{-\infty}^{\infty} \mathcal{P}_s (\{ v > t \}) \, dt
			=  \int_{-\infty}^{\infty} \mathcal{P}_s (\{ v \leq t \}) \, dt, 
		\end{split}
	\end{equation}
	for every $v \in L^{1}_{\textnormal{loc}} (\R^N)$.
	The following properties hold:
	
	\begin{enumerate}[label=(\roman*)]
		
		\item $J_r (\Chi{\emptyset}) = J_r (\Chi{\mathbb{R}^{N}}) = 0$ (so in particular $J$ is proper);
		
		\vspace{.1cm}
		
		\item $J_r ( \lambda u + c) = |\lambda | J_r (u)$ for every $u \in L^{1}_{\text{loc}} (\R^N)$ and for every $\lambda\, ,  c \in \R$;
		
				\vspace{.1cm}

		\item $J_r (u(\cdot + c)) = J _r(u (\cdot))$ for every $u \in L^{1}_{\text{loc}} (\R^N)$ and for every $c \in \R$; 
		
				\vspace{.1cm}

		\item  Subadditivity: 
		\begin{equation*} 
			J_r (u + v) \leq J_r (u) + J_r (v) \qquad \text{ for every } u, v \in  L^{1}_{\text{loc}} (\R^N);
		\end{equation*}
		
		\item Submodularity: 
		\begin{equation*} 
			J_r (u \vee v) + J_r (u \wedge v) \leq J_r (u) + J_r (v) \qquad \text{ for every } u, v \in 
			L^{1}_{\text{loc}} (\R^N);
		\end{equation*}
		
		\item  \( J_r ((u \land \mu ) \lor \lambda ) \leq J_r (u) \) for every $\lambda, \mu \in \overline{\R}$
		with $\lambda < \mu$, for every $u \in L^1_{loc} (\R^N)$. 
		
	\end{enumerate}
	The same properties hold with $J_r$ replaced by $\mathcal{J}_s$. 
	 Moreover, $J_r$ and $\mathcal{J}_s$ are convex, lower semicontinuous with respect to strong $L^p$
	convergence, and thus also with respect to weak $L^p$ convergence, for every $1\leq p < \infty$. The strong lower semicontinuity follows
	observing that $v \mapsto \osc_{B_r (w)} v$ is lower semicontinuous
	with respect to the strong convergence in $L^p_{loc} (\R^N)$, and then using Fatou's Lemma.
	
For the reader's convenience, let us show the submodularity for $J_r$. 
First of all, we claim that the coarea formula above, together with the convexity of $J_r$, implies the submodularity of $\P (\cdot)$:
\begin{equation} \label{submodularity P}
\P (E \cup F) + \P (E \cap F) \leq  \P (E) + \P (F) \qquad \text{ for every  measurable } E, F \subset \R^N.
\end{equation}
Indeed, observing that 
$$
\{ \chi_{E \cup F} > s \} = E \cup F = \left\{ \frac{\chi_{E} + \chi_{F}}{2} > s \right\}
\quad \text{ for every } s \in [ 0, 1/2), 
$$
and 
$$
\{ \chi_{E \cap F} > s \} = E \cap F = \left\{ \frac{\chi_{E} + \chi_{F}}{2} > s \right\}
\quad \text{ for every } s \in [1/2, 1), 
$$
we have
\begin{align*}
&\frac12\P (E \cup F) + \frac12\P (E \cap F) \\
&= \int_0^{1/2} \P ( E \cup F ) \, ds
+ \int_{1/2}^{1} \P ( E \cap F) \, ds \\
&= \int_0^{1/2} \P \Bigl(\Bigl\{\frac{\chi_E + \chi_F}{2} >s \Bigr\} \Bigr) \, ds
+ \int_{1/2}^1 \P \Bigl(\Bigl\{\frac{\chi_E + \chi_F}{2} >s \Bigr\} \Bigr) \, ds \\
&= J_r \Bigl( \frac{\chi_E + \chi_F}{2}  \Bigr) \leq 
\frac12 J_r ( \chi_E ) + \frac12  J_r ( \chi_F )
= \frac12 \P (E) + \frac12 \P (F),
\end{align*}
where in the third equality we used the coarea formula. This shows \eqref{submodularity P}. 
Now, thanks to coarea formula, we have  
\begin{align*}
J_r (u \vee v) =   \int_{-\infty}^{\infty} \P (\{ u \vee v > s \}) \, ds 
=  \int_{-\infty}^{\infty} \P (\{ u > s \} \cup \{ v > s \}) \, ds,
\end{align*}
and 
\begin{align*}
J_r (u \wedge v) =   \int_{-\infty}^{\infty} \P (\{ u \wedge v > s \}) \, ds 
=  \int_{-\infty}^{\infty} \P (\{ u > s \} \cap \{ v > s \}) \, ds.
\end{align*}
The last two identities, together with \eqref{submodularity P}, imply the submodularity of $J_r$.	
	%
		
	\section{The Minkowski--type mean curvature flow}\label{sec:2}

	\subsection{The generalised ROF problem}\label{sec:ROF}
	
	The aim of this section is to study, for $h > 0$ fixed,  the minimum problem 
	\beq \label{nlROF2}
	\min_{v\in L^{1}_{\text{loc}}(\R^N)}\biggl( \J (v)+ \frac{1}{2 h}\int_{\R^N}(v-g)^2\, dw\biggr),
	\eeq
	and the associated Euler-Lagrange equation.  
	We start by considering the case
	$g\in \R + L^2 (\R^N)$.  Recall  \eqref{eq:Cr} for the definition of $C_r$.
	
	\begin{proposition} \label{existence minimum}
		Let $h>0$. Then, for every $g\in \R + L^2 (\R^N)$ there exists a unique solution $u_h^g$ to the problem \eqref{nlROF2}, with
		$u_h^g \in \R + L^2 (\R^N)$.
		Moreover, if $g_1, g_2 \in \R + L^2 (\R^N)$ and  $g_1 \leq g_2$ a.e. in $\R^N$, then $u^{g_1}_h \leq u^{g_2}_h$ a.e. in $\R^N$.
		Finally, if $g \in \Lip (\R^N) \cap ( \R + L^2 (\R^N) )$, then $u_h^g \in \Lip (\R^N) \cap ( \R + L^2 (\R^N) )$ with 
		$\textnormal{Lip} (u_h^g) \leq \textnormal{Lip} (g)$, and 
		there exists $z\in L^\infty_{\mathbf{w}} \left(\R^N; \mathcal{M}^+(C_r)\right)$
		such that 
		\begin{itemize}
			\vspace{5pt}
			\item[(a)] for a.e. $w\in \R^N$ one has $ z_w \big(C_r\big) \leq 1$ and 
			\beq\label{(a)0}
			\int_{C_r}(u_h^g(w+x)-u_h^g(w+y))\, d z_w(x,y)= \osc_{B_r(w)}u_h^g\,;
			\eeq
			\item[(b)]  for every $\vphi\in C^0 (\R^N)\cap L^2(\R^N)$, with $\J(\vphi)<+\infty$
			$$
			\frac{h}{2r} \int_{\R^N} \int_{C_r}(\vphi(w+x)-\vphi(w+y))\, d z_w(x,y) dw = -  \int_{\R^N}\vphi(u_h^g-g)\, dw,
			$$ 
			\vspace{5pt}
		\end{itemize}
		where $z_w$ stands for the evaluation of $z$ at the point $w$, and $d z_w(x,y)$ for the integration in $(x,y)$ with respect to the measure $z_w$. 
	\end{proposition}
	\begin{proof}
		Since $g\in \R + L^2 (\R^N)$, there exists $M \in \R$ such that $g - M \in L^2 (\R^N)$. Testing with $v \equiv M$ and
		 observing that properties (i) and (ii) imply $J_r(M)=0$, we infer that the infimum of the functional in \eqref{nlROF2} is finite. Thus, letting  $(v_k) \subset L^{1}_{\text{loc}}(\R^N)$ be a minimising sequence, we have that the sequences $( v_k - g )$ and   $(v_k - M)$
		are bounded in $L^2(\R^N)$.
		Therefore, there exist $v \in \R + L^2 (\R^N)$ such that, up to subsequences,
		$$
		v_k - M \wto v - M \qquad \text{ weakly  in } L^2 (\R^N).
		$$
		Then, by lower semicontinuity 
		\begin{align*}
			&\liminf_{k \to \infty} \left( \J (v_k) +  \frac{1}{2 h} \int_{\R^N}(v_k -g)^2\, dw \right) \\
			&= \liminf_{k \to \infty} \left( \J (v_k - M) +   \frac{1}{2 h}\int_{\R^N}(v_k - M - (g - M))^2\, dw \right) \\
			&\hspace{.5cm}\geq \J (v-M) +  \frac{1}{2 h} \int_{\R^N}(v-g)^2\, dw
			= \J (v) +  \frac{1}{2 h} \int_{\R^N}(v-g)^2\, dw. 
		\end{align*}
		From the inequality above it follows that $v$ is a minimum point. 
		By strict convexity, the minimiser is unique.

		We now split  the remaining part of the proof into several steps.

		\smallskip
		
		\noindent
		\textbf{Step 1} (Comparison principle) The following argument is rather standard, but we include it for the reader's convenience.
		Let  $g_1, g_2 \in \R + L^2 (\R^N)$ and let us first assume that $g_1 < g_2$. 
		To ease the notation, we write $u_1$ and $u_2$ in place of $u_h^{g_1}$ and $u_h^{g_2}$, respectively.
		By minimality of $u_1$ and $u_2$ we have
		\begin{align*}
			\J (u_1)+  \frac{1}{2 h} \int_{\R^N}(u_1-g_1)^2\, dw
			\leq \J (u_1 \wedge u_2)+  \frac{1}{2 h} \int_{\R^N}(u_1 \wedge u_2-g_1)^2\, dw,
		\end{align*}
		and
		\begin{align*}
			\J (u_2)+  \frac{1}{2 h} \int_{\R^N}(u_2-g_2)^2\, dw
			\leq \J (u_1 \vee u_2)+  \frac{1}{2 h} \int_{\R^N}(u_1 \vee u_2-g_2)^2\, dw.
		\end{align*}
		Adding up the two inequalities above and using the submodularity of $J_r$ 
		it follows that 
		\begin{align*}
			0 &\leq \int_{\R^N} \left[ (u_1 \wedge u_2-g_1)^2 + (u_1 \vee u_2-g_2)^2 
			- (u_1-g_1)^2 - (u_2-g_2)^2 \right]  \, dw \\
			&= 2 \int_{\{ u_1 > u_ 2 \}} (u_1 - u_2) (g_1 - g_2)   \, dw.
		\end{align*} 
		Since we assumed $g_1 < g_2$ the above inequality implies that $| \{ u_1 > u_ 2 \} | = 0$, that is, 
		$u_1 \leq u_2$ a.e. in $\R^N$. 
		The general case $g_1 \leq g_2$ a.e. in $\R^N$ can be obtained
		replacing $g_2$ by $g_2 + \varepsilon$, so that  $u_1 \leq  u_2 + \varepsilon$ a.e. in $\R^N$, 
		since $u_2 + \varepsilon$  is the solution corresponding to $g_2 + \varepsilon$,
		and taking the limit $\varepsilon \to 0$. 
		

		\smallskip 
		
		\noindent
		\textbf{Step 2} (The case of a Lipschitz continuous datum $g$)
		Here we show that, if in addition $g \in \Lip (\R^N)$, 
		then $u_h^g  \in \Lip (\R^N)$ with $\textnormal{Lip} (u_h^g) \leq \textnormal{Lip} (g)$.
		This is also standard, and follows by the comparison principle and by the translation
		invariance. We sketch the argument for the reader's convenience.
		Let $v \in \R^N$ be fixed, and set $(\tau_v g ) (\cdot) : = g (\cdot + v)$.
		By translation invariance and by uniqueness, we have that 
		$u^{\tau_v g}_h  = \tau_v u^{g}_h$.
		On the other hand, setting $L:= \textnormal{Lip} (g)$, we have  
		$\tau_v g \leq g + L | v|$. Since trivially $u^{g + L |v|}_h = u_h^g + L |v|$, 
		by the comparison principle established above we have $\tau_v u^{g}_h = u^{\tau_v g}_h \leq u_h^g + L |v|$.
		Exchanging the role of $\tau_v g$ and $g$, by a similar argument one obtains also that
		$u_h^g - L |v| \leq u^{\tau_v g}_h = \tau_v u^{g}_h$, thus allowing us to conclude.

		\smallskip
		
		\noindent
		\textbf{Step 3}  (Euler-Lagrange equation for $u_h^g$)
		We now compute the Euler-Lagrange equation satisfied by $u_h^{g}$, under the additional  assumption $g \in \Lip (\R^N)$.

		For all $\vphi\in C^0 (\R^N)\cap L^2(\R^N)$ with $\J (\varphi) < \infty$ and $t \in (0, 1)$ one has by minimality
		\begin{align*}
			& \frac{1}{2 h} \int_{\R^N}(u_h^g-g)^2 \, dw
			- \frac{1}{2 h} \int_{\R^N}(u_h^g+t\vphi-g)^2\, dw \\
			&\leq \J(u_h^g+ t \vphi) - \J(u_h^g) 
			\leq t \bigg[  \J(u_h^g+ \vphi) - \J(u_h^g)  \bigg],
		\end{align*}
		where in the last inequality we have used the fact that the function 
		$t \mapsto \J (u + t \vphi)$ is convex. 
		Dividing the previous inequality by $t$, letting $t \to 0^+$ and using the subadditivity of $\J$, 
		we obtain that
		\begin{equation}  \label{**}
			\begin{aligned}
				- \frac{1}{h}  \int_{\R^N}\vphi(u_h^g-g)\, dw 
				\leq  \J(u_h^g+\vphi)-  \J(u_h^g)  \leq  \J(\vphi)  \\
				= \int_{\R^N}\|\Psi_{\vphi}(w; \cdot)\|_{C^{0}\left(C_r\right)}\, dw
				=\|\Psi_{\vphi}\|_{L^1\left(\R^N; C^{0}\left(C_r\right)\right)}\,, 
			\end{aligned}
		\end{equation}
		where we have set $\Psi_{\vphi}(w; x,y):=\frac{1}{2r}(\vphi(w+x)-\vphi(w+y))$.  
		Let us now consider the subspace $\mathcal{N}~\subset~L^1\left(\R^N; C^{0}\big(C_r\big)\right)$ defined as
		$$
		\mathcal{N}:=\left\{\Psi_\vphi:\, \vphi\in C^0(\R^N)\cap L^2(\R^N)\text{ with } \J(\vphi)<+\infty\right\}
		$$ and the linear operator $L: \mathcal{N} \to \R$ given by
		$$
		L(\Psi_\vphi):= - \frac{1}{h}  \int_{\R^N}\vphi(u_h^g-g)\, dw\,.
		$$
		By \eqref{**} it turns out that $L$ is continuous with $\|L\| \leq 1$, thus by Hahn-Banach Theorem we may extend it to an element $\tilde L$ of $ \left(L^1\left(\R^N; C^{0}(C_r)\right)\right)'$ with the same norm. 
		But by Riesz Representation Theorem (see \cite[Theorem~2.112]{FonsecaLeoniBook}) the dual of $L^1\left(\R^N; C^{0}\big(C_r\big)\right)$ can be identified with $L^\infty_{\mathbf{w}} \left(\R^N; \mathcal{M} (C_r)\right)$. 
		Therefore there  exists $\tilde{z} \in L^\infty_{\mathbf{w}} \left(\R^N; \mathcal{M} (C_r)\right)$ such that 
		$\| \tilde{z} \|_{L^\infty\left(\R^N; \mathcal{M} (C_r)\right)}=\| \tilde L \| \leq  1$ and
		\begin{equation}  \label{***}
			\begin{aligned}
				- \frac{1}{h}  \int_{\R^N}\vphi(u_h^g-g)\, dw 
				=L(\Psi_\vphi)=\int_{\R^N} \int_{C_r} \Psi_{\vphi}(w; x, y)\, d \tilde{z}_w  (x, y) dw  \\
				= \frac{1}{2r}\int_{\R^N} \int_{C_r}(\vphi(w+x)-\vphi(w+y))\, d \tilde{z}_w (x, y) dw,
			\end{aligned}
		\end{equation}
		%
		%
		for all  $\vphi\in C^0(\R^N)\cap L^2(\R^N)$ with $\J(\vphi)<+\infty$. 
		Let now, for each $w \in \R^N$, write  $\tilde{z}_w= \tilde{z}^{+} _w - \tilde{z}^{-}_w$, 
		where $\tilde{z}^{+}_w$ and $\tilde{z}^{-}_w$ are the positive and negative parts 
		of $\tilde{z}_w$, respectively. 
		Note that we can write 
		\begin{equation} \label{ahboh}
			\begin{aligned}
				\int_{\R^N} &\int_{C_r}(\vphi(w+x)-\vphi(w+y))\,  d \tilde{z}_w  (x,y) dw \\
				=& \int_{\R^N} \int_{C_r}(\vphi(w+x)-\vphi(w+y))\, d \tilde{z}^+_w  (x,y) dw \\
				&+ \int_{\R^N} \int_{C_r}(\vphi(w+y)-\vphi(w+x))\,d \tilde{z}^-_w  (x,y) dw \\
				=& \int_{\R^N} \int_{C_r}(\vphi(w+x)-\vphi(w+y))\, d \tilde{z}^+_w (x,y) dw \\
				&+ \int_{\R^N} \int_{C_r}(\vphi(w+x)-\vphi(w+y))\, d \hat{z}^{-}_w (x,y) dw, 
			\end{aligned}
		\end{equation}
		where $\hat{z}^{-}_w$ is the positive measure on $C_r$ defined by 
		$$
		\int_{C_r} f (x, y)\, d \hat{z}^{-}_w (x,y):= \int_{C_r} f (y, x)\, d \tilde{z}^{-}_w (x, y),
		$$
		for all $f \in C^0 (C_r)$.
		At this point, we can consider $z \in L^\infty_{\mathbf{w}} \left(\R^N; \mathcal{M}^+(C_r)\right)$ 
		given by $z_w:= \tilde{z}^{+}_w + \hat{z}^{-} _w$ 
		for a.e. $w \in \R^N$.
		Note that, by construction, $\| z \|_{L^\infty_{\mathbf{w}} \left(\R^N; \mathcal{M} (C_r)\right)} 
		= \| \tilde{z} \|_{L^\infty_{\mathbf{w}} \left(\R^N; \mathcal{M} (C_r)\right)}$. 
		Moreover, by \eqref{***} and  \eqref{ahboh}, 
		$$
		\frac1{2r}\int_{\R^N} \int_{C_r}(\vphi(w+x)-\vphi(w+y))\, d z_w (x,y) dw
		= - \frac{1}{h}  \int_{\R^N}\vphi(u_h^g-g)\, dw,
		$$
		for all  $\vphi\in C^0(\R^N)\cap L^2(\R^N)$ with $\J(\vphi)<+\infty$.  This shows property (b).  We now prove (a). To this aim, 
		let $M\in \R$ be as  at the beginning of the proof.  Choosing $\vphi= - (u_h^g - M)$ 
		in the identity above and using \eqref{**}, we get
		\begin{align*}
			- \frac1{2r}\int_{\R^N} 
			&\int_{C_r}(u_h^g(w+x)-u_h^g(w+y))\, d z_w(x,y) dw
			= \frac{1}{h} \int_{\R^N}(u_h^g - M)(u_h^g-g)\, dw\\
			&\leq  \J(M)-\J(u_h^g) =- \frac1{2r}\int_{\R^N}\osc_{B_r(w)}u_h^g\, dw\,.
		\end{align*}
		Using the fact that $ z_w \big(C_r\big) \leq 1$ for a.e. $w \in \R^N$, we deduce in particular that
		\begin{equation*}  
			\begin{aligned}
				\int_{\R^N}\osc_{B_r }u_h^g (w+\cdot)\, dw 
				\leq  \int_{\R^N} \int_{C_r}(u_h^g(w+x)-u_h^g(w+y))\, d z_w (x,y) dw \\
				\leq \int_{\R^N} \osc_{B_r}u_h^g(w+\cdot) \, z_w (C_r ) \, dw\leq  \int_{\R^N}\osc_{B_r}u_h^g(w+\cdot)\, dw.
			\end{aligned}
		\end{equation*} 
				Therefore, all the inequalities above are equalities, and thus
		\begin{equation*}  
			\begin{aligned}
			\int_{\R^N} f (w) dw = 0, \quad \text{ where } \quad f (w) := \osc_{B_r }u_h^g (w+\cdot) - \int_{C_r}(u_h^g(w+x)-u_h^g(w+y))\, d z_w (x,y).
			\end{aligned}
		\end{equation*}
Since $f$ is non negative, we obtain		
		$$
		\int_{C_r}(u_h^g(w+x)-u_h^g(w+y))\, d z_w(x,y)= \osc_{B_r(w)}u_h^g,
		$$ 
		for a.e. $w\in \R^N$. 
	\end{proof}
	Motivated by condition (b) of Proposition~\ref{existence minimum}, 
	we introduce the notion of \textit{nonlocal divergence} 
	of  $z\in L^\infty_{\mathbf{w}}\left(\R^N; \mathcal{M}^+(C_r)\right)$
	as the operator $\mathfrak{Div}^r z \in \mathcal{M} (\R^N)$ acting on functions
	as
	\beq \label{def div}
	\int_{\R^N} \vphi (w)  \, d (\mathfrak{Div}^r z) (w) := 
	-  \frac1{2r}\int_{\R^N} \int_{C_r}(\vphi(w+x)-\vphi(w+y))\, d z_w(x,y) \, dw
	\eeq
	for every $\vphi \in C^{0}_c (\R^N)$. Here, $\mathcal{M} (\R^N)$ denotes the space of Radon measures on $\R^N$. 
	From now on, we will refer to any $z\in L^\infty_{\mathbf{w}} \left(\R^N; \mathcal{M}^+(C_r)\right)$ 
	satisfying condition (a) of Proposition~\ref{existence minimum} as a
	\emph{generalised Cahn-Hoffman field} for $u^g_h$.

	We will now consider the case $g \in \mathscr{L}^{\pm}$, where we recall that
	$$
	\mathscr{L}^{\pm} := \Big\{ g \in \Lip(\R^N) : 
	\lim_{|x| \to \infty } g (x) = \pm \infty \Big\}.
	$$
	Note that in this case the functional in \eqref{nlROF2} could be identically equal to $+\infty$. Nevertheless, the Euler-Lagrange equation 
	computed in Proposition~\ref{existence minimum}
	still makes sense.
	In the next theorem, which is the other main result of the section, we establish existence and uniqueness 
	of a solution of this equation.
	\begin{theorem} \label{give me a name}
		Let $g \in \mathscr{L}^+$ ($g \in \mathscr{L}^-$) and let $h$ be a positive constant. 
		Then, there exist $z\in L^\infty_{\mathbf{w}} \left(\R^N; \mathcal{M}^+(C_r)\right)$ and a unique Lipschitz function $u_h^g$ in $\mathscr{L}^+$ ($\mathscr{L}^-$),  such that 
		\begin{itemize}
			\vspace{5pt}
			\item[(a)] for a.e. $w\in \R^N$ one has $ z_w \big(C_r\big) \leq 1$ and 
			\beq\label{(a)}
			\int_{C_r}(u_h^g(w+x)-u_h^g(w+y))\, d z_w(x,y)= \osc_{B_r(w)}u_h^g\,;
			\eeq
			
			\item[(b)]  $\displaystyle - h \, \mathfrak{Div}^r z + u_h^g = g \qquad \text{ in } \mathcal{M} (\R^N).$
			\vspace{5pt}

		\end{itemize}
		Moreover, $\textnormal{Lip} (u_h^g) \leq \textnormal{Lip} (g)$. 
		In addition, for every $t \in \R$ the sets  $A_t := \{ u_h^g < t \}$ ($A_t := \{ u_h^g > t \}$) and $E_t := \{ u_h^g \leq t \}$ ($E_t := \{ u_h^g \geq t \}$) are  the minimal and maximal solution (up to $\mathcal{L}^N$-negligible sets),  respectively, of the problem 
		\beq\label{ATW0}
		\min_{F\in \mathscr{M}(\R^N)} \left\{  \P (F) + \frac1h\int_{F} (g (w) - t) \, d w \right\}\quad 
		\left(\, \min_{F\in \mathscr{M}(\R^N)} \left\{  \P (F) + \frac1h\int_{F} (t-g (w)) \, d w \right\}\,\right).
		\eeq
		Finally, the following comparison principle holds:
		if $g_1, g_2 \in \mathscr{L}^+\cup \mathscr{L}^-$ and  $g_1 \leq g_2$, then $u^{g_1}_h \leq u^{g_2}_h$.
	\end{theorem}
	
	
	\begin{remark} \label{richiamare}
		It is easy to check that for every $g \in \mathscr{L}^+ (\mathscr{L}^-)$ and for every $M \in \R$ we have that 
		$$
		u_h^{M + g} = M + u_h^{g}.
		$$
	\end{remark}

	Before proving Theorem~\ref{give me a name}, we need to  establish  several auxiliary results.
	We start with the following lemma, describing the structure of generalised Cahn-Hoffman fields.  
	
	\begin{lemma} \label{structure z}
		Let $(u, z) \in (\Lip(\R^N)  \times  \mathcal{M}^+(C_r))$, with $z(C_r)\leq 1$, be such that 
		\beq\label{eqstructure z}
		\int_{C_r}(u(x)-u(y))dz(x,y)= \osc_{B_r} u\,. 
		\eeq
		Then, 		
		\beq\label{equality on the support of z}
		u (x) = \max_{\overline{B_r}} u \qquad \text{ and } \qquad u (y) = \min_{\overline{B_r}} u
		\qquad \textnormal{for $z$-a.e. $(x, y) \in C_r$.}
		\eeq
	\end{lemma}
	
	
	\begin{proof}[Proof of Lemma~\ref{structure z}]
		We have
		\begin{align*}
			\osc_{B_r} u &= \int_{C_r}(u(x)-u(y))\, dz(x,y)  
			\leq z(C_r)\,\osc_{B_r } u    \leq \osc_{B_r} u.
		\end{align*}
		Thus all the inequalities are equalities, and, in particular, $\displaystyle u (x) - u (y)= \osc_{B_r} u$ for $z$-a.e. $(x,y)\in C_r$.
		%
		%
		%
	\end{proof}
	\begin{remark}\label{rm:monotonicity}
		Note that if $(u, z) \in (\Lip(\R^N)  \times  \mathcal{M}^+(C_r))$, with $z(C_r)\leq 1$, satisfy \eqref{eqstructure z}, then by Lemma~\ref{structure z} identity \eqref{equality on the support of z} holds with $u$ replaced by $u\land \lambda$ and  $u\lor \lambda$, for any $\lambda\in \R$. 
	\end{remark}
	\begin{remark}\label{rm:monotonicity2}
		Note that the assumptions of Lemma~\ref{structure z}  are equivalent to saying that $z$ is in the subdifferential of the convex functional 
		\[
		C^0 (C_r) \ni f \longmapsto \| f \|_{C^0 (C_r)}
		\]
		at the function $f (x, y) = u (x) - u (y)$. Thus, 
		if for  $i=1,2$  $(u_i, z_i) \in (\Lip(\R^N)  \times  \mathcal{M}^+(C_r))$, with $z_i(C_r)\leq 1$, satisfy \eqref{eqstructure z} 
		(with $(u_i, z_i)$ in place of $(u,z)$),  then by convexity
		$$
		\int_{C_r} \left[  (u_1 - u_2)(x)- (u_1 - u_2)(y) \right]\, d (z_1- z_2) (x,y)\geq 0\,.
		$$
	\end{remark}
	
	In fact, the following stronger monotonicity property holds:
	\begin{lemma}\label{lm:monotonicity}
		Let $(u_i, z_i) \in \Lip(\R^N)  \times  \mathcal{M}^+(C_r)$, $i=1,2$ such that $z_i(C_r)\leq 1$ and  
		$$
		\int_{C_r}(u_i(x)-u_i(y))\,dz_i(x,y)=\osc_{B_r} u_i.
		$$
		Let $\psi, g_1, g_2: \ \R \to \R $ be continuous monotone nondecreasing, and set 
		\( v_{12}:=\psi( g_1 (u_1) - g_2 (u_2))  \).
		Then,  
		$$
		\int_{C_r}(v_{12}(x)-v_{12}(y)) d(z_1-z_2)(x,y)\geq 0.
		$$		
	\end{lemma}
	\begin{proof}
		We have
		\begin{align}
			&\int_{C_r} \left[  (v_{12}(x)- v_{12}(y) \right]\, d (z_1 - z_2 ) (x,y) \nonumber \\
			&= \int_{C_r} \left[  (v_{12}(x)- v_{12}(y) \right]\, d z_1  (x,y)
			- \int_{C_r} \left[  (v_{12}(x)- v_{12}(y) \right]\, dz_2 (x,y). \label{eqmono1}
		\end{align}
		Let us now set, for $i = 1, 2$, $v_i := g_i (u_i)$, 
		$M_i := \max_{\overline{B}_r} v_i$  and  $m_i := \min_{\overline{B}_r} v_i$.
		We observe that, if $z_1(C_r) < 1$ and $z_2 (C_r) < 1$, then 
		$M_1 = m_1$ and $M_2 = m_2$, so that $v_{12}$ is constant 
		in $C_r$, and the statement is trivially satisfied.
		From now on, we will then assume that $\max\{ z_1 (C_r) , z_2 (C_r) \} = 1$.
		
		Thanks to Lemma~\ref{structure z}, we have 
		$$
		v_i (x)  = M_i \qquad \text{ and } \qquad v_i (y) = m_i
		\qquad \textnormal{for $z_i $-a.e. $(x, y) \in C_r$, \quad for $i = 1, 2$}.
		$$
		Let us focus on the first integral in the right-hand side of \eqref{eqmono1}.
		By the above considerations, and using the fact that $\psi$, $v_1$ and $v_2$
		are continuous functions, we have that there exists 
		$(\overline{x}_1, \overline{y}_1) \in C_r$ such that
		\begin{align*}
			\int_{C_r} \left[  (v_{12}(x)- v_{12}(y) \right]\, dz_1  (x,y) 
			&= \int_{C_r} \left[  (
			\psi ( v_1 (x) - v_2 (x)  )
			- \psi ( v_1 (y) - v_2 (y)  ) \right]\, d z_1  (x,y) \\
			&=  \left[  
			\psi ( M_1  - v_2 (\overline{x}_1)  )
			- \psi ( m_1 - v_2 (\overline{y}_1)  ) \right] z_1   (C_r). 
		\end{align*}
		Analogously, there exists $(\overline{x}_2, \overline{y}_2) \in C_r$ such that
			$$- \int_{C_r} \left[  v_{12}(x)- v_{12}(y) \right]\, dz_2   (x,y) \\
			=  - \left[  
			\psi ( v_1 ( \overline{x}_2)  - M_2  )
			- \psi ( v_1 (\overline{y}_2)  - m_2  ) \right] z_2  (C_r).
		$$
		Therefore, 
		\beq\label{eqmono2}
		\begin{split}
			&\int_{C_r} \left[  (v_{12}(x)- v_{12}(y) \right]\, d (z_1 - z_2 ) (x,y) \\
			&= \left[  
			\psi ( M_1  - v_2 ( \overline{x}_1)  )
			- \psi ( m_1 - v_2 (\overline{y}_1)  ) \right]\! z_1 (C_r) 
			\!-\! \left[  
			\psi ( v_1 (\overline{x}_2)  - M_2  )
			- \psi ( v_1 (\overline{y}_2)  - m_2  ) \right] \!z_2  (C_r). 
		\end{split}
		\eeq
		
		We now consider two cases. 
		
		\vspace{.2cm}
		
		\noindent
		\textbf{Case 1:}  $\min\{ z_1 (C_r) , z_2 (C_r)\} < 1$ and $\max \{ z_1 (C_r) , z_2 (C_r)\} = 1$.
		We will only consider the case $z_1 (C_r) < 1$ and $z_2 (C_r) = 1$, 
		since the other case is analogous. 
		In this situation, $M_1 = m_1 = :C$ and 
		\begin{align*}
			&\int_{C_r} \left[  (v_{12}(x)- v_{12}(y) \right]\, d (z_1 - z_2 ) (x,y) \\
			&= \left[  
			\psi ( C  - v_2 ( \overline{x}_1)  )
			- \psi ( C - v_2 (\overline{y}_1)  ) \right] z_1 (C_r) 
			- \left[  
			\psi ( C  - M_2  )
			- \psi ( C  - m_2  ) \right] \\
			&= \left[  
			\psi ( C  - v_2 ( \overline{x}_1)  )
			- \psi ( C - v_2 (\overline{y}_1)  ) \right] z_1   (C_r) 
			+ \osc_{B_r} \psi (v_1 - v_2) \\
			&\geq ( 1  -  z_1 (C_r)  ) \osc_{B_r} \psi (v_1 - v_2) \geq 0.
		\end{align*}

		%
		%
		
		\vspace{.2cm}
		
		\noindent
		\textbf{Case 2:} $z_1(C_r) =  z_2(C_r) = 1$. 
		Therefore \eqref{eqmono2} becomes
		\begin{align*}
			&\int_{C_r} \left[  (v_{12}(x)- v_{12}(y) \right]\, d (z_1 - z_2 ) (x,y) \\
			&= \psi ( M_1  - v_2 ( \overline{x}_1)  ) - \psi ( v_1 ( \overline{x}_2)  - M_2  ) + \psi ( v_1 (\overline{y}_2)  - m_2  ) - \psi ( m_1 - v_2 (\overline{y}_1)  ) \geq 0,
		\end{align*}
		where the last inequality follows from the monotonicity of $\psi$, together with the fact that   
		$M_1  - v_2 ( \overline{x}_1) \geq v_1 ( \overline{x}_2)  - M_2$
		and $v_1 (\overline{y}_2)  - m_2 \geq m_1 - v_2 (\overline{y}_1)$.
	\end{proof}
	
	The following property hinges on  the geometric structure of the problem.
	\begin{lemma}\label{lm:geometricity}
		Let $(u_i, z_i) \in \Lip(\R^N)\times  L^\infty_{\mathbf{w}}\left(\R^N; \mathcal{M}^+ (C_r\big)\right)$ 
		satisfy (a) and (b) of Theorem~\ref{give me a name} for a function $g_i\in \Lip(\R^N)$, $i=1,2$. Let $K \subset \R^N$ be a compact set, and suppose that for some $\lambda\in \R$ we have
		$$
		\{ u_i \leq \lambda \} \subset K  \subset \{ g_1 = g_2 \}\qquad \text{ for } i= 1, 2.
		$$
		Then, 
		$$
		u_1 \land \lambda =  u_2 \land \lambda.
		$$
	\end{lemma}
	
	\begin{proof}
		For $i = 1, 2$ we have by property (b)
		$$
		\frac{h}{2 \, r} \int_{\R^N} \int_{C_r}(\vphi(w+x)-\vphi(w+y))\,  d (z_{i})_w  (x,y) dw=-\int_{\R^N}\vphi(u_i-g_i)\, dw, \quad \forall \, 
		\vphi\in C^0_c(\R^N).
		$$
		Let now $v_i =  u_i\land \lambda$ ($i = 1, 2$),  choose $\vphi = v_1 - v_2$ and note that its support is contained in $K$.  
		Thanks to Remarks~\ref{rm:monotonicity} and ~\ref{rm:monotonicity2}  and observing that $(v_1 - v_2) (u_1-u_2)\geq 0$, we obtain 
		\begin{align*}
			0 &\leq \frac{h}{2 \, r} \int_{\R^N} \int_{C_r} \left[  (v_1 - v_2)(w+x)- (v_1 - v_2)(w+y) \right]\, d (z_{1}- z_{2})_w (x,y) dw \\
			&=- \int_{\R^N} (v_1 - v_2) (u_1-u_2)\, dw \leq 0.
		\end{align*}
		Thus, 
		$$
		(v_1 - v_2) (u_1-u_2) = 0 \qquad \mathcal{L}^N\text{-a.e. on } \R^N
		$$
		and the conclusion follows.
	\end{proof}
	
	The next lemma deals with  the geometric problem \eqref{ATW0}.  We recall that the symbols $\subset_{\mathcal{L}^N}$
	and $=_{\mathcal{L}^N}$ are introduced in Section~\ref{sec:preliminaries}. 
	
	\begin{lemma}\label{lm:pbgeo}
		Let $g\in L^1_{loc}(\R^N)$ be such that $g \geq M>0$ in $\R^N\setminus B_R$ for some $M, R > 0$.
		Then, the problem 
		\beq\label{pbgeo}
		\min_{F \in \mathscr{M}(\R^N)} \left\{  \P (F) + \int_{F}g (w)  \, d w \right\}
		\eeq
		admits a minimal and a maximal solution  (up to $\mathcal{L}^N$-negligible sets),  denoted by $E_g^-$ and $E_g^+$, respectively, which are bounded.   Moreover,  if $g_1, g_2$ are as above and satisfy $g_1\leq g_2$, then 
		 $E^\pm_{g_2}\subset_{\mathcal{L}^N} E^\pm_{g_1}$.  If  in addition $g_1<g_2$ and $E_i$ is any solution of \eqref{pbgeo}, with $g_i$ in place of $g$ for $i=1, 2$, then  $E_2\subset_{\mathcal{L}^N} \! E_1$. 
	\end{lemma}
	\begin{remark} \label{rem inf g < 0}
		Note that  if $ \operatorname*{ess\,inf} g\geq 0$, then trivially  $E^+_g =_{\mathcal{L}^N} E^-_g =_{\mathcal{L}^N} \emptyset$. 
	\end{remark}
	\begin{proof}[Proof of Lemma~\ref{lm:pbgeo}] 
 To ease the notation, in this proof we identify Lebesgue measurable sets with their precise representative. 
 Most of the proof can be obtained by adapting the arguments of \cite[Proposition~6.1 and Lemma~6.2]{CMP15}, 
		see also \cite[Remark~6.3]{CMP15}.
		Here we only show the boundedness of solutions.
		In order to prove the boundedness of $E^+_g$,  assume first that $g$ is essentially bounded from below.  Assume
		by contradiction that there exists a sequence of points $\{x_n \} \subset  (\R^N\setminus B_R) \cap E^+_g$ such that $|x_n|\to +\infty$. By extracting a subsequence if needed we may assume $|x_i-x_j|>2R$ for $i\neq j$ and $|x_n|>3R$  for all $n \in \mathbb{N}$. Define
		\beq\label{gtilde}
		\tilde g:=\begin{cases}
			\operatorname*{ess\,inf} g & \text{in }B_{2R}\,,\\
			M & \text{in }\R^N\setminus B_{2R}\,,
		\end{cases}
		\eeq
		and note that, since $g \geq M$ out of $B_R$, we have $\tilde g\leq \tau_v g:=g(\cdot+v)$ for $|v|\leq R$. Hence, by the comparison principle and by translation invariance, 
		$$
		E^+_g+v = E^+_{\tau_v g} \subset E^+_{\tilde g}\qquad\text{for all $v$ s.t. }|v|\leq R\,. 
		$$
		In particular, the disjoint balls $\{B_R(x_n)\}_n$ are all contained in  $E^+_{\tilde g}$.  
		Therefore, 
		\begin{align*}
			&\P (E_{\tilde g}^+)+\int_{E_{\tilde g}^+}\tilde g\, dw\geq  (  \operatorname*{ess\,inf} g ) |B_{2R}|+\int_{\cup_nB_R( x_n)}\tilde g\, dw \\
			&= (  \operatorname*{ess\,inf} g ) |B_{2R}|+M\Big|\bigcup_nB_R( x_n)\Big|=+\infty\,,
		\end{align*}
		which is impossible. 
		
		We now remove the extra assumption. Define $g_n:=g\lor (-n)$. Then, by the above argument and the comparison principle, we get that  $E^+_{g_n}$ is an increasing sequence of bounded sets.
		In order to conclude, it will be enough to show that they are uniformly bounded as it would then follow that their $L^1$-limit $\cup_n E^+_{g_n}$ coincides with $E^+_g$. To show the equiboundedness we argue again by contradiction assuming that (up to passing to a non relabelled subsequence) there exists  a diverging sequence $\{x_n\}$ such that  $x_n \in E^+_{g_n}$,  $|x_n|> 3 R$  for all $n$, with  $|x_i-x_j|>2R$ if $i\neq j$.  Define $\tilde g_n$ as in \eqref{gtilde}, with $g_n$ in place of $g$, and note that by the comparison principle the sets $\{ E^+_{\tilde g_n} \}_n$ form an increasing sequence. 
		Then, arguing as before, if $m \geq n$ we have $B_R(x_n)\subset E^+_{\tilde g_n}\subset E^+_{\tilde g_m}$; i.e., for for every $m\in \N$, $\{B_{R}(x_i)\}_{i=1,\dots, m}$ is a family of disjoint balls contained in $E^+_{\tilde g_m}$. Thus,
		\begin{align*}
			0 &\geq \P (E_{\tilde g_m}^+)+\int_{E_{\tilde g_m}^+}\tilde g_m\, dw 
			\geq -\|g^-\|_{L^1(\R^N)}+\int_{\cup_{n=1}^mB_R( x_n)}\tilde g_m\, dw \\
			&=-\|g^-\|_{L^1(\R^N)}+M |\cup_{n=1}^mB_R( x_n)|\,,
		\end{align*}
		which leads to a contradiction if $m$ is large enough. 
	\end{proof}

	We are now  ready  to give the proof of Theorem~\ref{give me a name}.
	\begin{proof}[Proof of Theorem~\ref{give me a name}] We start by assuming $g\in \mathcal{L}^+$, while  the  case $g\in \mathcal{L}^-$ will be considered in Step 5. 
		For every $M > 0$, define $g^M:= g \wedge M$ and for every $t < M$ let $E_t^M$
		be the maximal solution of the problem
		\begin{equation} \label{probl M}
			\min_{F\in \mathscr{M}(\R^N)} \left\{ \P (F)+\frac1h \int_{F }(g^M-t)\, dw \right\},
		\end{equation}
		which exists, thanks to Lemma~\ref{lm:pbgeo}.
		For $t \geq M$ we  set $E_t^M = \R^N$.
		Note that, due to the comparison principle stated in 
		Lemma~\ref{lm:pbgeo} we conclude that for any $M_0>0$
		\beq\label{sublevelestimate}
		E_t^{M}\subseteq_{\mathcal{L}^N} E^{M_0+1}_t \subseteq_{\mathcal{L}^N} E^{M_0+1}_{M_0}  \quad \text{for all } M\geq M_0+1\text{ and }t\leq M_0\,,
		\eeq
		with $E^{M_0+1}_{M_0}$ a bounded set.
		We now show that there exists $t_0>0$ such that  
		\beq\label{claim}
	 	| E_{t_0}^{M} |  \neq 0 \quad\text{for all }M>0\,. 
		\eeq
		Since, by definition, $E_{t_0}^{M}=\R^N$ whenever $M \leq t_0$, we only consider the case $M > t_0$.
		To this aim, let $s_0>0$ be such that $F_0:=\{g\leq s_0\}$  has positive measure.  
		Then, by the minimality  of  $E^M_{t_0}$
		$$
		\P(E^M_{t_0})+\frac1h\int_{E^M_{t_0}}(g^M-t_0)\, dw\leq \P(F_0)+\frac1h\int_{F_0}(g^M-t_0)\, dw\leq 
		\P(F_0)-\frac1h (t_0-s_0)|F_0|<0\,,
		$$
		provided that $t_0$ is sufficiently large, thus showing that $E^M_{t_0}$  has positive measure. 
		

		We now define $u^M$ as
		\[
		u^M (x) := \inf \{ t \in \R : x \in  (E_t^M)^{(1)}  \}.
		\]
		Note that, by Remark~\ref{rem inf g < 0} and by construction,  we have 
		\beq\label{bounduM}
		\inf g\leq u^M\leq M\,.
		\eeq
		We split the remaining part of the proof into three steps.
		
		\noindent
		\textbf{Step 1.} (Minimality of $u^M$) We show that $u^M$ is the unique minimiser of \eqref{nlROF2}, with $g$ replaced by $g^M$.
		We start by observing that by  construction,  $\{u^M<t\}\subset_{\mathcal{L}^N} E_t^M \subset_{\mathcal{L}^N} \{u^M\leq t\}$  for all $t$ 
		and thus  $ E_t^M =_{\mathcal{L}^N} \{u^M\leq t\}$ 
		for all  $t \in \mathbb{R}$ such that $s \mapsto \mathcal{L}^N (\{ u^M \leq s  \})$ is continuous at $t$; i.e, 
		for all but countably many~$t$. 
		In turn, by approximation and using the lower semicontinuity of the perimeter it follows that 
		\begin{equation} \label{dichiariamo la cosa}
			\{u^M\leq t\} \text{ is a solution of \eqref{probl M} for every } t < M.
		\end{equation}
		
		By the Coarea Formula, a straightforward application of Fubini Theorem,  \eqref{bounduM}, and the minimality of $ E_t^M$, we have
		\[
		\begin{split}
			& J_r(u^M)+\frac1{2h}\int_{\R^N}\big[(u^M-g^M)^2-(g^M-M)^2\big]\, dw \\
			&= \int_{-\infty}^{M}\biggl(\P (E^M_t)+\frac1h\int_{E^M_t}(g^M-t)\, dw\biggr)dt\\
			&\leq  \int_{-\infty}^{M}\biggl(\P (\{v\leq t\})+\frac1h\int_{\{v\leq t\}}(g^M-t)\, dw\biggr)dt\\
			&\leq J_r(v)+  \frac1h \int_{-\infty}^{M} \int_{\{v\leq t\}}(g^M-t)\, dw \, dt\\
			&= J_r(v)+\frac1{2h}\int_{\R^N}\big[(v^M -g^M)^2-(g^M-M)^2\big]\, dw \\
			&\leq J_r(v)+\frac1{2h}\int_{\R^N}\big[(v -g^M)^2-(g^M-M)^2\big]\, dw\,,
		\end{split}
		\] 
		for all $v\in L^1_{loc}(\R^N)$, where we set $v^M: = v \land M$. Thus the claim follows. In turn,  by Proposition~\ref{existence minimum}, since $g^M \in \Lip (\R^N) \cap (\R + L^2 (\R^N))$   we also have $u^M \in \Lip (\R^N) \cap (\R + L^2 (\R^N))$ with
		$\textnormal{Lip} (u^M) \leq \textnormal{Lip} (g^M) \leq \textnormal{Lip} (g)$ and   there exists  $z^M\in L^\infty_{\mathbf{w}}\left(\R^N; \mathcal{M}^+(C_r)\right)$ such that (a) and (b) of Proposition~\ref{existence minimum} hold true, with $u_h^g$, $z$, and $g$ replaced by $u^M$, $z^M$, and $g^M$, respectively. Moreover,    $u^{M_1}\leq u^{M_2}$ if $M_1\leq M_2$.
		

		\smallskip
		
		\noindent
		\textbf{Step 2.} (Passing to the limit as $M \to \infty$) Note that by \eqref{sublevelestimate} (choosing $M_0  = t  = t_0$) and \eqref{claim} for every $M\geq 1+t_0$ there exists a point $w_M$ in the bounded set 
		 $(E^{1+t_0}_{t_0})^{(1)}$  such that $u^M(w_M)\leq t_0$.
		Recalling that $\textnormal{Lip} (u^M) \leq \textnormal{Lip} (g)$  and \eqref{bounduM}, we deduce that the family $\{u^M\}_{M>0}$ is locally uniformly bounded and equi-Lipschitz continuous. Recall also  that $u^{M_1}\leq u^{M_{2}}$ if $M_2>M_1$. Therefore, there exists 
		$u \in \Lip (\R^N)$ with $\textnormal{Lip} (u) \leq \textnormal{Lip} (g)$ such that 
		\begin{equation} \label{crescente}
			u^M \nearrow  u   \qquad \text{ as } M\to +\infty.
		\end{equation}
		Now, for every $M_0>0$ since the set $E^{M_0+1}_{M_0}$ is bounded, we have that there exists
		$\overline M=\overline M(M_0)\geq M_0+1$ such that $g^{M_1}=g^{M_2}$ in $E^{M_0+1}_{M_0}$ for $M_1$, $M_2\geq \overline M(M_0)$. In turn, by \eqref{sublevelestimate} and Lemma~\ref{lm:geometricity} we deduce that
		\beq\label{stableuM}
		u^{M_1}\land M_0=u^{M_2}\land M_0\quad\text{ for all $M_1$, $M_2\geq \overline M(M_0)$;}
		\eeq
		or, equivalently  
		\beq\label{stablesub}
		 E^{M_1}_t =_{\mathcal{L}^N} E^{M_2}_t  \quad\text{ for all $M_1$, $M_2\geq \overline M(M_0)$ and $t\leq M_0$.}
		\eeq
		Setting $E_t : = \{ u \leq t \}$, thanks to \eqref{crescente},
		\eqref{stableuM} and \eqref{stablesub} imply that
		$$
		 E_t=_{\mathcal{L}^N} E^M_t  \quad\text{ for all $M\geq \overline M(M_0)$ and $t\leq M_0$}\,,
		$$
		and that
		\beq\label{utile}
		u\land M_0=u^M\land M_0\quad\text{ for all $M\geq \overline M(M_0)$.}
		\eeq
		By the arbitrariness of $M_0$, the fact  that $u\in \mathscr{L}^+$ readily follows.
		Recall that  $\|z^M\|_{L^\infty_{\mathbf{w}}\left(\R^N; \mathcal{M} (C_r)\right)} = 1$  for every  $M > 0.$
		Therefore, by Banach-Alaoglu Theorem there exist a sequence $M_k \to \infty$
		and $z \in L^\infty_{\mathbf{w}}\left(\R^N; \mathcal{M}^+(C_r)\right)$ such that
		\beq\label{convzeta}
		z^{M_k} \stackrel{*}{\rightharpoonup} z \quad \text{ weakly star in } L^\infty_{\mathbf{w}}\left(\R^N; \mathcal{M} (C_r)\right),
		\eeq
		and 
		$$
		\| z \|_{L^\infty_{\mathbf{w}}\left(\R^N; \mathcal{M} (C_r)\right)} 
		\leq \liminf_{k \to \infty} \| z^{M_k} \|_{L^\infty_{\mathbf{w}}\left(\R^N; \mathcal{M} (C_r)\right)} \leq 1.
		$$
		Let now $\varphi \in C^0_c (\R^N)$ be fixed. Then, 
		\begin{equation*} 
			\begin{aligned}
				-   \int_{\R^N}\vphi(u^{M_k}-g^{M_k})\, dw 
				=  \frac{h}{2r}  \int_{\R^N} \int_{C_r}(\vphi(w+x)-\vphi(w+y))\, d z^{M_k}_w  (x,y) dw \qquad \text{ for all } k \in \mathbb{N}.
			\end{aligned}
		\end{equation*}
		Passing to the limit as $k \to \infty$, we obtain that 
		\[
		\frac{h}{2r} \int_{\R^N} \int_{C_r}(\vphi(w+x)-\vphi(w+y))\, d z_w(x,y) dw = -  \int_{\R^N}\vphi(u-g)\, dw,
		\]
		for every $\varphi \in C^0_c (\R^N)$. Hence, property (b) follows.  
		
		 We now recall that, thanks to Step~1, for all $k \in \mathbb{N}$ 
		properties (a) and (b) of Proposition~\ref{existence minimum} hold true, with $u_h^g$, $z$, and $g$ replaced by $u^{M_k}$, $z^{M_k}$, and $g^{M_k}$, respectively. 
		Let now $\psi\in C^0_c(\R^N)$. Then, for all $k \in \mathbb{N}$, 
		$$
		\int_{\R^N}\int_{C_r}\psi(w)(u^{M_k}(w+x)-u^{M_k}(w+y))\,d  z^{M_k}_w  (x,y)dw=
		\int_{\R^N}\psi  (w) \osc_{B_r(w)}u^{M_k}\, dw.
		$$
		Set now $K:=\textrm{Supp\, }\psi+\overline{B_r(0)}$  and $M_0:=\|u\|_{L^\infty(K)}$.
		Then, thanks to \eqref{utile}, for $k$ sufficiently large we have that $u^{M_k}=u$ on 
		$K$ and thus the above integral identity reads
		$$
		\int_{\R^N}\int_{C_r}\psi(w)(u(w+x)-u(w+y))\,d  z^{M_k}_w  (x,y) dw=
		\int_{\R^N}\psi\osc_{B_r(w)}u\, dw
		$$
		for all sufficiently large $k$. Using again \eqref{convzeta} and letting $k\to\infty$ we deduce
		$$
		\int_{\R^N}\psi(w)\int_{C_r}(u(w+x)-u(w+y))\,d z_w(x,y)dw=
		\int_{\R^N}\psi\osc_{B_r(w)}u\, dw
		$$
		for all $\psi\in C^0_c(\R^N)$, whence condition \eqref{(a)} for $u$ follows.

		\smallskip
		
		\noindent
		\textbf{Step 3.}  (Minimality for the geometric problem) We now show that $\{ u < t \}$ and $\{ u \leq t \}$ are the minimal and maximal solution, respectively, 
		of \eqref{ATW0}.  
		Let $t \in \R$. For $k$ sufficiently large such that $M_k > t$, thanks to \eqref{dichiariamo la cosa}
		and \eqref{utile} it follows that $\{ u \leq t \} = \{ u^{M_k} \leq t \}$ is a solution of \eqref{probl M} (with $M_k$ in place of $M$).
		Letting $k \to \infty$ we infer that $\{ u \leq t \}$ is a solution of \eqref{ATW0}.
		In order to show that $\{ u \leq t \}$ is the maximal solution, let $E$ be any solution of \eqref{ATW0}
		and let $t'>t$. Then, since $g-t'<g-t$ we infer from the minimality of $\{ u \leq t' \}$ and   Lemma~\ref{lm:pbgeo} that  $E\subset_{\mathcal{L}^N} \{ u \leq t' \}$.  
		But then
		$$
	 	E\subset_{\mathcal{L}^N} \bigcap_{t'>t} \{ u \leq t' \} = \{ u \leq t \}\,, 
		$$ 
		and thus $\{ u \leq t \}$ is the maximal solution.
		The argument to show that $\{ u < t \}$ is the minimal solution is analogous and we omit it.

		\noindent
		\textbf{Step 4.} (Uniqueness) We show that  there exists  a unique function 
		$u\in \mathscr{L}^+$ satisfying (a) and (b) of the statement. To this aim, let 
		$(u_i, z_i)\in \Lip(\R^N)\times  L^\infty_{\mathbf{w}}\left(\R^N; \mathcal{M}^+(C_r\big)\right)$, with
		$u_i\in \mathscr{L}^+$, $i=1,2$, satisfy (a) and (b). Then,  we may apply Lemma~\ref{lm:geometricity} to obtain $u_1 \land \lambda =  u_2 \land \lambda$ for all $\lambda>0$.

		\noindent
		\textbf{Step 5.} (The case $g\in \mathscr{L}^-$)
		If $g\in \mathscr{L}^-$, then clearly $\tilde g:=-g\in \mathscr{L}^+$. 
		Let now $(u^{\tilde g}_h, \tilde z)\in \mathscr{L}^+\times  L^\infty_{\mathbf{w}} \left(\R^N; \mathcal{M}^+(C_r)\right)$ satisfy the conclusions of Theorem~\ref{give me a name} with $g$ replaced by $\tilde g$.
		One can check that, if we set $u_h^g:=- u^{\tilde g}_h$ and we define $z$ such that
		$$
		\int_{C_r} f (x, y)\, d z_w (x,y):= \int_{C_r} f (y, x)\, d \tilde{z}_w (x, y), \quad \text{ for all } f \in C^0 (C_r),
		$$
		then the pair $(u_h^g, z)$ satisfies (a) and (b) of Theorem~\ref{give me a name} and $u_h^g$ is unique in 
		$\mathscr{L}^-$.
		Finally, it is also clear that the sets $A_t := \{ u_h^g > t \}$ and $E_t := \{ u_h^g \geq t \}$ are  the minimal and maximal solution, respectively, of the problem 
		$$
		\min_{F \in \mathscr{M}(\R^N)} \left\{  \P (F) + \frac1h\int_{F} (t-g (w)) \, d w \right\}.
		$$
		
		\noindent
		\textbf{Step 6.} (Comparison principle) We start by considering the case $g_1\in  \mathscr{L}^-$ and $g_2\in  \mathscr{L}^+$.

		For $i = 1, 2$  let $(u^{g_i}_h, z_i)$ be the pair given by the previous steps corresponding to $g_i$. Then, we have
		\beq\label{boh}
		\frac{h}{2 \, r} \int_{\R^N} \int_{C_r}(\vphi(w+x)-\vphi(w+y))\,  d (z_{i})_w  (x,y) dw=-\int_{\R^N}\vphi(u^{g_i}_h-g_i)\, dw, \quad \forall \, 
		\vphi\in C^0_c(\R^N).
		\eeq
		
		Since $u^{g_1}_h\in \mathscr{L}^-$ and $u^{g_2}_h\in \mathscr{L}^+$, the test function 
		$\varphi:=\big(u^{g_1}_h-u^{g_2}_h\big)_+$ has compact support. Plugging it  in \eqref{boh} and subtracting the two resulting equations and using that $g_1\leq g_2$, we get 
		\begin{align*}
			0 &\leq \frac{h}{2 \, r} \int_{\R^N} \int_{C_r} \left[  \big(u^{g_1}_h-u^{g_2}_h\big)_+(w+x)- \big(u^{g_1}_h-u^{g_2}_h\big)_+(w+y) \right]\, d (z_{1}- z_{2})_w (x,y) dw \\
			&\leq -\int_{\R^N}\big(u^{g_1}_h-u^{g_2}_h\big)_+ (u^{g_1}_h-u^{g_2}_h)\, dw\leq 0,
		\end{align*}
		where the first inequality follows from  Lemma~\ref{lm:monotonicity}, while the last one from the fact that $\big(u^{g_1}_h-u^{g_2}_h\big)_+ (u^{g_1}_h-u^{g_2}_h)\geq 0$. Thus $\big(u^{g_1}_h-u^{g_2}_h\big)_+ (u^{g_1}_h-u^{g_2}_h)\equiv 0$, which easily implies   $u^{g_1}_h\leq u^{g_2}_h$.
		
		If both $g_1$ and $g_2$ are in $\mathscr{L}^+$ ($\mathscr{L}^-$), we choose $\varphi:=\big(u^{g_1}_h\land \lambda-u^{g_2}_h\land \lambda\big)_+$ ($\varphi:=\big(u^{g_1}_h\lor (-\lambda)-u^{g_2}_h\lor (-\lambda)\big)_+$) for any $\lambda>0$ and we argue similarly to obtain that $u^{g_1}_h\land \lambda\leq u^{g_2}_h\land \lambda$ ($u^{g_1}_h\lor (-\lambda)\leq u^{g_2}_h\lor (-\lambda)$) for all $\lambda>0$. 
	\end{proof}
	We conclude this section with an estimate of the solution $u_h^g$ provided by Theorem~\ref{give me a name} when $g$ is the Euclidean norm.  This, in turn, will be instrumental in estimating the dicrete evolutions of balls, see Section~\ref{sec:balls} below.  To this aim, if  $E \subset \R^N$ is closed  and either $E$ or $E^c$ is bounded, let us denote by $T_h (E)$ the  maximal solution of the problem
	\begin{equation} \label{se vuoi}
		\min_{F \in \mathcal{M} (\R^N)} \left\{ \P (F) + \frac{1}{h} \int_{F} d_E (w) \, dw \right\},
	\end{equation}
	where $d_E$ denotes the signed distance function from $\pa E$, that is,
	$$
	d_E(w):=
	\begin{cases}
		\,\dist(w, \pa E) & \text{if }w\in E^c\,,\\
		-\dist(w, \pa E) & \text{if }w\in E\,.
	\end{cases}
	$$
	
	\begin{proposition} \label{prop no name}
		For every $h > 0$, let $\phi^h := u_h^g$, with $g (w) = |w|$. 
		Then, for every $\eta > 0$ there exist $c_0 = c_0 (r, \eta) > 0$ and $h_0 = h_0 (r, \eta) > 0$ 
		such that for every $h \leq h_0$
		\begin{equation} \label{questa finalmente va bene}
			\phi^h (w) \leq (|w| + c_0 h) \vee \eta \qquad \text{ for every } w \in \R^N.
		\end{equation}
		Moreover, for every $\alpha \in (0, 1/(N+1))$ there exists $h_1 = h_1 (\alpha) > 0$ such that
		\begin{equation} \label{anche questa finalmente va bene}
			\phi^h (0) \leq 2 h^{\alpha}
		\end{equation}
		for every $h \leq h_1$. 
		%
		%
		%
	\end{proposition}
		\begin{proof}
		Let us set
		\[
		v (w):=  (|w| + c_0 h) \vee \eta.
		\]
		We will show that there exist $z\in L^\infty_{\mathbf{w}} \left(\R^N; \mathcal{M}^+(C_r)\right)$,
		$c_0 = c_0 (r, \eta) > 0$, and $h_0=h_0(r, \eta)>0$ such that identity
		(a) of Theorem~\ref{give me a name} is satisfied with $u^g_h$ replaced by $v$ 
		and such that 
		\begin{equation} \label{inequality there}
			- h \, \mathfrak{Div}^r z + v \geq g \qquad \text{ in } \mathcal{M} (\R^N),
		\end{equation}
		for all $0<h\leq h_0$, where we recall that $g (w) = |w|$. 
			 
Note that, once \eqref{inequality there} is established, one can repeat the very same argument of Step~6 of the proof of Theorem~\ref{give me a name}, 
and this will imply \eqref{questa finalmente va bene}. 
		
		We start by defining $z$. 
		We consider two different cases.
		
		\vspace{.2cm}
		
		\noindent
		\textbf{Case 1:} $r < \eta/4$. For $h$ sufficiently small we can assume 
		$\eta - c_0 h \geq \frac{3}{4}\eta + r$.  
		We set
		\[
		z_w: = \delta_{r t (|w|) \widehat{w} } \times \delta_{- r t (|w|) \widehat{w} }, 
		\]
		where $\widehat{w} := w/|w|$ and $t : (0, \infty) \to [0, \infty)$, with the notation used in \eqref{def misure},
		and where
		\[
		t (\rho) : =
		\begin{cases}
			0 & \text{ if } \rho \leq \frac{\eta}{2}, \\
			\frac{4}{\eta} \rho  - 2  & \text{ if } \frac{\eta}{2} < \rho < \frac{3}{4}\eta,  \\
			1 & \text{ if } \rho \geq  \frac{3}{4}\eta. 
		\end{cases}
		\]
		A direct calculation shows that property (a) 
		of Theorem~\ref{give me a name} is satisfied with $u^g_h$ replaced by $v$.
		
		 Let us now show \eqref{inequality there}. 
		 Let $\vphi \in C^0_c (\R^N)$ with $\vphi \geq 0$.
		We have 
		\begin{align*}
			& \frac{h}{2r}  \int_{\R^N} \int_{C_r}(\vphi(w+x)-\vphi(w+y))\, d z_w (x,y) dw \\
			&=  \frac{h}{2r}  \int_{\R^N}  (\vphi(w + r t (|w|) \widehat{w})-\vphi(w - r t (|w|) \widehat{w}))\,  dw.
		\end{align*}
		Note that
		\begin{align*}
			& \frac{h}{2r}  \int_{\R^N}  \vphi(w + r t (|w|) \widehat{w}) \,  dw 
			=  \frac{h}{2r}  \int_{0}^{+\infty} \int_{\partial B_{\rho}}  
			\vphi (w + r t (\rho) \widehat{w} ) 
			\, d \mathcal{H}^{N-1} (w) \,  d \rho \\
			&=  \frac{h}{2r}  \int_{0}^{+\infty} 
			\left(\frac{\rho}{\rho + r t (\rho)}
			\right)^{N-1} 
			\int_{\partial B_{\rho + r t (\rho)}}  
			\vphi (w ) 
			\, d \mathcal{H}^{N-1} (w) \,  d \rho \\ 
			&=  \frac{h}{2r}  \int_{0}^{\eta/2} 
			\int_{\partial B_{s}}  
			\vphi (w ) 
			\, d \mathcal{H}^{N-1} (w) \,  d s \\
			&+  \frac{h}{2r} 
			\frac{1}{\left(1+ \frac{4r}{\eta} \right)^N}
			\int_{\eta/2}^{\frac{3}{4}\eta+ r} 
			\left(1+ \frac{2 r}{s}
			\right)^{N-1} 
			\int_{\partial B_{s}}  
			\vphi (w ) 
			\, d \mathcal{H}^{N-1} (w) \,  d s \\
			&+  \frac{h}{2r}  \int_{\frac{3}{4}\eta+ r}^{\infty} 
			\left(1- \frac{r}{s}
			\right)^{N-1} 
			\int_{\partial B_{s}}  
			\vphi (w ) 
			\, d \mathcal{H}^{N-1} (w) \, ds,
		\end{align*}
		where we made the change of variables $\rho + r t(\rho) = s$.
		An analogous calculation shows that
		\begin{align*}
			& \frac{h}{2r}  \int_{\R^N}  \vphi(w - r t (|w|) \widehat{w}) \,  dw \\
			&=  \frac{h}{2r}  \int_{0}^{\eta/2} 
			\int_{\partial B_{s}}  
			\vphi (w ) 
			\, d \mathcal{H}^{N-1} (w) \,  d s \\
			&+  \frac{h}{2r}  \frac{1}{\left(1- \frac{4r}{\eta} \right)^N} \int_{\eta/2}^{\frac34\eta - r} 
			\left(1- \frac{2 r}{s}
			\right)^{N-1} 
			\int_{\partial B_{s}}  
			\vphi (w ) 
			\, d \mathcal{H}^{N-1} (w) \,  d s \\
			&+  \frac{h}{2r}  \int_{\frac{3}{4}\eta - r}^{\infty} 
			\left(1+ \frac{r}{s}
			\right)^{N-1} 
			\int_{\partial B_{s}}  
			\vphi (w ) 
			\, d \mathcal{H}^{N-1} (w) \, ds.
		\end{align*}
		Therefore, 
		\begin{align*}
			& \frac{h}{2r}  \int_{\R^N} \int_{C_r}(\vphi(w+x)-\vphi(w+y))\, d z_w (x,y) dw \\
			&=  \frac{h}{2r}  \,  \int_{\eta/2 }^{\frac34 \eta - r}  
			\left[  \frac{\left(1+ \frac{2 r}{s}
				\right)^{N-1}}{\left(1+ \frac{4r}{\eta} \right)^N}
			-  \frac{\left(1- \frac{2 r}{s}
				\right)^{N-1}}{\left(1- \frac{4r}{\eta} \right)^N} \right]
			\int_{\partial B_{s}}  
			\vphi (w ) 
			\, d \mathcal{H}^{N-1} (w) \,  d s \\
			&+  \frac{h}{2r}  \int_{\frac34\eta+ r}^{\infty} 
			\left[ \left(1- \frac{r}{s}
			\right)^{N-1} - \left(1+ \frac{r}{s}
			\right)^{N-1} \right]
			\int_{\partial B_{s}}  
			\vphi (w ) 
			\, d \mathcal{H}^{N-1} (w) \, ds \\
			&+ \frac{h}{2r}  \int_{\frac34\eta - r}^{\frac34\eta + r} 
			\left[ 
			\frac{\left(1+ \frac{2r}{s}
				\right)^{N-1} }{\left(1+ \frac{4r}{\eta} \right)^N}
			- \left(1+ \frac{r}{s}
			\right)^{N-1}\right]
			\int_{\partial B_{s}}  
			\vphi (w ) 
			\, d \mathcal{H}^{N-1} (w) \, ds
		\end{align*}
		Now,
		\begin{align*}
			& \frac{h}{2r}  \,  \int_{\eta/2 }^{\frac34 \eta - r}  
			\left[  \frac{\left(1+ \frac{2 r}{s}
				\right)^{N-1}}{\left(1+ \frac{4r}{\eta} \right)^N}
			-  \frac{\left(1- \frac{2 r}{s}
				\right)^{N-1}}{\left(1- \frac{4r}{\eta} \right)^N} \right]
			\int_{\partial B_{s}}  
			\vphi (w ) 
			\, d \mathcal{H}^{N-1} (w) \,  d s \\
			&=  \frac{h}{2r}  \frac{1}{\left(1+ \frac{4r}{\eta} \right)^N} \,  \int_{\eta/2 }^{\frac34 \eta - r} 
			\left[  \left(1+ \frac{2 r}{s}
			\right)^{N-1} -  \left(1- \frac{2 r}{s}
			\right)^{N-1} \right]
			\int_{\partial B_{s}}  
			\vphi (w ) 
			\, d \mathcal{H}^{N-1} (w) \,  d s \\
			&-  \frac{h}{2r}  \,  \Bigg[    \frac{1}{\left(1- \frac{4r}{\eta} \right)^N} 
			- \frac{1}{\left(1+ \frac{4r}{\eta} \right)^N} \Bigg]
			\int_{\eta/2 }^{\frac34 \eta - r} 
			\left(1- \frac{2r}{s}
			\right)^{N-1} 
			\int_{\partial B_{s}}  
			\vphi (w ) 
			\, d \mathcal{H}^{N-1} (w) \,  d s \\
			&\geq - h \, L_1 \, \frac{4}{\eta} 
			\int_{\eta/2 }^{\frac34 \eta - r}  
			\left(1- \frac{2r}{s}
			\right)^{N-1} 
			\int_{\partial B_{s}}  
			\vphi (w ) 
			\, d \mathcal{H}^{N-1} (w) \,  d s \\
			&\geq - h \, L_1 \, \frac{4}{\eta} 
			\int_{\eta/2 }^{\frac34 \eta - r} 
			\int_{\partial B_{s}}  
			\vphi (w ) 
			\, d \mathcal{H}^{N-1} (w) \,  d s 
			\geq - h \, L_1 \, \frac{4}{\eta} 
			\int_{\R^N}   \vphi (w )  \,  d w, 
		\end{align*}
		where $L_1 = L_1 (r, \eta)$ is the Lipschitz constant of the function $t \mapsto 1/(1+t)^{N}$ in the interval 
		$[-\frac{4r}{\eta}, \frac{4r}{\eta}]$.
		We observe that
		\begin{align*}
			& \frac{h}{2r}  \int_{\frac34\eta+ r}^{\infty} 
			\left[ \left(1- \frac{r}{s}
			\right)^{N-1} - \left(1+ \frac{r}{s}
			\right)^{N-1} \right]
			\int_{\partial B_{s}}  
			\vphi (w ) 
			\, d \mathcal{H}^{N-1} (w) \, ds \\
			&\geq - h \,    \frac{L_2}{\eta} 
			\int_{\frac34\eta+ r}^{\infty} 
			\int_{\partial B_{s}}  
			\vphi (w ) 
			\, d \mathcal{H}^{N-1} (w) \, ds
			\geq 
			-  h \,   \frac{L_2}{\eta} 
			\int_{\R^N}  
			\vphi (w ) 
			\, d  w,
		\end{align*}
		where $L_2$ is the Lipschitz constant of the function $t \mapsto (1+t)^{N-1}$ in the interval $[-1, 1]$.
		We also have 
		\begin{align*}
			&\frac{h}{2r}  \int_{\frac34\eta - r}^{\frac34\eta + r} 
			\left[ 
			\frac{\left(1+ \frac{2r}{s}
				\right)^{N-1} }{\left(1+ \frac{4r}{\eta} \right)^N}
			- \left(1+ \frac{r}{s}
			\right)^{N-1}\right]
			\int_{\partial B_{s}}  
			\vphi (w ) 
			\, d \mathcal{H}^{N-1} (w) \, ds \\
			&\geq \frac{h}{2r}  \int_{\frac34\eta - r}^{\frac34\eta + r} \left(1+ \frac{r}{s}
			\right)^{N-1} \left[ 
			\frac{1}{\left(1+ \frac{4r}{\eta} \right)^N}
			- 1\right]
			\int_{\partial B_{s}}  
			\vphi (w ) 
			\, d \mathcal{H}^{N-1} (w) \, ds \\
			&\geq -  \frac{h}{2r} 2^{N-1} \int_{\frac34\eta - r}^{\frac34\eta + r}  \left[ 
			1 - \frac{1}{\left(1+ \frac{4r}{\eta} \right)^N} \right]
			\int_{\partial B_{s}}  
			\vphi (w ) 
			\, d \mathcal{H}^{N-1} (w) \, ds \\
			&\geq -  h \,   \frac{L_1 2^N}{\eta} \int_{\R^N}  
			\vphi (w ) 
			\, d  w.
		\end{align*}
		We have thus obtained that
		\begin{align*}
			& \frac{h}{2r}  \int_{\R^N} \int_{C_r}(\vphi(w+x)-\vphi(w+y))\, d z_w (x,y) dw 
			\geq  - h  \, \frac{\, L}{\eta} 
			\int_{\R^N}   \vphi (w )  \,  d w, 
		\end{align*}
		where $L =  (4+ 2^N) L_1 +  L_2$.
		Choosing $c_0 =  \frac{L}{\eta}$ and $h_0 = 1/c_0(\eta/4-r)$ we obtain \eqref{inequality there}.
		
		\vspace{.2cm}
		
		\noindent
		\textbf{Case 2:} $r \geq \eta/4$.
		For $h$ sufficiently small we can assume $\eta - c_0 h \geq \tfrac{3}{4} \eta$.
		We set
		\[
		z_w: =
		\begin{cases}
			\delta_{r\widehat{w} } \times \delta_{- r\widehat{w} } & \text{ if } |w| \geq \frac{\eta}{2}  + r, \\
			\delta_{r\widehat{w} } \times \frac{1}{|D (w)|} \mathcal{L}^N \res D (w) & \text{ if }   |w| < \frac{\eta}{2}  + r, 
		\end{cases}
		\]
		where we set \(D (w) :=  B_r (w) \cap B_{\eta- c_0 h}\) and $\widehat{w} = w / |w|$.
		A direct calculation shows that property (a) 
		of Theorem~\ref{give me a name} is satisfied with $u^g_h$ replaced by $v$.
		Let $\vphi \in C^0_c (\R^N)$ with $\vphi \geq 0$.
		We have 
		\begin{align*}
			& \frac{h}{2r}  \int_{\R^N} \int_{C_r}(\vphi(w+x)-\vphi(w+y))\, d z_w (x,y) dw \\
			&=  \frac{h}{2r}  \int_{\{ |w| \geq r + \frac{\eta}{2} \}}  (\vphi(w + r \widehat{w})-\vphi(w - r  \widehat{w}))\,  dw \\
			&+  \frac{h}{2r}  \int_{\{  |w| < r+ \frac{\eta}{2} \}} \Bigl( \vphi(w+r \widehat{w})  - \medint_{D (w)} 
			\vphi(w+y) \, d y\Bigr) \, dw.
		\end{align*}
		Arguing as above,
		\begin{align*}
			& \frac{h}{2r}  \int_{\{ |w| \geq r + \frac{\eta}{2} \}}  (\vphi(w + r \widehat{w})-\vphi(w - r  \widehat{w}))\,  dw \\
			&=  \frac{h}{2r}  \int_{2 r + \frac{\eta}{2}}^{+ \infty} \left(1-\frac{r}{s} \right)^{N-1} \int_{\partial B_{s}}  \vphi(w) \, d \mathcal{H}^{N-1} (w) \, d s \\
			&-  \frac{h}{2r}  \int_{\frac{\eta}{2}}^{+ \infty} \left(1+\frac{r}{s} \right)^{N-1} \int_{\partial B_{s}}  \vphi(w) \, d \mathcal{H}^{N-1} (w) \, d s \\
			&\geq -  \frac{h}{2r}  \int_{\frac{\eta}{2}}^{+ \infty} \left(1+\frac{r}{s} \right)^{N-1} \int_{\partial B_{s}}  \vphi(w) \, d \mathcal{H}^{N-1} (w) \, d s \\ 
			&\geq -  \frac{h}{2r}   \left(1+\frac{2 r}{\eta} \right)^{N-1} \int_{\R^N} \vphi(w) \, dw. 
		\end{align*}
		Moreover, 
		\begin{align*}
			& \frac{h}{2r}  \int_{\{  |w| < r+ \frac{\eta}{2} \}} \Bigl( \vphi(w+r \widehat{w})  - \medint_{D (w)} 
			\vphi(w+y) \, d y\Bigr) \, dw \\
			&\geq -  \frac{h}{2r}  \int_{\{  |w| < r+ \frac{\eta}{2} \}}
			\medint_{D (w)} 
			\vphi(w+y) \, d y  \, dw \\
			&\geq -  \frac{h}{2r}  \, C (r, \eta) \int_{\R^N}  \vphi(w)  \, dw,
		\end{align*}
		where we set 
		\[
		C (r, \eta)  = \left(| B_r (  p  ) \cap B_{3\eta/4 } | \right)^{-1} | B_{r+ \eta/2 }|,
		\]
		and  $p$  is any vector with $| p | = \tfrac\eta2 + r$.
		
		The conclusion follows by setting \( c_0 = c_0 (r, \eta) =  \frac{1}{2r}  C (r, \eta) +  \frac{1}{2r}  \left(1+\frac{2 r}{\eta} \right)^{N-1}\)
		and choosing $h_0$ such that $\eta - c_0 h_0 = \tfrac{3}{4} \eta$.

		We conclude the proof by showing \eqref{anche questa finalmente va bene}.
		
Let $\alpha$ be as in the statement. Recall (see \eqref{se vuoi}) that \( T_h (\overline{B}_{h^{\alpha}}) \) is the maximal solution of
\[
\min_{F \in \mathcal{M} (\R^N)} \left\{ \P (F) + \frac{1}{h} \int_{F} ( |w| - h^{\alpha}) \, dw \right\}.
\]
We first claim that 
\begin{equation} \label{intersection for h small}
T_h (\overline{B}_{h^{\alpha}}) \cap \overline{B}_{h^{\alpha}} \neq \emptyset \quad \text{for $h$ small.}
\end{equation}
Indeed, if \( T_h (\overline{B}_{h^{\alpha}}) \cap \overline{B}_{h^{\alpha}} = \emptyset \), then its energy is non-negative. However,
\[
\P (\overline{B}_{h^{\alpha}/2}) + \frac{1}{h} \int_{\overline{B}_{h^{\alpha}/2}} ( |w| - h^{\alpha}) \, dw 
< \P (\overline{B}_{h^{\alpha}/2}) - \frac{h^{\alpha-1}}{2} |\overline{B}_{h^{\alpha}/2} |,
\]
and the right-hand side diverges to $-\infty$ as $h \to 0^+$ by the assumption on $\alpha$. By minimality, \eqref{intersection for h small} follows.

By Theorem~\ref{give me a name}, \( T^h (\overline{B}_{h^{\alpha}}) = \{ \phi^h \leq h^{\alpha} \} \). 
Using \eqref{intersection for h small}, pick $\overline{x} \in \{ \phi^h \leq h^{\alpha} \} \cap \overline{B}_{h^{\alpha}}$. Since $\phi^h$ is $1$-Lipschitz,
\[
\phi^h (0) \leq \phi^h (\overline{x}) + |\overline{x}| \leq 2 h^{\alpha}.
\]
	\end{proof}

	\subsection{Weak formulation of the flow}\label{sec:weakformulation}
	In this section we introduce the weak formulation of the Minkowski mean curvature flow (formally) given by
	\beq\label{oee}
	V=-\mathcal{K}_r\,,
	\eeq
	with $\mathcal{K}_r$ denoting the nonlocal curvature associated with $\P$, see \eqref{our curvature}. Such a formulation extends to the present nonlocal setting the distributional formulation introduced in \cite{CMP17, CMNP19}.
	
	To this aim, given $z \in   L^\infty_{\mathbf{w}}    \left(\R^N  \times (0,T); \mathcal{M}^+(C_r)\right)$ we still denote by $\mathfrak{Div}^r z$
	its nonlocal spatial divergence, that is the measure  on $\R^N  \times (0,T)$ satisfying 
	\beq \label{def div spatial}
	\int_{\R^N  \times (0,T)} \varphi (w, t) \, d (\mathfrak{Div}^r z ) (w, t)
	=-\frac{1}{2r} \int_{0}^{T}\! \int_{\R^N}\!\int _{C_r}(\varphi(w+x,t)- \varphi(w+y,t))\, d z_{(w,t)} (x,y)\, dw dt
	\eeq
	for every $\varphi \in C^0_c (\R^N  \times (0,T))$, where $z_{(w,t)}$ denotes the evaluation at $(w, t)$ of the map $z$. Note that the above integrals are also well defined for 
	$\varphi\in L^1(0,T; C^0 (K))$ for any compact set $K\subset \R^N$.
	We recall that given a set $E$ and  $\delta>0$, we denote by  $(E)_\delta$ its $\delta$-neighbourhood $E+B_\delta$.  We also recall that a sequence of closed sets  $\{E_n\}$ of $\R^d$ converges to a closed set $E\subset \R^d$ in the Kuratowski sense, and we write $E_n \stackrel{\mathscr K}{\longrightarrow} E$, if the following two conditions hold: 
	\begin{enumerate}[label=(\roman*)]
		\item  if  $x_n\in E_n$ for all $n$, then any cluster point of $\{x_n\}_n$ belongs to $E$;
		\item any $x\in E$ is the limit of a sequence $\{x_n\}_n$, with $x_n\in E_n$ for all $n$.
	\end{enumerate}
	In the following, given a subset $E$ of $\R^N\times [0,+\infty)$, with $\text{Int}(E)$ we denote the interior of $E$ in the relative topology of $\R^N\times [0,+\infty)$.
	
	\begin{definition}\label{Defsol}
		Let $E^0\subset\R^N$ be 
		a closed set with compact boundary. 
		Let $E$ be a closed set in $\R^N\times [0,+\infty)$ and
		for each $t\geq 0$ denote $E(t):=\{x\in \R^N:\, (x,t)\in E\}$. 
		We say that $E$ is a {\em weak superflow} of the curvature flow \eqref{oee} with
		initial datum $E^0$ if
		\begin{itemize}
			\item[(a)] $E(0)\subseteq {E}^0$ and  either $E(t)$ is bounded for all $t\geq 0$,  or $\R^N\setminus E(t)$  is bounded for all $t\geq 0$;
			\item[(b)] 
			For every $t \geq 0$ and for every $\delta > 0$ there exists $\tau(\delta) > 0$ such that 
			$ E (s) \subset (E (t))_{\delta}$\footnote{We recall that $(E (t))_{\delta}$ denotes the $\delta$-neighborhood of $E(t)$} for every $s \in (t, t + \tau (\delta))$;
			\item[(c)] 
			$E(s)\stackrel{\mathscr K}{\longrightarrow} E(t)$ as $s\nearrow t$ for all $t>0$;

\item[(d)]   setting $ T^{\textnormal{ext}}:=\inf\{t>0:\, E(s)=\emptyset \text{ for $s\geq t$}\}$
			and $d(x,t):=\dist (x, E(t))$ for every $(x,t)\in \R^N\times (0,T^{\textnormal{ext}})$, 
			there exists $z\in L^\infty_{ \mathbf{w} } \left(\R^N  \times (0,T^{\textnormal{ext}}); \mathcal{M}^+(C_r)\right)$ 
			 such that the following holds:
			
\begin{itemize}
			
\vspace{.1cm}
	
\item[(d1)] we have 
\begin{equation} \label{eq:supersol}
\partial_t d \geq \mathfrak{Div}^r z \quad \text{ as measures in } (\R^N\times (0,T^{\textnormal{ext}}))\setminus E;
\end{equation}
	
\vspace{.1cm}
	
\item[(d2)] for  a.e.  $(w, t) \in \R^N\times (0,T^{\textnormal{ext}})$, we have  $z_{(w,t) } (C_r)\leq 1$ and
\beq \label{corriere}
\int_{C_r}(d (w+x, t)-d(w+y, t))\, d z_{(w,t)} (x,y) = \osc_{B_r(w)} d (\cdot, t);
\eeq 
 
 \vspace{.1cm}
	
\item[(d3)] for every $\lambda > 0$
 \beq\label{divpositivepart}
(\mathfrak{Div}^r z)^+ 
\in L^\infty(\{(w,t)\in\R^N\times (0,T^{\textnormal{ext}}):\, d(w,t)\geq \lambda\}).
\eeq
\end{itemize}

\end{itemize}
		
		Let $A$ be a relatively open set in $\R^N\times [0,+\infty)$. 
		We say that $A$ is a {\em weak subflow} of the curvature flow \eqref{oee} with
		initial datum $E^0$ if $\big( \R^N\times [0,+\infty) \big) \setminus A$ is a superflow with initial datum  $\R^N \setminus  \text{Int}_{\R^N}(E^0)$, 
		where   $\text{Int}_{\R^N}(E^0)$  denotes 
		the interior part of $E^0$  in the topology of $\R^N$.
		
		Finally, let $E^0 \subset \R^N$ be a closed set, 
		coinciding with the closure of its interior.
		Then, we say that a closed set $E \subset \R^N\times [0,+\infty)$
		is a {\em weak flow} with initial datum $E^0$ if $E = \overline{\text{Int}(E)}$, 
		$E$ is a weak superflow with initial datum $E^0$, and
		$\text{Int}(E)$ is a weak subflow
		with initial datum $E^0$.
	\end{definition}
	
	\begin{remark} \label{consequences of (b)}
		Note that, in particular, taking into account that $(\emptyset)_{\delta} = \emptyset$ for every $\delta > 0$, condition (b)  above implies the following: 
		
		\begin{itemize}
			
			\item[(b')]  if $t\ge 0$ and ${E}(t)=\emptyset$, then $E(s)=\emptyset$ for all $s > t$.
		\end{itemize} 
		
		We also note that condition (b) could likely be derived from (b') and the other properties, as done for instance for anisotropic mean curvature flows in \cite{CMP17, CMNP19}. However, since it follows directly from the minimizing movements construction used below (see Remark~\ref{per passare al limite nella formulazione debole}), we choose to include it in the definition rather than derive it from the other axioms, in order to simplify the presentation.

		
	\end{remark}
	We now define the corresponding weak level set flows.
	
	\begin{definition}[Level set subsolutions and supersolutions]\label{deflevelset1}
		Let  $u^0$ be a  uniformly continuous function on $\R^N$ such that $\{ u^0 \leq \lambda\}$ has compact boundary for all $\lambda\in \R$. We will say that a  lower semicontinuous   function $u:\R^N\times [0, +\infty)\to \R$ is a {\em  level set supersolution} corresponding to \eqref{oee}, with initial datum $u^0$,  if $u(\cdot, 0)\geq u^0 (\cdot)$ and if for a.e. $\lambda \in \R$ the closed sublevel set $\{(x,t)\,:\, u(x,t)\leq \lambda\}$  is a superflow {of}  \eqref{oee} in the sense of Definition~\ref{Defsol}, with initial datum $\{u^0\leq \lambda\}$.
		
		We will say that an  upper-semicontinuous   function $u:\R^N\times [0, +\infty)\to \R$ is a {\em  level set subsolution} corresponding to \eqref{oee}, with initial datum $u^0$, if~$-u$ is a level set supersolution  in the previous sense, with initial datum $-u^0$. 
		
		Finally, we will  say that a continuous   function $u:\R^N\times [0, +\infty)\to \R$ is a {\em  level set solution}  to \eqref{oee} if it is both a level set subsolution and  supersolution. 
	\end{definition}

	\subsection{Minimizing movements}\label{sec:atw}
	We aim at establishing an existence result for the weak flows introduced in  Definition~\ref{Defsol}, by adapting to the nonlocal setting a variant of the
	Almgren-Taylor-Wang minimizing movement scheme (\cite{ATW}) introduced in  \cite{Chambolle, CaCha}.

	Let $E^0\subset \R^N$ be closed with compact boundary. 
	Note that, in this case, we have either $\dd_{E^0} \in \mathscr{L}^+$ (if $E^0$ is bounded), 
	or  $\dd_{E^0} \in \mathscr{L}^-$ (if $\R^N \setminus E^0$ is bounded).
	Fix a time-step $h>0$ and set
	$E^0_h=E^0$. We then inductively define $E_h^{k+1}$ (for all $k\in \N \cup \{ 0 \}$) according to the following  procedure: 
	If $E_h^{k}\neq \emptyset$, $\R^N$, then let   $(u_h^{k+1},z_h^{k+1}) 
	\in \Lip(\R^N)\times   L^\infty_{\mathbf{w}}  \left(\R^N; \mathcal{M}^+(C_r)\right)$ satisfy
		\beq \label{eq:iterk}
	\begin{cases}
		\displaystyle \int_{C_r}(u_h^{k+1} (w+x)-u_h^{k+1}(w+y))\, d ( z_h^{k+1})_w (x,y) 
		= \osc_{B_r(w)}u_h^{k+1},\quad  \text{for a.e. }w; \\[10pt]
		(z_h^{k+1})_w(C_r)\leq 1, \quad  \text{for a.e. }w;\\[5pt]
		- h \, \mathfrak{Div}^r z_h^{k+1} + u_h^{k+1} = \dd_{E_h^k} \qquad \text{ in } \mathcal{M} (\R^N);
	\end{cases}
	\eeq	
	and  set  $E_h^{k+1}:=\{u_h^{k+1}\le 0\}$. If either $E_h^{k}=\emptyset$ or $E_h^{k}=\R^N$, then set $E_h^{k+1}:=E_h^{k}$. 
	We recall that the definition of nonlocal divergence $\mathfrak{Div}^r z$ 
	for $z \in L^\infty_{\mathbf{w}} \left(\R^N; \mathcal{M}^+ (C_r)\right)$ is given in \eqref{def div}.
	
	We denote by $T^{\,\textnormal{ext}}_h$ the \textit{ discrete  extinction time}, i.e. the first discrete time $hk$ such that $E_h^k=\emptyset$, if such a time exists; otherwise we set  $T^{\,\textnormal{ext}}_h =+\infty$.
	Analogously, we denote by $T^{\, \textnormal{all}}_h$ the 
	 discrete  extinction time of the complement of $E_h^k$, that is, 
	the first discrete time $hk$ such that $E_h^k=\R^N$, if such a time exists; otherwise we set  
	$T^{\, \textnormal{all}}_h=+\infty$.
	Finally, we define the \textit{discrete final time} $T^*_h$
	as $T^*_h: = \min \{ T^{\, \textnormal{ext}}_h , T^{\, \textnormal{all}}_h \}$.
	
	By Theorem~\ref{give me a name}, for every $k \in \N$ such that $h k < T^*_h$
	and if $\dd_{E_h^k} \in \mathscr{L}^{\pm}$,  there exists a unique $u_h^{k+1} \in \mathscr{L}^{\pm}$ 
	such that  \eqref{eq:iterk} is satisfied for some $z_h^{k+1} \in  L^\infty_{\mathbf{w}}  \left(\R^N; \mathcal{M}^+(C_r)\right)$.
	Thus, by induction, $E_h^{k}$ is closed for all $k \in \N$.
	Moreover, if $E^0$  (resp. $\R^N\setminus E^0$) is bounded, then $E_h^{k}$ (resp. $\R^N\setminus E_h^{k}$) is bounded for all $k \in \N$, and  therefore the construction above is well defined.
	
	\begin{remark} \label{robe di kappa}
		Note that, thanks to Theorem~\ref{give me a name}, $E_h^{k+1}$ is the maximal solution 
		of \eqref{se vuoi} with $E$ replaced by $E_h^{k}$. 
	\end{remark} 
	
	Note also that
	by Theorem~\ref{give me a name} it follows that $| \nabla u_h^{k+1} | \le  | \nabla \dd_{E_h^k} | = 1$~a.e.~in $\R^N$.
	Thus one deduces, in particular, that
	\begin{equation}\begin{array}{ll}\label{eq:ineqd}
			u_h^{k+1} \le d_{E_h^{k+1}} & \textup{ in } \{x \in \R^N : \, \dist(x, E_h^{k+1})>0\}\,,\\
			u_h^{k+1} \ge d_{E_h^{k+1}} & \textup{ in } \{x \in \R^N : \, \dist(x, E_h^{k+1})<0\}\,.
		\end{array}
	\end{equation}
	We are now in  a position to define the time discrete evolutions. 
	For every $h > 0$ we set
	\begin{equation}\label{discretevol}
		\begin{array}{l}
			E_h:=\{(w,t) \in \R^N \times [0, \infty) : w \in E_h^{[t/h]}\}, \vspace{2pt}\\
			E_h(t):=E_h^{[t/h]}=\{ w \in \R^N: (w,t)\in E_h\}, \quad t \geq 0, \vspace{2pt}\\
			d_h(\cdot,t):=\dd_{E_h(t)}, \quad
			d_h: \R^N \times [0, \infty) \to \R,\vspace{2pt}\\
			u_h(\cdot,t):=u_h^{[t/h]}, \quad u_h: \R^N \times [h, T^*_h) \to \R,\vspace{2pt} \\
			z_h(\cdot, t):= z_h^{[t/h]},  \quad z_h: \R^N \times [h, T^*_h) \to \mathcal{M}^+ (C_r), 
		\end{array}
	\end{equation}
	where $[\cdot]$ stands for the integer part of its argument.
	\begin{remark}[Discrete comparison and avoidance principle]\label{rm:dcp}
		The comparison principle stated in Theorem~\ref{give me a name}  implies that the  scheme is monotone, that is, the discrete evolutions satisfy the comparison principle. More precisely, if $E^0\subseteq F^0$ are closed sets and if we denote by $E_h$ and $F_h$ the  discrete evolutions with initial datum $E^0$ and $F^0$, respectively, then $E_h(t)\subseteq F_h(t)$ for all $t\geq 0$.  
		
		Moreover, if $\textnormal{dist} (E^0, F^0) = : \Delta > 0$, then $\textnormal{dist} (E_h(t), F_h(t)) \geq \Delta$
		for every $t \geq 0$.
		To see this, observe that the assumption implies that 
		$\dd_{E^0} \geq \Delta - \dd_{F^0}$. 
		Therefore, by the comparison principle stated in Theorem~\ref{give me a name} and thanks to  Remark~\ref{richiamare} we have that 
		$u^{\dd_{E^0}}_h \geq u_h^{\Delta -  \dd_{F^0}} = \Delta - u^{\dd_{F^0}}_h$.
		In turn, this inequality, together with \eqref{eq:ineqd} (with $k = 0$) implies that  
		$\textnormal{dist} (E_h^1, F_h^1) \geq \Delta$.
		The conclusion then follows by induction.
	\end{remark}

	\subsection{Comparison with the evolution of balls}\label{sec:balls} In this subsection, we exploit Remark~\ref{rm:dcp} to compare the discrete evolutions  \eqref{discretevol} with the minimizing movements  of the balls and derive an estimate, which will be useful in the convergence analysis.
	
	Let $\lambda >0$ and fix $\eta \in (0,\lambda)$. 
	Recalling the definition of $T_h(\overline{B}_\lambda)$ in \eqref{se vuoi} and Remark~\ref{richiamare}, we have that $T_h(\overline{B}_\lambda)=\{u_h^g \leq \lambda\}$, where $u_h^g$ is the solution provided by Theorem \ref{give me a name} with $g(w)=|w|$. Thanks to Proposition~\ref{prop no name}, there exist $c_0=c_0(r,\eta)>0$ and $h_0=h_0(r,\eta)>0$ such that for every $h \leq h_0$ we have 
	$$
	u_h^g(w) \leq (|w| + c_0 h) \vee \eta, \quad \text{ for every } w \in \R^N.
	$$
	Thus, for $h \leq h_0$ such that $c_0 h < \lambda$, the following inclusion holds
	\begin{equation}\label{1}
		T_h(\overline{B}_\lambda)=\{u_h^g \leq \lambda\} \supseteq \{|w|+c_0 h \leq \lambda\}=\overline{B}_{\lambda-c_0 h},
	\end{equation}
	where  we used the fact that $\eta < \lambda$.
	Iterating the argument of \eqref{1}, by Remark~\ref{robe di kappa} and Remark~\ref{rm:dcp} it follows that if $E^0= \overline{B}_{\lambda}$ one has 
	\beq \label{errehstima}
	E_h(s) \supseteq \overline{B}_{\lambda - c_0  \left[ \tfrac{s}{h} \right] h} \supseteq \overline{B}_{\lambda - c_0  \, s}, 
	\eeq
	for every $h \leq h_0$, as long as $0 < s < (\lambda - \eta)/ c_0)$.
	
	Now we return to  the motion from an arbitrary closed set $E^0$, with compact boundary.
	If for some $(w,t)\in\R^N\times [0, T_h^*)$ we have $d_h(w,t)> \lambda$, 
	then $\textnormal{dist} ( \overline{B}_\lambda (w) , E_h(t) ) > 0$.
	%
	Hence, by the comparison principle stated in Remark~\ref{rm:dcp} and by  \eqref{errehstima} 
	we have
	$$
	d_h(w,t+s)\geq \lambda - c_0 (s+h)
	$$
	for $s>0$ and $s+h < (\lambda - \eta)/ c_0$, and for $h \leq h_0 $. 
	

	By letting $\lambda \nearrow d_h(w,t)$ we obtain that for any $\eta \in (0, d_h(w,t))$
	there exist $c_0 = c_0 (\eta, r) > 0$ and $h_0 = h_0 (\eta, r) > 0$ such that 
	\beq\label{straponzina21}
	d_h(w,t+s) \geq d_h(w,t) - c_0 (s+h)
	\eeq
	for $s>0$ and $s+h < (d_h(w,t) - \eta)/ c_0$, and for $h \leq h_0 $. 
	By a similar argument, if $d_h(w,t)< 0 $ for some $(w,t)\in\R^N\times [0,T_h^*)$, 
	we obtain that for any $\eta \in (0, - d_h(w,t))$
	there exist $c_0 = c_0 (\eta, r) > 0$ and $h_0 = h_0 (\eta, r) > 0$ such that 
	\beq\label{straponzina2100}
	d_h(w,t+s) \leq d_h(w,t) + c_0 (s+h)
	\eeq
	for $s>0$ and $s+h < (-d_h(w,t) - \eta)/ c_0$, and for $h \leq h_0 $.

	\subsection{Convergence of the scheme} \label{subsec:minmov}
	In this section we show the (subsequential) convergence of the minimizing movement scheme introduced in  Section~\ref{sec:atw} 
	 to a weak superflow and subflow in the sense of Definition~\ref{Defsol}.
	We start by observing that there exists a subsequence $(h_l)_l$, 
	a closed set $E$ and open set $A\subset  E$,  such that
	\beq \label{E and A}
	\overline E_{h_l} \stackrel{\mathscr K}{\longrightarrow} E\qquad\text{and}\qquad 
	{(\text{Int}(E)_{h_l})}^c\stackrel{\mathscr K}{\longrightarrow} A^c,
	\eeq
	where closure, interior and complementation are meant relative to space-time $\R^N\times [0,+\infty)$.
	For every $t \geq 0$, we set 
	\[
	E(t) := \{ w \in \R^N :  (w, t) \in E \} \quad \text{ and } \quad 
	A(t) := \{ w \in \R^N :  (w, t) \in A \}.
	\]
	We have the following result.
	
	\begin{proposition}\label{prop:E}
		Possibly extracting a further subsequence (not relabeled) 
		from $(h_l)_{l}$, there exists a countable set $\mathcal N\subset (0, +\infty)$ 
		and a function $\td: \R^N \times ((0, \infty) \setminus \mathcal N)\to \overline{\R}$ such that 
		for all $t\in (0, +\infty) \setminus \mathcal{N}$ we have that
		${d_{h_l}}(\cdot, t) \to {\td}(\cdot, t)$ locally uniformly and
		\[
		{\td}(\cdot, t)^+ = \dist(\cdot, E(t)), \qquad {\td}(\cdot, t)^- = \dist(\cdot,A^c( t )).
		\]
		%
		Moreover, for every $w \in \R^N$ the functions $\dist(w,E(\cdot))$ and  $\dist(w,A^c(\cdot))$ are left~continuous
		and  right lower semicontinuous, respectively.
		Finally, $E(0)= E^0$ and  $A(0)= \textnormal{Int}_{\R^N}(E^0)$.
	\end{proposition}
	
	\begin{proof}
		We can argue exactly as in the the proof of \cite[Proposition~4.4]{CMP17}, using \eqref{straponzina21} and \eqref{straponzina2100} in place of \cite[Equations~(4.8)~and~(4.9)]{CMP17}.
	\end{proof}
	
	\begin{remark} \label{per passare al limite nella formulazione debole}
		Note that, using Proposition~\ref{prop:E} and passing to the limit in \eqref{straponzina21}
		and \eqref{straponzina2100} we obtain that for every 
		$(w, t) \in E^c$ with $ t \in (0, +\infty) \setminus \mathcal{N}$
		and for every $\eta \in (0, \td (w,t))$, there exists $c_0 = c_0 (\eta,r) > 0$ such that 
		\beq \label{ci piace}
		\td (w,t+s)\geq \td (w,t) - c_0 s
		\eeq
		for $0< s < (\td(w,t) - \eta)/ c_0$. 
		Analogously, for every 
		$(w, t) \in A$ with $ t \in (0, +\infty) \setminus \mathcal{N}$
		and for every $\eta \in (0, - \td (w,t))$ there exists $c_0 = c_0 (\eta,r) > 0$ such that  
		\beq \label{ci piace anche lei}
		\td (w,t+s)\leq \td (w,t) + c_0 s
		\eeq
		for $0< s < (- \td(w,t) - \eta)/c_0$. 
		In particular, \eqref{ci piace} implies that $E$ satisfies property (b) 
		of Definition~\ref{Defsol}.
		Condition \eqref{ci piace anche lei} implies that an analogous property holds for $A^c$.
	\end{remark}
	
	\begin{remark} \label{scemo chi legge}
		As already mentioned in Remark~\ref{consequences of (b)}, 
		if $E(t)=\emptyset$ for some $t\ge 0$,
		then 
		$E(s)=\emptyset$ for all $s\ge t$,
		so that we can define the extinction time $T^{\textnormal{ext}}$ of $E$.
		Analogously, from \eqref{ci piace anche lei} it follows that, 
		if $A (t) = \R^N$ for some $t \ge0$, then $A (s) = \R^N$ for all $s > t$,
		so we can define the extinction  time 
		${T}^{\textnormal{all}}$ of $A^c$. 
	\end{remark}
Motivated by the previous remark, in the following we set $T^*:= \min \{ {T}^{\textnormal{ext}}$, ${T}^{\textnormal{all}} \}$.
Note that
\[
T^* \leq \liminf_{l \to + \infty} T_{h_l}^*,
\]
where $(h_l)_l$ is the subsequence provided by Proposition~\ref{prop:E}.
\begin{theorem}\label{themthm}
	Let $E$ and $A$ be the sets given in \eqref{E and A}.
	Then, $E$ is a weak superflow in the sense of Definition~\ref{Defsol}
	with initial datum $E^0$, while $A$ is a weak subflow with initial datum ${E}^0$.
\end{theorem}
\begin{proof}
	Properties \textit{(a), (b)} and \textit{(c)} of Definition~\ref{Defsol} follow from Proposition~\ref{prop:E}
	and Remark~\ref{per passare al limite nella formulazione debole}, so we will only show \textit{(d)}.
	If $T^* < T^{\textnormal{ext}}$, then $E (t) = \R^N$ for every $t \ge T^*$.
	Therefore, if will be enough to check condition (d) in the time interval $(0, T^*)$. We adapt the arguments of \cite[Theorem~4.5]{CMP17}.
	
	Let $(h_l)_{l \in \N}$ be the subsequence given by Proposition~\ref{prop:E}. 
	Up to a further subsequence, and setting 
	$z_{h_l}(\cdot, t):=0$ for $t \in (0, T^*) \setminus [h_l, T_{h_l}^* )$,
	we have that $z_{h_l}$ converges  weakly-$*$ in $L^\infty_{ \mathbf{w} }  \left(\R^N  \times (0,T^*); \mathcal{M}^+(C_r)\right)$  to some $z$
	satisfying $\| z \|_{L^\infty_{ \mathbf{w} }  ( \R^N \times (0,T^*)  ; \mathcal{M}^+(C_r) )} \le 1$.
	Thanks to \eqref{eq:ineqd}, we have $u_h^{k+1} \le \dd_{E_h^{k+1}}$
	whenever $\dd_{E_h^{k+1}}\ge 0$ and thus, from~\eqref{eq:iterk}, 
	\begin{equation}\label{eq:ineqdisc}
		 h_l  \, \mathfrak{Div}^r z^{k+1}_{ h_l } \leq \dd_{E^{k+1}_{ h_l }} - \, \dd_{E^k_{ h_l }} 
		\qquad \text{ in } \mathcal{M} (\R^N \setminus E^{k+1}_{ h_l }),
	\end{equation}
	for every $ l , k \in \N$ such that $k  h_l  < T^*_{ h_l }$.
	
	Let $\varphi \in C_c^\infty((\R^N\times (0,T^*))\setminus E)$ be nonnegative.
	By the Kuratowski convergence of $( E_{h_l} )_{l}$ to $E$, for $l$ sufficiently large, 
	the distance of the support of $\varphi$ from $E_{h_l}$ is bounded away from zero. 
	In particular, there exist two positive constants $m$ and $M$ such that $m< d_{h_l} < M$
	on $\text{supp} \, \varphi$ for $l$ large enough. Thanks to~\eqref{eq:ineqdisc}, 
	\beq\label{oscineq}
	\begin{split}
		&\int_{0}^{T^*}\int_{\R^N} \biggl( \varphi(w,t) \frac{d_{h_l}(w,t+{h_l})-d_{h_l}(w,t)}{h_l} \\ 
		&\hspace{2cm}+ \frac{1}{2 \, r}  \int_{C_r}(\varphi(w+x, t)-\varphi(w+y, t))\, d (z_{h_l})_{(w, t + h_l)} (x,y)  
		\biggr) dw \, dt \\ 
		&= \int_{0}^{T^*}\int_{\R^N} \biggl( - \frac{\varphi(w,t)-\varphi(w,t-h_l)}{h_l} \, d_{h_l}(w,t) \\ 
		&\hspace{2cm}+ \frac{1}{2 \, r}  \int_{C_r}(\varphi(w+x, t)-\varphi(w+y, t))\, d (z_{h_l})_{(w, t + h_l)} (x,y) 
		\biggr) dw \, dt 
		\ge 0.
	\end{split}
	\eeq
	Using the convergence properties of $d_{h_l}$ to $\td$ stated in Proposition~\ref{prop:E} and the weak$^*$-convergence of $z_{h_l}$ to $z$, we can pass to the limit as $l\to\infty$  in the above expression to obtain 
	\begin{align*}
		&\int_{0}^{T^*}\int_{\R^N} \biggl( - \partial_t \varphi (w,t) \td(w,t)\, dwdt \\ 
		&\hspace{2cm}+ \frac{1}{2 \, r}  \int_{C_r}(\varphi(w+x, t)-\varphi(w+y, t))\, d z_{(w, t)} (x,y)  
		\biggr) dw \, dt 
		\ge 0,
	\end{align*}
	for every 
	$\varphi\in C_c^\infty((\R^N\times (0,T^*))\setminus E)$ with $\varphi \geq 0$. It follows
	that 
	\[
	\partial_t \td \geq \mathfrak{Div}^r z \quad \text{ as measures in } \mathcal{M} ((\R^N\times (0,T^*))\setminus E),
	\]
	where, we recall,  $\mathfrak{Div}^r z$ is the Radon measure defined in \eqref{def div spatial}.  This shows property (d1) of Definition~\ref{Defsol}. 
	
	We now estimate the measure $\mathfrak{Div}^r z_{h_l}$ from above, away from $E_{h_l}$. 
	We start by observing that
	\[
	\dd_{E^k_h} (w) \leq \min_{\tau \in E_h^k} | w - \tau | \quad \text{ for every } w \in \R^N, 
	\]
	so that, by  \eqref{eq:iterk} and the comparison principle in Theorem~\ref{give me a name}, 
	\[
	u_h^{k+1} (w) \le \min_{\tau \in E_h^k} \phi^h (w - \tau) \quad \text{ for every } w \in \R^N,
	\]
	where $\phi^h$ is as in Proposition~\ref{prop no name}.
	Thus,  if $\dd_{E^k_h}(w)\geq \lambda >0$,
	then, applying \eqref{questa finalmente va bene} with $\eta = \lambda$ we have that there exists 
	$c_0 = c_0 (\lambda, r)>0$ and $h_0 = h_0 (\lambda, r)> 0$ such that
	\beq \label{from above}
	u_h^{k+1}(w)\le \min_{\tau \in E_h^k} | w - \tau|+ c_0  h
	= \dd_{E^k_h}(w) + c_0  h , 
		\eeq
	provided $h \le h_0$.
	In turn, by~\eqref{eq:iterk}, we get
	\begin{equation}\label{numero}
		\mathfrak{Div}^r z_h^{k+1} \le c_0  \qquad\text{ in }\{\dd_{E^k_h} \geq \lambda \} 
		\quad \text{ as measures,} 
	\end{equation}
	and passing to the limit 
	\beq\label{bounddiv}
	\mathfrak{Div}^r z  \le c_0
	\qquad\text{in }\{(w,t)\in\R^N\times (0, T^*):\, \td(w,t)\geq \lambda \} 
	\eeq
	as measures.
	We have thus shown  that $(\mathfrak{Div}^r z)^+ 
	\in L^\infty(\{(w,t)\in\R^N\times (0,T^*):\, \td(w,t)\geq \lambda\})$ for every $\lambda>0$, that is, property (d3) of Definition~\ref{Defsol}. 
	
	We now note that if $\dd_{E^k_h}(w) >0$, then 
	$\dd_{E^k_h} (\cdot) \geq \dd_{E^k_h}(w) - |\cdot-w|$. Therefore, once again by  comparison principle 
	and thanks to \eqref{anche questa finalmente va bene}, we have that 
	$$
	u_h^{k+1}(w)\geq \dd_{E^k_h}(w) -\phi^h(0) \geq \dd_{E^k_h}(w)  - 2h^{1/(2 (N+1))} \,,
	$$
	for $h$ sufficiently small.
	The above inequality, together with \eqref{from above}, implies 
	that for all $t\in [h_l, T^*_{h_l})$ and any $\lambda>0$
	\beq \label{speed of conv}
	\|u_{h_l}(\cdot, t)-d_{h_l}(\cdot, t-h_l)\|_{L^{\infty}(\{w  :d_{h_l}(w,t-h_l)\geq \lambda \})}
	\leq 2 h_l^{1/(2 (N+1))},
	\eeq
	for $l$ sufficiently large.  An analogous estimate holds in the region where $d_{h_l}(w,t-h_l)\leq -\lambda$.
	Hence, recalling  the convergence properties
	of $E_{h_l}$ and $d_{h_l}$ (see also \cite[formula (4.13)]{CMP17}), we have that 
	for a.e. $t \in (0, T^*)$ 
	\begin{equation}\label{elleuno}
		u_{h_l} (\cdot, t) \to \td (\cdot, t) \qquad \text{ locally uniformly in } \R^N,
	\end{equation}
	with the sequence $(u_{h_l})_l$ locally uniformly bounded in space--time. 
	
	We are now going to show identity \eqref{corriere} with $\td$ in place of $d$.
	Let $\varphi \in C_c^\infty(\R^N\times (0,T^*))$. 
	From \eqref{eq:iterk} it follows that 
	\begin{align}
		& \int_{0}^{T^*} \int_{\R^N} \varphi (w, t)  \int_{C_r}  (u_{h_l} (w+x, t)-u_{h_l}  (w+y, t))\, d ( z_{h_l})_{(w, t)} (x,y) \, dw \, dt
		\nonumber  \\
		&= \int_{ h_l }^{T^*_{h_l}}  \int_{\R^N}  \varphi (w, t)  \osc_{B_r(w)}u_{h_l}  (\cdot, t) \, dw \, dt, \label{we will pass to the limit}
	\end{align}
	for every $l \in \N$, where we also used the fact that $z_{h_l}(\cdot, t):=0$ for $t \in (0, T^*) \setminus [h_l, T_{h_l}^* )$.
	
	Let us first consider the right hand side of the expression above.
	From \eqref{elleuno} it follows that for a.e. $t$
	$$
	\osc_{B_r(w)}u_{h_l}  (\cdot, t) \to \osc_{B_r(w)} \td  (\cdot, t) \qquad \text{ for every } w \in \R^N.
	$$
	Since the sequence $(u_{h_l})_l$ is locally (in space and time) uniformly bounded, this allows to apply the Lebesgue Dominated 
	Convergence Theorem, so that 
	\beq \label{RHS}
	\int_{ h_l }^{T^*_{h_l}}  \int_{\R^N}  \varphi (w, t)  \osc_{B_r(w)}u_{h_l}  (\cdot, t) \, dw \, dt
	\to \int_{0}^{T^*}\int_{\R^N}  \varphi (w, t)  \osc_{B_r(w)} \td (\cdot, t) \, dw \, dt \quad \text{ as } l \to + \infty.
	\eeq
	Concerning the left hand side of \eqref{we will pass to the limit}, we first observe that, 
	since $z_{h_l}$ converges to $z$ weakly-$*$ in $L^\infty_{ \mathbf{w} }  \left(\R^N  \times (0,T^*); \mathcal{M}^+(C_r)\right)$,
	we have  
	\begin{align}
		&\lim_{l \to \infty} \int_{0}^{T^*}\int_{\R^N} \varphi (w, t) \int_{C_r} 
		(\td (w+x, t)- \td  (w+y, t))\, d ( z_{h_l})_{(w, t)} (x,y)  \, dw \, dt \nonumber \\
		&= \int_{0}^{T^*}\int_{\R^N} \varphi (w, t) \int_{C_r} 
		(\td (w+x, t)- \td  (w+y, t))\, d z_{(w, t)} (x,y) \, dw \, dt. \label{LHS1}
	\end{align}
	Note also that
	\begin{align*}
		& \left| \int_{0}^{T^*}\int_{\R^N} \varphi (w, t)  \int_{C_r} 
		(u_{h_l} (w+x, t)- u_{h_l}  (w+y, t))\, d ( z_{h_l})_{(w, t)} (x,y) \, dw \, dt \right. \\
		&- \left. \int_{0}^{T^*}\int_{\R^N} \varphi (w, t) \int_{C_r} 
		(\td (w+x, t)- \td  (w+y, t))\, d ( z_{h_l})_{(w, t)} (x,y) \, dw \, dt \right| \\
		&\leq  \int_{0}^{T^*}\int_{\R^N} | \varphi (w, t)| \int_{C_r} 
		| u_{h_l} (w+x, t) - \td (w+x, t) | \, d ( z_{h_l})_{(w, t)} (x,y) \, dw \, dt  \\
		&+ \int_{0}^{T^*}\int_{\R^N} |\varphi (w, t)| \int_{C_r} 
		| u_{h_l} (w+y, t) - \td (w+y, t) | \, d ( z_{h_l})_{(w, t)} (x,y) \, dw \, dt,
	\end{align*}
	where the last two integrals vanish as $l \to \infty$ due to \eqref{elleuno} and an appeal to the Dominated Convergence Theorem. Taking into account \eqref{RHS} and \eqref{LHS1}, we obtain
	\begin{align*}
		& \int_{0}^{T^*}\int_{\R^N} \varphi (w, t) \int_{C_r} (\td (w+x, t)- \td  (w+y, t))\, d z_{(w, t)} (x,y) \, dw \, dt
		\nonumber  \\
		&= \int_{0}^{T^*}\int_{\R^N}  \varphi (w, t)  \osc_{B_r(w)} \td  (\cdot, t) \, dw \, dt.
	\end{align*}
	By the arbitrariness of $\varphi$, we have that identity \eqref{corriere} holds for the function $\td$. The fact that it holds also for $d=\td^+$ is now a consequence of Remark~\ref{rm:monotonicity}. This concludes the proof of condition (d2)
	and shows that $E$ is a weak superflow. The fact that $A$ is a weak subflow can be proved 
	in a similar way.
\end{proof}

\subsection{Comparison principle and well-posedness of the level set flow}
In this section we will establish a comparison principle between weak superflows and subflows, which will imply (by standard arguments) the generic uniqueness (up to fattening) of weak flows and the uniqueness of the corresponding  level set flow. 

We are now ready to state and prove  the the main result of this section.
\begin{theorem} \label{thm: wfcomparison} Let $E^0, F^0 \subset \R^N$ be two closed sets with compact boundary and let $E$, $F$ be, respectively, a superflow with initial datum $E^0$ and   a subflow with initial datum $F^0$, according to Definition~\ref{Defsol}.
	If  $\dist(E^0,\R^N\setminus F^0)~=:~\Delta>0$, then 
	$$
	\dist(E(t),\R^N\setminus F(t)) \geq \Delta \ \ \ \forall \ t \geq 0,
	$$
	with the convention that $\dist(A,\emptyset)=+\infty$ for any $A$.
\end{theorem}
\begin{proof} We adapt the proof of \cite[Theorem~3.3]{CMP17} to the present setting.   
	Let $T^{\,\textnormal{ext}}_E$ and $T^{\,\textnormal{ext}}_{F^c}$ be the extinction times for $E$ and $F^c:=\big(\R^N\times[0,\infty)\big)\setminus F$, respectively. For all $t>\min\left\lbrace T^{\,\textnormal{ext}}_E,T^{\,\textnormal{ext}}_{F^c} \right\rbrace =:T^*$, we have that either $E(t)$ or $\R^N\setminus F(t)$ is empty. For all such $t$'s the thesis is trivially true. Hence, we may assume, without loss of generality, that $T^* >0$ and we consider the case $t \leq T^*$. By continuation (using the left continuity of $d$) it is enough to show the conclusion of the theorem for a time interval $(0,t^*)$ for some $0<t^* \leq T^*$. We start by assuming that both $E(t)$ and $F(t)$ are bounded for all $t\geq 0$.
	
	Let us fix $0<\eta_1<\eta_2<\eta_3<\Delta/2$. Let $z_E$, $z_{F^c}\in L^\infty_{ \mathbf{w} }\left(\R^N  \times (0,T^*); \mathcal{M}^+(C_r)\right)$ the generalised Cahn-Hoffman fields  as in the definition of superflow, Definition~\ref{Defsol}-(d), corresponding to $E$ and $F^c$, respectively. We also denote
	$$
	d_E(x,t):=\dist(x, E(t))\quad\text{and}\quad d_{F^c}(x,t):=\dist(x, R^N\setminus F(t))
	$$
	and set 
	$$
	S:=\left\lbrace x \in \R^N: \ d_E(x,0)>\eta_1 \right\rbrace \cap  \left\lbrace x \in \R^N: \ d_{F^c}(x,0)>\eta_1 \right\rbrace.
	$$ 
	Note that $S$ is  a bounded open region  between the two boundaries $\pa E^0$ and $\pa F^0$. 
	
	Denote
	$$
	\td_E:=d_E \vee (\eta_2+Ct) \ \ \mbox{and} \ \ \td_{F^c}:=d_{F^c} \vee (\eta_2+Ct)
	$$
	where $C>0$ is a constant that will be chosen later.
	
	By our assumption we have $(\td_E+\td_{F^c})(\cdot,0) \geq \Delta$,  and  by construction 
	$$
	\td_E+\td_{F^c} \geq \Delta + (\eta_2-\eta_1) \ \ \ \mbox{in } (\R^N\setminus S) \times \left\lbrace 0 \right\rbrace.
	$$
	It follows from  Definition~\ref{Defsol}-(b) that there exists $t^* \in (0, T^*)$ such that
	\begin{equation}\label{w=0}
		\td_E+\td_{F^c} \geq \Delta+\frac{\eta_2-\eta_1}{2} \ \ \ \mbox{in } (\R^N\setminus S) \times (0,t^*),
	\end{equation}
	\begin{equation}\label{compaboh}
		S \subset \left\lbrace x \in \R^N: \ d_E(x,t)>\frac{\eta_1}{2} \right\rbrace \cap  \left\lbrace x \in \R^N: \ d_{F^c}(x,t)>\frac{\eta_1}{2} \right\rbrace\,\quad \text{for all }t\in (0, t^*),
	\end{equation}
	$$
	E(t)\subset \subset F(t) \ \ \ \mbox{for all } t \in (0,t^*)
	$$
	and 
	\begin{equation}\label{S''}
		\td_E=d_E \ \ \mbox{and} \ \ \td_{F^c}=d_{F^c} \ \ \ \mbox{in } S' \times (0,t^*),
	\end{equation}
	where 
	$$
	S':=\left\lbrace x \in \R^N: \ d_E(x,0)>\eta_3 \right\rbrace \cap  \left\lbrace x \in \R^N: \ d_{F^c}(x,0)>\eta_3 \right\rbrace.
	$$
	Since $d_E$ is Lipschitz continuous in space and $\partial_t d_E$ is a measure in $\{d_E>0\}$, it follows that $d_E$ (and thus $\td_E$) is of class $BV_{loc}(S \times (0,t^*))$, with time distributional derivative taking the form  
	$$
	\partial_t d_E= \sum_{t \in J} (d_E(\cdot,t+0)-d_E(\cdot,t-0)) d \mathcal{H}^N \res  \ (\R^N \times \left\lbrace t \right\rbrace)  +\partial_t^d d_E.
	$$
	Here $J$ denotes the (countable) set of times where $d_E$ jumps, while  $\partial_t^d $ is the diffuse part of $\partial_t d_E$. Clearly, the same holds for $\partial_t d_{F^c}$.
	It turns out that (arguing for instance as in  \cite[Lemma~2.4]{CMP17}) $d_E(\cdot,t+0)-d_E(\cdot,t-0)\geq 0$ for each $t \in J$. 
	Then, the jump part of $\partial_t d_E$ is a \textit{positive} Radon measure in $S \times (0,t^*)$.
	Since by \eqref{divpositivepart} the positive part of $\divr z_E$ is absolutely continuous, equation \eqref{eq:supersol} implies that 
	\beq\label{compa1}
	\partial_t^d d_E \geq \divr z_E\quad\text{as measures in $S \times (0,t^*)$.}
	\eeq
	An analogous inequality holds for $d_{F^c}$.
	By the chain rule for BV functions (see \cite{ADM} and \cite[Theorem~3.69]{AFP}), in $S \times (0,t^*)$ we have 
	\begin{equation}\label{compa2}
		\partial_t^d \td_E=
		\begin{cases}
			C & \text{a.e. in  }\left\lbrace (x,t): \eta_2+Ct > d_E(x,t) \right\rbrace \\
			\partial_t^d d_E & |\partial_t^d d_E|\text{-a.e. in }\left\lbrace (x,t): \eta_2+Ct \leq d_E(x,t) \right\rbrace.
		\end{cases}
	\end{equation}
	and similarly for $\partial_t^d \td_{F^c}$. Recalling again that $(\divr z_E)^+$ and $(\divr z_{F^c})^+$ belong to $L^\infty(S \times (0,t^*))$, using (\ref{compa1}) and (\ref{compa2}) (and the analogs for $d_{F^c}$) we deduce
	\begin{equation}\label{dimuniqflow}
		\partial_t^d \td_E \geq \divr z_E \ \ \mbox{and} \ \ \partial_t^d \td_{F^c} \geq  \divr z_{F^c} \quad\text{as measures in $S \times (0,t^*)$,}
	\end{equation}
	provided that we choose 
	$$
	C \geq \Vert \big(\divr z_E\big)^+ \Vert_{L^\infty(S \times (0,t^*))} + \Vert \big(\divr z_{F^c}\big)^+ \Vert_{L^\infty(S \times (0,t^*))}. 
	$$
	Note that the finiteness of the above $L^\infty$-norms is guaranteed by \eqref{divpositivepart} in Definition~\ref{Defsol}, together with \eqref{compaboh}.
	
	Let $\Psi \in C^1(\R)$ be a convex function such that $\Psi(t)=0$ for $t \leq 0$ and $\Psi(t)>0$ for $t>0$. Set $\Phi:=\Psi(\Delta-(\td_E+\td_{F^c}))$. By (\ref{w=0}) and the properties of $\Psi$ we have $\Phi(\cdot, t), \, 
	\Psi' \big( \Delta-(\td_E+\td_{F^c}) (\cdot, t) \big) \in C^0_c(S)$ for all  $t\in  (0,t^*)$. Thus, using as before the chain rule for BV functions, recalling that the jump parts  of $\partial_t \td_E$ and $\partial_t \td_{F^c}$ are positive Radon measures, taking into account \eqref{dimuniqflow} and the fact that $\Phi(\cdot, 0)=0$, we have for any $s\in (0, t^*)$
	\beq\label{compafinal}
	\begin{split} 
	 & \int_S \Phi(w,s) \, dw = \int_{\R^N} \Phi(w,s) \, dw   = \int_{\R^N\times (0,s)} d (\partial_t \Phi ) \\
	&
	\leq -\int_{\R^N\times (0,s)} \Psi'(\Delta-(\td_E+\td_{F^c})) \, d(\partial_t^d \td_E+\partial_t^d \td_{F^c}) \\
		& \leq - \int_{\R^N\times (0,s)} \Psi'(\Delta-(\td_E+\td_{F^c})) \, d(\divr (z_E+z_{F^c}))  \\
		& =- \frac{1}{2 r} \int_0^s \int_{\R^N} \int_{C_r}\Big(\Psi'((\Delta-\td_E)-\td_{F^c}))(w+x,t)\\		
		&\qquad\qquad\qquad\qquad-\Psi'((\Delta-\td_E)-\td_{F^c}))(w+y,t)\Big)d[(-z_E-z_{F^c})_{(w,t)}](x,y) dwdt\,.
	\end{split}
	\eeq
	Recall now that by \eqref{corriere}
	\[
	\int_{C_r}\big((\Delta -d_E (w+x, t))-(\Delta-d_E(w+y, t))\big)\, d ( - (z_E)_{ (w,t) }  ) (x,y) = \osc_{B_r(w)} (\Delta-d_E) (\cdot, t)\,,\ 
	\]
	and 
	\[
	\int_{C_r}\big(d_{F^c} (w+x, t))-d_{F^c}(w+y, t))\big)\, d (    ( z_{F^c} )_{(w,t)}  ) (x,y) = \osc_{B_r(w)} d_{F^c} (\cdot, t)\,,\ 
	\]
	for a.e. $(w, t) \in \R^N\times (0,T^*)$.  
	Thus, noticing also that 
	\beq\label{noticing}
	(\Delta-\td_E)= (-d_E+\Delta)\land (-\eta_2-C t+\Delta),
	\eeq
	by Lemma~\ref{lm:monotonicity} we get
	$$
	\int_{C_r} \big(\Psi'((\Delta-\td_E)-\td_{F^c}))(w+x,t)-\Psi'((\Delta-\td_E)-\td_{F^c}))(w+y,t)\big)d(-z_E-z_{F^c})(w,t)\geq 0
	$$
	for a.e. $(w, t) \in \R^N\times (0,T^*)$.	Thus, \eqref{compafinal}  yields
	$$
	\int_S \Phi(w,s) \, dw\leq 0
	$$ 
	for any $s\in(0, t^*)$, which by definition of $\Phi$ is equivalent to  $\td_E+\td_{F^c} \geq \Delta$ in $S \times (0,t^*)$. In particular, using \eqref{S''}, we obtain that $d_E+d_{F^c} \geq \Delta$ in $S' \times (0,t^*)$. We can now argue as at the end of \cite[Proof of Theorem~3.3]{CMP17} to infer from the latter inequality that $\dist(E(t),F^c(t)) \geq \Delta$ for every $t \in (0,t^*)$. 
	
	The case where both $E^0$ and $F^0$ are unbounded is completely analogous.  Finally, the case where $E^0$ is bounded, while $F^0$ is not, requires some additional observations. Indeed, $S$ and $S'$ are now unbounded. However, one can choose a ball $B$ sufficiently large that $d_{E}(\cdot, 0)$, $d_{F^c}(\cdot, 0)\geq 2\Delta$ in $B^c$ and enforce, by taking $t^*$ smaller if needed, $d_{E}(\cdot, t)$, $d_{F^c}(\cdot, t)\geq \Delta$ in $B^c$ for all $t\in (0, t^*)$. Thus the supports of $\Phi(\cdot, t), \, \Psi' \big( \Delta-(\td_E+\td_{F^c}) (\cdot, t) \big)$ are contained in $B$ for all $t\in (0, t^*)$ and all the above arguments apply also in this case. 
\end{proof}
An immediate corollary is the following comparison principle between level set sub- and super-solutions. 

\begin{theorem}\label{th:lscomp}
	Let $u^0$, $v^0$ be uniformly continuous functions on $\R^N$ such that $\{u_0\leq \lambda\}$, $\{v_0\leq \lambda\}$ have compact boundary for all $\lambda\in \R$.  Let $u$, $v$ be respectively a level set subsolution with initial datum $u^0$ and a level set supersolution with initial datum $v^0$, according to Definition~\ref{deflevelset1}. If $u^0\leq v^0$, then $u\leq v$.
\end{theorem}
\begin{proof}
	It follows easily from Theorem~\ref{thm: wfcomparison}. See \cite[Proof of Theorem~2.8]{CMNP19}.
\end{proof}

\medskip
At this point, having established existence and comparison principle for weak subflows and superflows, the existence and uniqueness for the level set flow follows by standard arguments. Let us summarise the results. We introduce the following {\em level set minimising movements}. 

Let $u^0:\R^N\to \R$ be a uniformly continuous function such that $\{u_0\leq \lambda\}$ has compact boundary for all $\lambda\in \R$.  Let $E^{0,\lambda}:= \{u^0 \le \lambda\}$ and let $E_{\lambda,h}$ be the corresponding time discrete evolutions, defined according with \eqref{discretevol} with $E^0$ replaced by $E^{0,\lambda}$. 

We introduce the level set discrete evolution $w_h:\R^N\times \R \to \R$ defined by
\begin{equation}\label{defls}
	w_h(x,t):= \inf\{\lambda\in\R\,: \, x\in E_{\lambda,h}(t)  \}.
\end{equation}
Arguing exactly as in \cite[Theorem~4.8]{CMNP19} one can finally prove the following main result concerning the well-posedness of the level set flow. 
\begin{theorem}\label{th:phiregularlevelset} 
	Let  $u^0$ be a uniformly continuous function on $\R^N$ such that $\{u_0\leq \lambda\}$ has compact boundary for all $\lambda\in \R$. Then the following holds: 
	
	{\rm (i) (Existence and uniqueness)} There exists a unique solution $u$ to the level set flow   with initial datum  $u^0$,  in the sense of Definition~\ref{deflevelset1}.
	
	{\rm (ii) (Approximation via minimizing movements)} The solution  $u$ is the locally uniform limit in $\R^N\times [0,+\infty)$, as $h\to 0^+$,  of the    level set minimizing movements  $w_h$ defined in \eqref{defls}.
	
	{\rm (iii) (Properties of the level set flow)}  For all but countably many $\lambda \in \R$, the fattening phenomenon does not occur;   that is,
		\beq \label{eq:nonfattening}
			\overline {\{(x,t)\,:\, u(x, t) < \lambda\}} = \{(x,t)\,:\, u(x, t) \le \lambda\}\,.
		\eeq

		Moreover,  for every  $\lambda$ such that \eqref{eq:nonfattening} holds,  the  sublevel set 
		$\{(x,t)\,:\, u(x, t) \le \lambda\}$ is the unique  solution
		to  \eqref{oee} in the sense of Definition~\ref{Defsol}, with initial datum $E^{0,\lambda}$, and 
		\beq\label{eq:nonfatteningbiz}
		E_{\lambda, h}\stackrel{\mathscr K}{\longrightarrow} \{(x,t)\,:\, u(x, t) \le \lambda\}\quad\text{and}\quad
		{ (\mathrm{Int\, }E_{\lambda, h})^c}\stackrel{\mathscr K}{\longrightarrow} \{(x,t)\,:\, u(x, t) \geq \lambda\}
		\eeq
		as $h\to 0^+$.
		Finally, for all $\lambda\in \R$ the sets $\{(x,t)\,:\, u(x, t) \le \lambda\}$ and $\{(x,t)\,:\, u(x, t) < \lambda\}$ are respectively the maximal superflow and minimal sublow with initial datum $E^{0,\lambda}$.
	\end{theorem}

	\subsection{Regularised  Minkowski type flows}\label{sec:variant}
	As in \cite{CMP12, CMP15} one could also consider  a  regularised version of the functional $J_r$ introduced in Section~\ref{sec:preliminaries}. To this aim,  note that  $J_r$ may be rewritten as
	\[
	J_r(E)= \frac{1}{2 r} \int_{\R^N}  \Chi{(-r,r)}(d_{E^{(1)}}(w))\, dw 
	\]
	where $E^{(1)}$ is defined by \eqref{def: A1}, and $d_{E^{(1)}}$ denotes  the signed distance function to $\pa E^{(1)}$. The idea is to smoothen up the step function 
	$ \frac{1}{2r} \Chi{(-r,r)}$; that is, to  replace it with 
	$f:\R\to [0,+\infty)$  of class $C^1$,  even,  nonincreasing in  $[0,+\infty)$ and  with compact support in $[-r, r]$. The associated regularised  nonlocal perimeter then reads
	\[
	\textnormal{Per}^f(E):=\int_{\R^N} f(d_{E^{(1)}}(w))\, dw,
	\]
	and notice that by the coarea formula (see \cite{CMP12} for the details) it may be rewritten also as 
	\[
	\textnormal{Per}^f(E)\ 
	= \int_0^{r} (-2s f'(s)) \textnormal{Per}_s(E) \,ds=  \int_0^{r} g(s) \textnormal{Per}_s(E) \,ds=\int_{\R^N}\int_0^{r} \frac{g(s)}{2 s} \osc_{B_s (w)} \Chi{E}\, ds \,dw
	\]
	where  $\textnormal{Per}_s$ is defined by \eqref{def: Per_r}, and where  we set \(g(s):= -2s f'(s)\). Note that $g$ is continuous, non-negative and with compact support. Hence, $\textnormal{Per}^f(\cdot)$ is a sort of weighted average of  $\textnormal{Per}_s(\cdot)$ over a bounded range of radii $s$.   The corresponding nonlocal total variation is given by 
	\[
	J^f(u)= \int_0^{r} g(s) J_s(u)\, ds =
	\int_{\R^N}\int_0^{r} \frac{g(s)}{2 s} \osc_{B_s (w)} u\, ds \,dw= 
	\int_{-\infty}^{\infty} \textnormal{Per}^f (\{ u \leq t \}) \, dt,
	\] 
	where the last inequality is due to the generalised coarea formula, see \cite[formula (2.9)]{CMP12}.	The associated nonlocal curvature is then formally (and rigorously for $C^2$ sets, see \cite[formulae (2.18) and (2.19)]{CMP12}) given by 
	\[
	\mathcal K^f (w, E)=\int_0^{r} g(s) \mathcal K_s(w,E)\, ds
	\]
	for $w\in \pa E$, where $K_s = \mathcal{K}_s^+ +\mathcal{K}_s^-$, see \eqref{our curvature}. One could then consider the generalised ROF problem  \eqref{nlROF2}, with datum $g\in \R+L^2(\R^N)$ and show the existence of a unique minimiser $u^g_h$ arguing as in Proposition~\ref{existence minimum}. In order  to derive the Euler-Lagrange equation, one can argue as for \eqref{**} to obtain that for all $\vphi\in C^0 (\R^N)\cap L^2(\R^N)$ with $J^f (\varphi) < \infty$ 
	\[
	\begin{aligned}
		- \frac{1}{h}  \int_{\R^N}\vphi(u_h^g-g)\, dw 
		\leq  J^f(u_h^g+\vphi)-  \J(u_h^g)  \leq  J^f(\vphi)  \\
		= \int_{\R^N}\int_0^r\|\Psi_{\vphi}(w, s; \cdot)\|_{C^{0}\left(C_1\right)}\, ds\,dw
		=\|\Psi_{\vphi}\|_{L^1\big(\R^N\times(0, r); C^{0}\left(C_1\right)\big)}\,, 
	\end{aligned}
	\]
	where we have set $\Psi_{\vphi}(w,s; x,y):=\frac{g(s)}{2s}(\vphi(w+sx)-\vphi(w+sy))$ and, as usual, $C_1$ stands for
	$\overline{B_1} \times \overline{B_1}$. Using Hahn-Banach Theorem to extend the linear continuous functional 
	\[
	L(\Psi_\vphi):= - \frac{1}{h}  \int_{\R^N}\vphi(u_h^g-g)\, dw
	\]
	to an element of $ \left(L^1\left(\R^N\times(0,r); C^{0}(C_1)\right)\right)'$  and  Riesz Representation Theorem \cite[Theorem~2.112]{FonsecaLeoniBook}, we can argue as in Step 3 of Proposition~\ref{existence minimum} to infer the existence of $z \in L^\infty_{\mathbf{w}} \left(\R^N\times(0,r); \mathcal{M}^+(C_1)\right)$ such that 
	\[
	z_{(w,s)} \big(C_1\big) \leq 1 \qquad\text{and}\qquad \int_{C_1}(u_h^g(w+sx)-u_h^g(w+sy))\, d z_{(w,s)}(x,y)= \osc_{B_s(w)}u_h^g
	\]
	for a.e. $(w, s)\in \R^N\times(0,r)$ and 
	\[
	h \int_{\R^N} \int_{0}^r\frac{g(s)}{2s}\int_{C_1}(\vphi(w+sx)-\vphi(w+sy))\, d z_{(w,s)}(x,y) \, ds\, dw = -  \int_{\R^N}\vphi(u_h^g-g)\, dw,
	\] 
	for all $\vphi\in C^0 (\R^N)\cap L^2(\R^N)$, with $J^f(\vphi)<+\infty$. Here $z_{(w,s)}$ denotes as usual the evaluation at 
	$(w,s)$ of the map $z$. This result suggests that the right {\em nonlocal divergence} operator $\mathfrak{Div}^f $, acting on any field  $z \in   L^\infty_{\mathbf{w}}    \left(\R^N \times(0,r) \times (0,T); \mathcal{M}^+(C_1)\right)$,
	is given by 
	\begin{multline*}
		\int_{\R^N  \times (0,T)} \varphi (w, t) \, d (\mathfrak{Div}^f z ) (w, t)\\
		:=-\frac{1}{2r} \int_{0}^{T}\! \int_{\R^N}\!\!\int_0^r\frac{g(s)}{2s}\int _{C_1}(\varphi(w+sx,t)- \varphi(w+sy,t))\, d z_{(w,s,t)} (x,y)\,ds\, dw\, dt
	\end{multline*}
	for every $\varphi \in C^0_c (\R^N  \times (0,T))$. With this definition at hand, one can define the notion of {\em weak superflow} for the geometric evolution $V=-\mathcal K^f$ as a  closed tube in space-time satisfying properties (a), (b), (c) of Definition~\ref{Defsol} and with  (d) replaced by 
	\begin{itemize}

\item[(d)]   setting $ T^{\textnormal{ext}}:=\inf\{t>0:\, E(s)=\emptyset \text{ for $s\geq t$}\}$
			and $d(x,t):=\dist (x, E(t))$ for every $(x,t)\in \R^N\times (0,T^{\textnormal{ext}})$, 
			there exists $z\in L^\infty_{ \mathbf{w} } \left(\R^N  \times (0,r)\times (0,T^{\textnormal{ext}}); \mathcal{M}^+(C_1)\right)$
			 such that the following holds:
				
\begin{itemize}
			
\vspace{.1cm}
	
\item[(d1)]  we have 
\begin{equation*} 
\partial_t d \geq \mathfrak{Div}^r z \quad \text{ as measures in } (\R^N\times (0,T^{\textnormal{ext}}))\setminus E;
\end{equation*}
	
\vspace{.1cm}
	
\item[(d2)] for  a.e. $(w,s, t) \in \R^N\times (0,r)\times (0,T^{\textnormal{ext}})$, we have  $z_{(w,s, t) } (C_1)\leq 1$ and
\[
				\int_{C_1}(d (w+sx, t)-d(w+sy, t))\, dz_{(w,s,t)}   (x,y) = \osc_{B_s(w)} d (\cdot, t);
		\]
 
 \vspace{.1cm}
	
\item[(d3)] for every $\lambda > 0$
 \[ 
(\mathfrak{Div}^r z)^+ 
\in L^\infty(\{(w,t)\in\R^N\times (0,T^{\textnormal{ext}}):\, d(w,t)\geq \lambda\}).
\]
\end{itemize}

	\end{itemize}
	The notions of weak subflow and weak flow are then modified accordingly.  All the results of Section~\ref{sec:2}  extend, mutatis mutandis,  to this weak formulation for $V=-\mathcal K^f$. We leave the details to the interested reader.

	\section{The fractional mean curvature flow}\label{sec:fractional}
	
	In this section we show how to apply the distributional approach  to the  {\em fractional mean curvature flow} (formally) given by 
	\beq\label{oeef}
	V=-\kappa_s\,,
	\eeq 
	where $\kappa_s$ stands for the fractional $s$-curvature; i.e., the first variation of \eqref{pesse}, defined (for sufficiently regular sets) as 
	\beq\label{scurvature}
	\kappa_s(x,E):= c_s\, P.V. \left( \int_{\R^N} \frac{\chi_{E^c}(y)-\chi_E(y)}{|x-y|^{N+s}} \, dy\right)\,, \qquad  x \in \partial E, 
	\eeq
	with $P.V.$ denoting the principal value and $c_s$ the positive constant defined in \eqref{pesse}, 
	see for instance \cite[Theorem~6.1]{FFMM}. Since the methods and the proofs are similar to those of the previous section, we will state the results and only indicate the main changes in the proofs.
	
	\subsection{The fractional ROF problem}
	
	Here we consider the problem
	\beq \label{nlROF2f}
	\min_{v\in L^{1}_{\text{loc}}(\R^N)}\biggl( \mathcal{J}_s (v)+ \frac{1}{2 h}\int_{\R^N}(v-g)^2\, dx\biggr),
	\eeq
	and the associated Euler-Lagrange equation, with $\mathcal{J}_s$ defined in \eqref{TVf}.
	We have 
	\begin{proposition} \label{existence minimum f}
		Let $h$ be a positive real number. Then, for every $g\in \R + L^2 (\R^N)$ there exists a unique solution $u_h^g$ to the problem \eqref{nlROF2f}, and 
		$u_h^g \in \R + L^2 (\R^N)$. In addition, 
		there exists $z\in L^\infty(\R^N \times \R^N)$ 
		such that 
		\begin{itemize}
			\vspace{5pt}
			\item[(a)] $\|z\|_{L^\infty(\R^N \times \R^N)} \leq 1$ and
			\beq\label{(a)0f}
			z(x,y)( u_h^g (x) - u_h^g (y) )=|  u_h^g (x) - u_h^g (y)  |\qquad\text{ for a.e. $(x,y)$;}
			\eeq
			\item[(b)]  for every $\varphi \in L^2(\R^N)$ with $\mathcal{J}_s(\varphi)<\infty$ 
			$$
			h \,  \frac{c_s}{2} \int_{\R^N} \int_{\R^N} z(x,y) \frac{\varphi(x)-\varphi(y)}{|x-y|^{N+s}} \, dxdy  = -  \int_{\R^N}\vphi(u_h^g-g)\, dx\,.
			$$ 
			\vspace{5pt}
		\end{itemize}
		Moreover, if $g_1, g_2 \in \R + L^2 (\R^N)$ and  $g_1 \leq g_2$ a.e. in $\R^N$, then $u^{g_1}_h \leq u^{g_2}_h$ a.e. in $\R^N$.
		Finally, if $g \in \Lip (\R^N) \cap ( \R + L^2 (\R^N) )$, then $u_h^g \in \Lip (\R^N)$ with 
		$\textnormal{Lip} (u_h^g) \leq \textnormal{Lip} (g)$.
	\end{proposition}
	\begin{proof}
		The existence of a unique solution $u^{g}_h$, which obeys the comparison principle, and the fact that $u^{g}_h$ inherits the Lipschitz continuity of $g$ (with $\textnormal{Lip} (u_h^g) \leq \textnormal{Lip} (g)$) can be proved exactly as in Proposition~\ref{existence minimum}. Also the proof of (a) and (b) is similar, however we sketch it for the reader's convenience.
		Arguing as for \eqref{**}, we arrive at 
		\begin{equation}\label{eulag}
			\mathcal{J}_s(\varphi) \geq \mathcal{J}_s(  u^g_h +\varphi)-\mathcal{J}_s(  u^g_h ) \geq -\frac{1}{h}\int_{\R^N} (  u^g_h -g)\varphi \, dx
		\end{equation}
		for all $\varphi \in L^2(\R^N)$ such that $\mathcal{J}_s(\varphi)<\infty$.
		Thus, defining the subspace of $L^1(\R^N \times \R^N)$
		$$
		\mathcal{M}:=\left\lbrace \frac{c_s}{2} \frac{\varphi(x)-\varphi(y)}{|x-y|^{N+s}} : \varphi \in L^2(\R^N), \ \mathcal{J}_s(\varphi)<\infty \right\rbrace 
		$$
		and the map $T:\mathcal{M} \to \R$ as
		$$
		T \left(  \frac{c_s}{2} \frac{\varphi(x)-\varphi(y)}{|x-y|^{N+s}} \right):=-\frac{1}{h}\int_{\R^N} (  u^g_h -g)\varphi \, dx,
		$$
		we have that the operator $T$ is continuous, with $\Vert T \Vert \leq 1$. Then, by Hahn-Banach Theorem, we may extend it to an element $\hat{T}$ of $(L^1(\R^N \times \R^N))'$ with the same norm. Finally, by Riesz Representation Theorem, there exists a function $z \in L^\infty(\R^N \times \R^N)$ with $\Vert z \Vert_{L^\infty(\R^N \times \R^N)}  = \Vert \hat{T} \Vert_{(L^1)'} \leq 1$ such that 
		$$
		-\frac{1}{h}\int_{\R^N} (  u^g_h -g)\varphi \, dx=   \frac{c_s}{2}  \int_{\R^N} \int_{\R^N} z(x,y) \frac{\varphi(x)-\varphi(y)}{|x-y|^{N+s}} \, dxdy
		$$
		and (b) is proved.
		Finally we take $\varphi=-(  u^g_h -M)$, where $M \in \R$ is such that $g - M \in L^2 (\R^N)$. 
		Then from (\ref{eulag}) we have
		\begin{align*}
			-\mathcal{J}_s(  u^g_h )=\mathcal{J}_s(M)-\mathcal{J}_s(  u^g_h )& \geq \frac{1}{h}\int_{\R^N} (  u^g_h -g)(  u^g_h -M) \, dx= \\
			&= -  \frac{c_s}{2}  \int_{\R^N} \int_{\R^N} z(x,y) \frac{  u^g_h (x)-  u^g_h (y)}{|x-y|^{N+s}} \, dxdy.
		\end{align*}
		Thus, recalling that $\|z\|_{L^\infty(\R^N \times \R^N)} \leq 1$, we obtain
		\begin{align*}
			\mathcal{J}_s(  u^g_h ) & \leq  \frac{c_s}{2}  \int_{\R^N} \int_{\R^N} z(x,y) \frac{  u^g_h (x)-  u^g_h (y)}{|x-y|^{N+s}} \, dxdy 
			 \leq  \frac{c_s}{2}  \int_{\R^N} \int_{\R^N} \frac{|  u^g_h (x)-  u^g_h (y)|}{|x-y|^{N+s}} \, dxdy=\mathcal{J}_s(  u^g_h ),
		\end{align*}	
		from which \eqref{(a)0f} follows.
	\end{proof}
	Motivated by the above proposition and following  \cite{GO} we introduce the following nonlocal operators:
	Given a function $\vphi : \R^N \to \R$, we define the nonlocal gradient $\nabla_s \vphi : (\R^N \times \R^N)\setminus D  \to \R$ as 
	$$
	(\nabla_s \vphi) (x,y):= \frac{c_s}{2} \frac{\vphi(x)-\vphi(y)}{|x-y|^{N+s}}
	$$
	where $D:=\left\lbrace x=y\right\rbrace $. 
	Given a {\em nonlocal vector} $w \in L^\infty( \R^N \times \R^N)$, we define its nonlocal divergence $\divs w$ as the operator 
	\begin{equation} \label{divr0f}
		\langle  \divs w, \vphi \rangle:= -\langle \nabla_s \vphi,w \rangle=- \frac{c_s}{2}   \int_{\R^N} \int_{\R^N} w(x,y) \frac{\varphi(x)-\varphi(y)}{|x-y|^{N+s}} \, dxdy 
	\end{equation}
	for every $\vphi$ such that $\mathcal{J}_s(\vphi)<+\infty$. In particular,  $\divs w$ defines a distribution.  Thus, condition (b) of Proposition~\ref{existence minimum f} is equivalent to 
	$$
	-h\,\divs z+u^g_h=g
	$$
	as distributions in $\R^N$.
	We can now state the main result of the section, which is the analogue of Theorem~\ref{give me a name}.

	\begin{theorem} \label{give me a name f}
		Let $g \in \mathscr{L}^+$ ($g \in \mathscr{L}^-$) and let $h$ be a positive constant. 
		Then, there exist $z\in L^\infty\left(\R^N\times \R^N \right)$ and a unique Lipschitz function $u_h^g$
		in $\mathscr{L}^+$ ($\mathscr{L}^-$),  such that 
		\begin{itemize}
			\vspace{5pt}
			\item[(a)] $\|z\|_{L^\infty\left(\R^N\times \R^N \right) }\leq 1$  and
			\beq\label{(a)f}
						z(x,y)( u_h^g (x) - u_h^g (y) )=|  u_h^g (x) - u_h^g (y)  |\qquad\text{ for a.e. $(x,y)$;}
			\eeq
			
			\item[(b)]  $\displaystyle-h\,\divs z+u^g_h=g \qquad \text{ in } \mathcal{D}' (\R^N).$
			\vspace{5pt}

		\end{itemize}
		Moreover, $\textnormal{Lip} (u_h^g) \leq \textnormal{Lip} (g)$. 		
		In addition, for every $t \in \R$ the sets  $A_t := \{ u_h^g < t \}$ ($A_t := \{ u_h^g > t \}$) and $E_t := \{ u_h^g \leq t \}$ ($E_t := \{ u_h^g \geq t \}$) are  the minimal and maximal solution, respectively, of problem 
		\eqref{ATW0} with $ \P$ substituted by $\Ps$.
Finally, the following comparison principle holds: if $g_1, g_2 \in \mathscr{L}^+\cup \mathscr{L}^-$ and  $g_1 \leq g_2$, then $u^{g_1}_h \leq u^{g_2}_h$.
	\end{theorem}
	
	\begin{remark} \label{test f}
	If $w \in L^\infty( \R^N \times \R^N)$, 
	from \eqref{divr0f} it follows that $\langle  \divs w, \vphi \rangle$ is well defined also for functions $\vphi \in \Lip_c(\R^N)$. Then, identity (b) can be tested also with functions $\vphi \in \Lip_c(\R^N)$.
	\end{remark}

		As in Section~\ref{sec:ROF} we  state  some preliminary results and observations.
	\begin{remark}\label{rm:monotonicity2f}
	For $i=1,2$, let $z_i\in   [-1,1]$ and $b_i, c_i \in \R$ be such that
		$z_i (b_i - c_i )= |b_i - c_i|$. 
		Then, 
				\[
				(z_1-z_2)\big( (b_1 - b_2) - (c_1 - c_2) \big) \geq 0\,.
		\]
	\end{remark}
	In fact, a stronger monotonicity property holds, which is the analogue of Lemma~\ref{lm:monotonicity}, 
	whose elementary proof is left to the reader.
	\begin{lemma}\label{lm:monotonicityf}
		For $i=1,2$, let $z_i\in   [-1,1]$ and $b_i, c_i \in \R$ be such that
		$z_i (b_i - c_i )= |b_i - c_i|$. 
		Let $\psi, g_1, g_2: \ \R \to \R $ be monotone nondecreasing.
		Then 
		$$
		(z_1-z_2) \Big( \psi \big( g_1 (b_1) - g_2 (b_2) \big) - \psi \big( g_1 (c_1) - g_2 (c_2) \big) 	
		\Big) \geq 0.
		$$	
	\end{lemma}
	We can now prove:
	\begin{lemma}\label{lm:geometricityf}
		Let $(u_i, z_i) \in \Lip(\R^N)\times  L^\infty\left(\R^N\times\R^N\right)$ 
		satisfy (a) and (b) of Theorem~\ref{give me a name f} for a function $g_i\in \Lip(\R^N)$, $i=1,2$. Let $K \subset \R^N$ be a compact set, and suppose that for some $\lambda\in \R$ we have
		$$
		\{ u_i \leq \lambda \}  \subset K \subset \{ g_1 = g_2 \}  \qquad \text{ for } i= 1, 2.
		$$
		Then, 
		$$
		  u_1 \land \lambda  =  u_2 \land \lambda. 
		$$
	\end{lemma}
	\begin{proof}
		The proof is similar to that of Lemma~\ref{lm:geometricity}. For $i = 1, 2$, thanks to Remark~\ref{test f}, we have 
		$$
		h \, \frac{c_s}{2} \int_{\R^N} \int_{\R^N} z_i(x,y) \frac{\varphi(x)-\varphi(y)}{|x-y|^{N+s}} \, dxdy =- \int_{\R^N}\vphi(u_i-g_i)\, dw, \quad \forall \, 
		\vphi\in \Lip_c(\R^N).
		$$
		Let now $v_i =   u_i \land \lambda$ ($i = 1, 2$), and choose $\vphi = v_1 - v_2$. 
		Thanks to Lemma~\ref{lm:monotonicityf} we obtain that
		\begin{align*}
			0 &\leq h \, \frac{c_s}{2} \int_{\R^N} \int_{\R^N} (z_1(x,y)-z_2(x,y)) \frac{(v_1-v_2)(x)-(v_1-v_2)(y)}{|x-y|^{N+s}} \, dxdy \\
			&=- \int_{\R^N} (v_1 - v_2) (u_1-u_2)\, dw \leq 0\,.
		\end{align*}
		We may now conclude as in the proof of Lemma~\ref{lm:geometricity}.
	\end{proof}
	We  also have: 
	\begin{remark}\label{lm:pbgeof}
		One can readily check that the statement of Lemma~\ref{lm:pbgeo}  holds true also with $\P$ replaced by $\Ps$, with  the same proof.
	\end{remark}

	We are ready to prove Theorem~\ref{give me a name f}.
	\begin{proof}[Proof of Theorem~\ref{give me a name f}]
		We define $g^M$, $u^M$, $z^M$ as in the proof of Theorem~\ref{give me a name} (with $\P$ replaced by $\Ps$) and we proceed in the same way (using Proposition~\ref{existence minimum f},  Lemma~\ref{lm:geometricityf} and Remark~\ref{lm:pbgeof} in place of  Proposition~\ref{existence minimum},  Lemmas~\ref{lm:geometricity} and \ref{lm:pbgeo}, respectively) to establish all the properties \eqref{sublevelestimate}--\eqref{utile}. We now show how to pass to the limit (as $M\to \infty$) in the equations satisfied by $(u^M, z^M)$. 
		Recall now that $\|z^M\|_{L^\infty\left(\R^N\times\R^N\right)}  \leq 1$. Thus, there exists a sequence $M_k \to \infty$ and $z \in L^\infty(\R^N \times \R^N)$ such that
		\begin{equation}\label{T6}
			z^{M_k} \stackrel{*}{\rightharpoonup} z \ \ \mbox{weakly--$*$ in} \  L^\infty(\R^N \times \R^N),
		\end{equation}
		with $\|z\|_{L^\infty\left(\R^N\times\R^N\right)} \leq 1$.
		Let now $\varphi \in C^{\infty}_c(\R^N)$ be fixed and recall that 
		$$
		-\int_{\R^N} \varphi(u^{M_k}-g^{M_k}) \, dx=h \, \frac{c_s}{2} \int_{\R^N} \int_{\R^N} z^{M_k}(x,y) \frac{\varphi(x)-\varphi(y)}{|x-y|^{N+s}} \, dxdy \ \ \ \forall k \in \N.
		$$
		Recall that, by $\eqref{crescente}$, there exists a Lipschitz function $u$ such that
		$u^{M_k} \nearrow  u$. Then, passing to the limit  as $k \to \infty$ and using \eqref{T6}, we get
		$$
		-\int_{\R^N} \varphi(u-g) \, dx=h \, \frac{c_s}{2} \int_{\R^N} \int_{\R^N} z(x,y) \frac{\varphi(x)-\varphi(y)}{|x-y|^{N+s}} \, dxdy 
		$$
		for all $\varphi \in C^{\infty}_c(\R^N)$ . This establishes property (b) of the statement.
		
		Let now $\psi\in C^{\infty}_c(\R^N\times\R^N)$ and recall that $z^{M_k}(x,y)(u^{M_k}(x)-u^{M_k}(y))=|u^{M_k}(x)-u^{M_k}(y)|$ for a.e. $(x,y)\in \R^N\times\R^N$. Let $K\subset \R^N$ be a compact set such that $\textrm{Supp\,}\psi\subset K\times K$ and set $M_0:=\|u\|_{L^\infty(K)}$. Then, thanks to \eqref{utile}, for $k$ sufficiently large we have that $u^{M_k}=u$ on 
		$K$ and thus
		$$
		\int_{\R^N}\int_{\R^N}z^{M_k}(x,y)\psi(x,y)(u(x)-u(y)) dxdy=
		\int_{\R^N}\int_{\R^N}\psi(x,y)|u(x)-u(y)| dxdy
		$$
		for all sufficiently large $k$. Passing to the limit and using again \eqref{T6} we obtain 
		$$
		\int_{\R^N}\int_{\R^N}z(x,y)\psi(x,y)(u(x)-u(y)) dxdy=
		\int_{\R^N}\int_{\R^N}\psi(x,y)|u(x)-u(y)| dxdy
		$$
		for all $\psi\in C^\infty_c(\R^N\times\R^N)$. This clearly implies property (a) of the statement. 
		The minimality of the sets $\{ u < t \}$ and $\{ u \leq t \}$ for the geometric problem can be obtained 
		by repeating Step 3 of the proof of Theorem~\ref{give me a name}.
		 The uniqueness of the function $u^g_h$ follows by the arguments of Step 4 
		 of the proof of Theorem~\ref{give me a name}, and using Lemma~\ref{lm:geometricityf}
		 in place of Lemma~\ref{lm:geometricity}.
		 The case of $g\in \mathscr{L}^-$ can be treated by using the arguments of Step 5 
		 of the proof of Theorem~\ref{give me a name}. More precisely, let us set $\tilde g:=-g$. 
		 Then, $\tilde g \in \mathscr{L}^+$.
		Let now $(u^{\tilde g}_h, \tilde z)\in \mathscr{L}^+\times  L^\infty \left(\R^N \times \R^N\right)$ satisfy the conclusions of Theorem~\ref{give me a name f} with $g$ replaced by $\tilde g$. Then, if we set
		$u_h^g:=- u^{\tilde g}_h$ and $z (x, y) = \tilde z (y, x)$, then the pair $(u_h^g, z)$ satisfies (a) and (b) and $u_h^g$ is unique in $\mathscr{L}^-$.

Let us now prove the comparison principle. 
Suppose first that $g_1\in  \mathscr{L}^-$ and $g_2\in  \mathscr{L}^+$ and, for $i = 1, 2$, let $(u^{g_i}_h, z_i)$ 
be the pair given by the previous steps corresponding to $g_i$.		 
Arguing as in Step~6 of the proof of Theorem~\ref{give me a name} and using Lemma~\ref{lm:monotonicityf} in place of Lemma~\ref{lm:monotonicity}, we obtain that
\begin{align*}
			0 &\leq h \, \frac{c_s}{2} \int_{\R^N} \int_{\R^N} (z_1(x,y)-z_2(x,y)) \frac{\big(u^{g_1}_h-u^{g_2}_h\big)_+ (x)
			- \big(u^{g_1}_h-u^{g_2}_h\big)_+ (y)}{|x-y|^{N+s}} \, dxdy \\
			&=- \int_{\R^N} \big(u^{g_1}_h-u^{g_2}_h\big)_+ \big(u^{g_1}_h-u^{g_2}_h\big) \, dw \leq 0\,,
		\end{align*}		 
	which implies   $u^{g_1}_h\leq u^{g_2}_h$.	 	 
The other cases follow as in the last part of the proof of Theorem~\ref{give me a name}.
	\end{proof}
	
	In what follows,  $T_h(\overline{B}_R)$  will denote the maximal solution  of the problem
	\[
		\min_{F \in \mathcal{M} (\R^N)} \left\{ \mathcal{P}_s (F) + \frac{1}{h} \int_{F} (|x|-R) \, dx \right\},
	\]
	while $\kappa_s^R$ is the (constant)  $s$-fractional curvature of the ball $B_R$. Finally, we will let  $\phi^h$ denote the solution $u^g_h$ provided by Theorem~\ref{give me a name f}, with $g(x):=|x|$.
	The analogue of Proposition~\ref{prop no name} now reads: 
	\begin{lemma}\label{prop no name f}
	For every  $\delta>0$ there exists $h_0=h_0(s, \delta)>0$ such that for every $h\in (0, h_0)$ and $R \geq \delta$
	it holds
	\beq\label{DR}
		T_h(\overline{B}_R) \supseteq \overline{B}_{R-\kappa_s^{\delta/2}h}.
	\eeq
	Moreover, setting $c_0 = \kappa_s^{\delta/2}$, we have $\phi^h(x) \leq |x|+c_0 h$ whenever $|x|+c_0h \geq \delta$  and  $h\in (0, h_0)$. 
	Finally, for every $\alpha \in (0, 1/(N+1))$ there exists $h_1 = h_1 (\alpha) > 0$ such that
		\[
			\phi^h (0) \leq 2 h^{\alpha}
		\]
		for every $h \leq h_1$. 
	\end{lemma}
	\begin{proof}
	Given $\delta>0$, we apply \cite[Lemma 6.6 (ii)]{CMP15} with $R_0=\delta$ and $\sigma=2$, to infer the existence of $h_0=h_0(s, \delta)>0$ such that for every $h \leq h_0$ and $R \geq \delta$ it holds
		\begin{equation}\label{2}
			T_h(\overline{B}_R) \supseteq \overline{B}_{R-\kappa_s^{R/2}h}.
		\end{equation}
			Recall now that 		 
	 $$
		\kappa_s^R:=C(N,s)R^{-s}, 
		$$
		where $C(N,s)>0$ is a constant depending only on the dimension $N$ and $s$ (see \cite[Lemma 2]{SaVal19}), so that in particular $R\mapsto \kappa_s^R$ is decreasing. Thus, \eqref{DR} follows from \eqref{2}
		provided  $R \geq \delta$.
		By Theorem \ref{give me a name f}, for every $R \geq \delta$ and $h \leq h_0$ we thus get
		$$
		\{\phi^h \leq R \}=T_h(\overline{B}_R) \supseteq \{ |x|  +c_0h \leq R \}=\overline{B}_{R-c_0h}.
		$$
		Therefore, for every $h \leq h_0$, we conclude that $\phi^h(x) \leq |x|+c_0 h $ whenever $|x|+c_0h \geq \delta$. 
		Finally, the proof of last part of the statement follows the same lines of the proof of \eqref{anche questa finalmente va bene},  and we omit it.
	\end{proof}
		
		
	\subsection{The fractional weak flow}\label{sec:weakformulationf}

	In this section we introduce the distributional weak formulation of the fractional mean curvature flow \eqref{oeef}.
	To this aim, given $z \in L^\infty\left(\R^N\times\R^N  \times (0,T)\right)$ we will still denote by $\divs z$ its nonlocal ``spatial'' divergence, as the distribution acting as  
	\begin{equation} \label{divss}
		\langle  \divs z, \vphi \rangle:= - \frac{c_s}{2} \int_0^T\int_{\R^N} \int_{\R^N} z(x,y,t) \frac{\varphi(x,t)-\varphi(y,t)}{|x-y|^{N+s}} \, dxdydt 
	\end{equation}
	for all $\varphi\in C_c^{\infty}\left(\R^N \times (0,T)\right)$. Note that the above definition extends to all $\varphi$ such that $\int_0^T\mathcal{J}_s(\varphi(\cdot, t))\, dt<\infty$. 
	In the following definition, $T^{\textnormal{ext}}$ and $d$ are as in Definition~\ref{Defsol}.
	\begin{definition}\label{Defsolf}
		Let $E^0\subset\R^N$ be 
		a closed set with compact boundary. 
		Let $E$ be a closed set in $\R^N\times [0,+\infty)$ and
		for each $t\geq 0$ denote $E(t):=\{x\in \R^N:\, (x,t)\in E\}$. 
		We say that $E$ is a {\em weak superflow} of the curvature flow \eqref{oeef} with
		initial datum $E^0$ if  conditions  (a), (b), and (c) of Definition~\ref{Defsol} are satisfied, and		
		\begin{itemize}
		\item[(d)] there exists $z\in L^\infty\left(\R^N  \times\R^N\times (0,T^{\textnormal{ext}})\right)$ with 
$\|z\|_{L^\infty\left(\R^N  \times\R^N\times (0,T^{\textnormal{ext}})\right)} \leq 1$  such that the following holds:
		
\begin{itemize}
			
\vspace{.1cm}
	
\item[(d1)]  setting $d(x,t):=\dist (x, E(t))$, we have
\begin{equation}\label{eq:supersolf}
				\partial_t d \geq \divs z \quad \text{ in } \mathcal{D}' ((\R^N\times (0,T^{\textnormal{ext}}))\setminus E);
			\end{equation}
		
\vspace{.1cm}
	
\item[(d2)] for  a.e. $(x,y, t) \in \R^N\times\R^N\times (0,T^{\textnormal{ext}})$ 
\beq \label{corrieref}
z(x,y,t)(d(x,t)-d(y,t)) = |d(x,t)-d(y,t)|;
\eeq
	\item[(d3)] for every $\lambda > 0$	
\beq\label{divpositivepartf}
			(\divs z)^+ 
			\in L^\infty(\{(x,t)\in\R^N\times (0,T^{\textnormal{ext}}):\, d(x,t)\geq \lambda\}).
			\eeq
\end{itemize}			
%
%
%
%
		\end{itemize}
		  {\em Weak subflows} and  {\em weak flows} for \eqref{oeef} are defined accordingly, 
		as at the end of Definition~\ref{Defsol}.
	\end{definition}
		
	\begin{remark}\label{in the sense of measures}
		The bound \eqref{divpositivepartf} implies that $\divs z$ is a measure in 
		$(\R^N\times (0,T^{\textnormal{ext}}))\setminus E$ and  thus \eqref{eq:supersolf} holds also in the sense of measures.
		 We also notice that condition (d2) forces 
		\[
		z(x,y,t)=\mathrm{sign}\big( d(x,t)-d(y,t) \big)\qquad \text{ for  a.e. $(x,y, t) \in \R^N\times\R^N\times (0,T^{\textnormal{ext}})$. }
		\]
		This follows from the fact that for all $t\in (0,T^{\textnormal{ext}})$ the positive level sets of the distance function $d(\cdot, t)$ have vanishing Lebesgue measure and thus $d(t,x)-d(t,y)\neq 0$ for a.e. $(x,y)\in \{d(\cdot, t)>0\}$.		
	\end{remark}
	Clearly, all the observations of Remark~\ref{consequences of (b)} apply also in this case. 
	We can also  define the {\em level set supersolutions (subsolutions)} coresponding to  \eqref{oeef}  as in Definition ~\ref{deflevelset1}, by requiring that a.e. closed (open) sublevel set is a weak superflow (subflow) of \eqref{oeef}, according to Definition~\ref{Defsolf}.
	
	The existence of weak subflows and superflows is achieved also in this case by implementing the minimizing movement scheme of Section~\ref{sec:atw}, specified for the fractional flow. Precisely: Given $E^0\subset \R^N$  with compact boundary,
	we set $E^0_h=E^0$. We then inductively define $E_h^{k+1}$ (for all $k\in \N$) according to the following  procedure: 
	If $E_h^{k}\neq \emptyset$, $\R^N$, then let   $(u_h^{k+1},z_h^{k+1}) 
	\in \Lip(\R^N)\times  L^\infty\left(\R^N\times\R^N\right)$, with 
	$\|z_h^{k+1}\|_{L^\infty\left(\R^N\times\R^N\right)} \leq 1$, satisfy 
	\beq \label{eq:iterkf}
	\begin{cases}
		\displaystyle z_h^{k+1}(x,y)(u_h^{k+1} (x)-u_h^{k+1}(y))   = |u_h^{k+1} (x)-u_h^{k+1}(y)|&\text{for a.e. }(x,y)\in \R^N\times\R^N,\vspace{.3cm} \\
		- h \, \divs z_h^{k+1} + u_h^{k+1} = \dd_{E_h^k} \qquad \text{ in } \mathcal{D}' (\R^N),
	\end{cases}
	\eeq
	and  set  $E_h^{k+1}:=\{u_h^{k+1}\le 0\}$. If either $E_h^{k}=\emptyset$ or $E_h^{k}=\R^N$, then set $E_h^{k+1}:=E_h^{k}$. 
	We recall that the definition of nonlocal divergence $\divs$ is given in \eqref{divr0f}. As in Section~\ref{sec:atw} we can define $T^{\, \textnormal{all}}_h$, $T^{\, \textnormal{ext}}_h$, $T^*_h$ and observe that $E_h^{k}$ is closed for all $k$  and if $E^0$  ($(E^0)^c$) is bounded, then $E_h^{k}$ ($(E_h^{k})^c$) is bounded for all $k$. Moreover, \eqref{eq:ineqd} holds true also in this case. Finally, we may define $E_h$, 
	$E_h(t)$, $d_h(\cdot,t)$, $u_h(\cdot,t)$ as in \eqref{discretevol}, and $z_h(\cdot, \cdot, t):= z_h^{[t/h]}$,  $z_h: \R^N \times\R^N\times [h, T^*_h) \to \R$.
	In view of Theorem~\ref{give me a name f}, we can prove the discrete comparison principle of Remark~\ref{rm:dcp} also in the present case. In turn, all the considerations of Section~\ref{sec:balls} are still valid, thanks also to Remark~\ref{prop no name f}.
	
	Hence, we can extract a subsequence $h_\ell$ such that 
	\beq \label{E and A f}
	\overline E_{h_l} \stackrel{\mathscr K}{\longrightarrow} E\qquad\text{and}\qquad 
	{(\mathring{E}_{h_l})}^c\stackrel{\mathscr K}{\longrightarrow} A^c,
	\eeq
	and the conclusions of Proposition~\ref{prop:E} (and Remark~\ref{per passare al limite nella formulazione debole}) hold. Also  in this case, as in  Remark~\ref{scemo chi legge}, we can define 
	$T^{\textnormal{ext}}$, ${T}^{\textnormal{all}}$, and $T^*= \min \{ {T}^{\textnormal{ext}}$, ${T}^{\textnormal{all}} \}$.
	
	\begin{theorem}\label{themthmf}
		Let $E$ and $A$ be the sets given in \eqref{E and A f}.
		Then, $E$ (resp. $A$) is a weak superflow (resp. subflow)
		with initial datum $E^0$ in the sense of Definition~\ref{Defsolf}.
	\end{theorem}
	\begin{proof}
		Also in this case the proof is very similar to that of Theorem~\ref{themthm} and we only indicate the main changes. We only check condition (d) in $(0,T^*)$.  Possibly extracting a new subsequence not relabeled and setting $z_{h_l}(\cdot,\cdot,t):=0$ for $t \in (0,T^*) \setminus [h_l,T_{h_l}^*)$, we can assume that the sequence $z_{h_l}$ converges weakly-$*$ in $L^\infty(\R^N \times \R^N\times (0,T^*)) $ to some function $z$, with  $\Vert z \Vert_{L^\infty(\R^N \times \R^N\times (0,T^*))} \leq 1.$
		
		Arguing in the same way as for \eqref{oscineq}, we get
		\begin{align*}
			&0\leq \int_0^{T^*}  \left( \int_{\R^N}  \eta(x,t) \frac{d_{h_l}(x,t+h_l)-d_{h_l}(x,t)}{h_l} \, dx
			+  \frac{c_s}{2} \int_{\R^N} \int_{\R^N} \frac{(\eta(x,t)-\eta(y,t))}{|x-y|^{N+s}} z_{h_l}(x,y,t) \, dx \, dy   \right)  \!dt\\
			& =\int_0^{T^*}  \left( \int_{\R^N} -\frac{\eta(x,t)-\eta(x,t-h_l)}{h_l} d_{h_l}(x,t) \, dx
			+  \frac{c_s}{2} \int_{\R^N} \int_{\R^N} \frac{(\eta(x,t)-\eta(y,t))}{|x-y|^{N+s}} z_{h_l}(x,y,t) \, dx \, dy   \right)   \, dt.
		\end{align*}
		Letting $l\to \infty$ we obtain
		\begin{align*}
			\int_0^{T^*}  \left( \int_{\R^N} -\partial_t \eta(x,t) d(x,t) \, dx 
			+ \frac{c_s}{2} \int_{\R^N} \int_{\R^N} \frac{(\eta(x,t)-\eta(y,t))}{|x-y|^{N+s}} z(x,y,t) \, dx \, dy   \right)   \, dt \geq 0,
		\end{align*}
		for every nonnegative  $\eta \in C^{\infty}_c ((\R^N \times (0,T^*))\setminus E)$.  Thus, 
		$\partial_t d\geq \divs z$  in $ \mathcal{D}'((\R^N \times (0,T^*)) \setminus E).$
		
		Now one can prove \eqref{bounddiv}, with $\mathfrak{Div}^r z$ replaced by $\divs z$,  \eqref{speed of conv} and \eqref{elleuno} for a.e. $t \in (0, T^*)$  arguing in the same way as in Theorem~\ref{themthm} (using Lemma~\ref{prop no name f} in place of Proposition~~\ref{prop no name}).
		
		Finally, we show \eqref{corrieref} with $\td$ in place of $d$, where $\td$ is the function given by  Proposition \ref{prop:E}.
		To this aim, consider  $\eta\in C_c^\infty(\R^N\times\R^N\times  (0,T^*))$. 
		From \eqref{eq:iterkf} it follows that 
		\begin{align}
			& \int_{0}^{T^*}\int_{\R^N}\int_{\R^N} \eta (x,y, t)   z_{h_l} (x,y,t)  (u_{h_l} (x, t)-u_{h_l}  (y, t))\, dxdydt
			\nonumber  \\
			&= \int_{0}^{T^*}\int_{\R^N}\int_{\R^N} \eta (x,y, t)  |u_{h_l} (x, t)-u_{h_l}  (y, t)|\, dxdydt, \label{we will pass to the limit f}
		\end{align}
		for every $l \in \N$.
		
		Let us first consider the right hand side of the expression above.
		Recalling \eqref{elleuno}, by Dominated 
		Convergence  Theorem  we have  
		\beq
			\int_{0}^{T^*}\!\!\!\!\int_{\R^N}\!\!\int_{\R^N} \eta (x,y, t)  |u_{h_l} (x, t)-u_{h_l}  (y, t)|\, dxdydt \to 
			\int_{0}^{T^*}\!\!\!\!\int_{\R^N}\!\!\int_{\R^N} \eta (x,y, t)  |\td (x, t)-\td  (y, t)|\, dxdydt \label{RHSf}
		\eeq
		as $l \to + \infty$.
		Concerning the left hand side of \eqref{we will pass to the limit f}, we first observe that, 
		since $z_{h_l}$ converges to $z$ weakly-$*$ in $L^\infty(\R^N \times \R^N\times (0,T^*))$,
		we have  
		\begin{align}
			&\lim_{l \to \infty} \int_{0}^{T^*}\int_{\R^N}\int_{\R^N} \eta (x,y, t)   z_{h_l} (x,y,t)  (\td (x, t)-\td  (y, t))\, dxdydt \nonumber \\
			&= \int_{0}^{T^*}\int_{\R^N}\int_{\R^N} \eta (x,y, t)   z (x,y,t)  (\td (x, t)-\td  (y, t))\, dxdydt. \label{LHS1f}
		\end{align}
		Note also that
		\begin{align*}
			& \left| \int_{0}^{T^*}\int_{\R^N}\int_{\R^N} \eta (x,y, t)   z_{h_l} (x,y,t)  (u_{h_l} (x, t)-u_{h_l}  (y, t))\, dxdydt \right. \\
			&- \left. \int_{0}^{T^*}\int_{\R^N}\int_{\R^N} \eta (x,y, t)   z_{h_l} (x,y,t)  (\td (x, t)-\td  (y, t))\, dxdydt \right| \\
			&\leq  \int_{0}^{T^*}\int_{\R^N}\int_{\R^N} | \eta (x,y, t)|  
			| u_{h_l} (x, t) - \td (x, t) |  \, dxdydt  \\
			&+ \int_{0}^{T^*}\int_{\R^N}\int_{\R^N} | \eta (x,y, t)| 
			| u_{h_l} (y, t) - \td (y, t) |  \, dxdydt.
		\end{align*}
		where the last two integrals vanish as $l \to \infty$ again by  \eqref{elleuno} and  Dominated Convergence Theorem.  Taking into account \eqref{RHSf} and \eqref{LHS1f}, we obtain
		\[ \int_{0}^{T^*}\!\!\!\!\int_{\R^N}\!\!\int_{\R^N} \eta (x,y, t)   z (x,y,t)  (\td (x, t)-\td  (y, t))\, dxdydt
		= \int_{0}^{T^*}\!\!\!\!\int_{\R^N}\!\!\int_{\R^N} \eta (x,y, t)  |\td (x, t)-\td  (y, t)|\, dxdydt.
		\]
		By the arbitrariness of $\eta$,  \eqref{corrieref} follows for $\td$ and thus also for  $d=\td^+$. This concludes the proof of condition (d)
		and shows that $E$ is a weak superflow for \eqref{oeef}. The fact that $A$ is a weak subflow can be proved 
		in a similar way.
	\end{proof}
	We are now ready to extend  the comparison principle (Theorem~\ref{thm: wfcomparison}, with exactly the same statement) to weak superflows  and  subflows  in the sense of Definition~\ref{Defsolf}.
	\begin{proof}[Proof of the comparison principle] The proof is essentially the same as  that of Theorem~\ref{thm: wfcomparison}. In the bounded case, we can proceed exactly in the same way to define $\td_E$ and $\td_{F^c}$ and obtain all the properties \eqref{w=0}--\eqref{dimuniqflow}  (with of course $\divr z_E$, $\divr z_{F^c}$ replaced by  $\divs z_E$, $\divs z_{F^c}$). Let now $\Psi$ and $\Phi$ be defined as in the proof of  Theorem~\ref{thm: wfcomparison}.  	Then, arguing similarly to the case of \eqref{compafinal},  we get
		\beq\label{compafinalf}
		\begin{split} 
		& \int_S \Phi(x,s) \, dx 
		= \int_{\R^N} \Phi(x,s) \, dx
		= \int_{\R^N\times (0,s)} d (\partial_t \Phi )\\
			& \leq -  \int_{\R^N \times (0, s)} \Psi'(\Delta-(\td_E+\td_{F^c})) 
			\, d(\partial_t^d \td_E+\partial_t^d \td_{F^c}) \\
			& \leq - \int_{\R^N\times (0,s)} \Psi'(\Delta-(\td_E+\td_{F^c})) \, d(\divs (z_E+z_{F^c}))  \\
			& =- \frac{c_s}{2}\int_0^s \int_{\R^N} \int_{\R^N}(-z_E-z_{F^c})(x,y,t)\frac{\Psi'((\Delta-\td_E)-\td_{F^c})(x,t)-\Psi'((\Delta-\td_E)-\td_{F^c})(y,t)}{|x-y|^{N+s}} dxdydt\,.
		\end{split}
		\eeq
		Recall now that 
		$$
		-z_E(x,y, t)\big((\Delta -d_E (x, t))-(\Delta-d_E(y, t))\big)=|(\Delta -d_E (x, t))-(\Delta-d_E(y, t))|
		$$
		and 
		$$
		z_{F^c}(x,y, t)\big(d_{F^c} (x, t)-d_{F^c}(y, t)\big)=|d_{F^c} (x, t)-d_{F^c}(y, t)|
		$$
		for a.e. $(x,y, t) \in \big(\R^N\times\R^N\times (0,T^*)\big)$. Thus, taking into account also \eqref{noticing}, we conclude by Lemma~\ref{lm:monotonicityf} that 
		$$
		(-z_E-z_{F^c})(x,y,t)\Big(\Psi'((\Delta-\td_E)-\td_{F^c})(x,t)-\Psi'((\Delta-\td_E)-\td_{F^c})(y,t)\Big)\geq 0
		$$
		for a.e. $(x,y, t) \in \R^N\times\R^N\times (0,T^*)$. Thus, \eqref{compafinalf} implies
		$\int_S \Phi(w,s) \, dw\leq 0$  for any $s\in(0, t^*)$. We may now continue exactly as in Theorem~\ref{thm: wfcomparison}.
	\end{proof}
	The  above comparison theorem yields the well posedness of the fractional weak level set flow, defined analogously to  Definition~\ref{deflevelset1}. In other words, Theorems~\ref{th:lscomp}~and~\ref{th:phiregularlevelset} extend also to the fractional weak level set flow.
	\begin{remark}[Comparison with viscosity solutions]
		In \cite{Imbert_2009} a  unique level set flow for the fractional mean curvature motion is constructed by solving  the problem 
		\beq\label{imbo}
		\begin{cases} 
			\partial_t u(x,t)= |Du(x,t)|\kappa_s\big(x,\{y:  u(y,t) \leq u(x,t) \}\big) \\
			u(x,\cdot)=u_0\,,
		\end{cases}
		\eeq
		in the viscosity sense,  where $\kappa_s(x, E)$
		 is given by \eqref{scurvature}.
		On the other hand, in  \cite{CMP15} it is shown that the above viscosity flow is the unique limit of the level set minimising movement scheme recalled above. Therefore, the distributional level set flow constructed in this section must coincide with the unique viscosity solution of \eqref{imbo}. 

		A similar remark applies to the regularised Minkowski flow (see Section~\ref{sec:variant}), for which a viscosity theory is available (see \cite{CMP15}), ensuring consistency between the distributional and viscosity level-set formulations.
		
	\end{remark}

\section*{Acknowledgments}
 The authors would like to thank Lucia De Luca for useful discussions, and the anonimous referees for their comments, which helped improving the first version of the manuscript. 
F.~C. and  M.~M.  have been partially supported by PRIN 2022 Project “Geometric Evolution Problems and Shape Optimization (GEPSO)”, PNRR Italia Domani, financed by European Union via the Program NextGenerationEU, CUP D53D23005820006. F.~C. was partially supported also by the EPSRC under the Grant EP/P007287/1 “Symmetry of Minimisers in Calculus of Variations”. D.~R. is supported by the Cluster of Excellence EXC 2044-390685587, Mathematics M\"unster: Dynamics-Geometry-Structure funded by the Deutsche
Forschungsgemeinschaft (DFG, German Research Foundation). The three authors are also members of the Gruppo Nazionale per l’Analisi
Matematica, la Probabilit\`a e le loro Applicazioni (GNAMPA), which is part of the Istituto
Nazionale di Alta Matematica (INdAM).

\section*{Data Availability Statement} 
Data sharing is not applicable to this article as no new data were created or analyzed in this study.

\bibliographystyle{plain}
\bibliography{references}

\end{document}